\documentclass[11pt, leqno]{amsart}

\usepackage[utf8]{inputenc}
\usepackage{graphicx}
\usepackage[T1]{fontenc}
\usepackage[all]{xy}
\usepackage{t1enc} 
\usepackage{amsfonts,amsmath,amssymb,amsthm,color,amsbsy}
\usepackage{color}
\newtheorem{thm}{Theorem}[section]
\newtheorem{lem}[thm]{Lemma}

\newtheorem{asm}[thm]{Assumption}
\theoremstyle{definition}
\newtheorem{df}[thm]{Definition}
\newtheorem{rem}[thm]{Remark}

\definecolor{labelkey}{rgb}{0,0,1}
\usepackage{dsfont}
\def\pr{\right )}
\def\le{\left (}

\def\R{\mathbb{R}}

\def\mm{\kern +2pt \raisebox{+0.5 pt}{$\shortmid$}\kern -2pt\hbox{$\multimap$}\kern +2pt}
\def\H1{H^1\le\Omega,\mathbb{R}^N\pr}
\def\h2{H^2\le\Omega,\mathbb{R}^N\pr}
\def\h10{H^1_0\le\Omega,\mathbb{R}^N\pr}

\begin{document}
\title{Stability of phase portrait for a gradient ODE with memory}

\author{Piotr Kalita, Piotr Zgliczyński}
\address{Faculty of Mathematics and Computer Science, Jagiellonian University,
	ul. Łojasiewicza 6, 30-348 Kraków, Poland}
\email{piotr.kalita@ii.uj.edu.pl, umzglicz@cyf-kr.edu.pl}

\begin{abstract}
We consider the problem governed by the gradient ODE $x'=\nabla F(x)$ in $\mathbb{R}^d$ on which we assume that it has a finite number of hyperbolic equilibria whose stable and unstable manifolds intersect transversally. This problem is perturbed by the memory term 
$$x'(t)=\nabla F(x(t))+\varepsilon\int_{-\infty}^t M(t-s)x(s)\, ds$$ where $\varepsilon>0$ is a small constant. The key result is that the structure of connections between the equilibria of the unperturbed problem is exactly preserved for a small $\varepsilon>0$.
\end{abstract}

%%
%%  LaTeX will not make the title for the paper unless told to do so.
%%  This is done by uncommenting the following.
%%

 \maketitle
 
% \today
 
\tableofcontents

\section{Introduction.}
\noindent This paper deals with the gradient ordinary differential equation in $\mathbb{R}^d$
\begin{equation}\label{unperturbed}
x'(t) = \nabla F(x(t))\ \ \textrm{for}\ \ F\in C^3(\mathbb{R}^d),
\end{equation}
and its perturbation by the linear memory term though which the derivative of the unknown solution depends not only on the instantaneous value of the solution but also on its past values
\begin{equation}\label{perturbed}
x'(t) = \nabla F(x(t)) + \varepsilon\int_{-\infty}^tM(t-s)x(s)\, ds.
\end{equation}
Because of the presence of the distributed delay term, if we study the flow governed by \eqref{perturbed}, we need to consider it in an infinite dimensional space, containing functions defined in the time interval from minus infinity to the current time instance.

We make the following assumptions on the functions $F:\R^d\to \R$ and $M:[0,\infty)\to \R^{d\times d}$:
\begin{itemize}
	\item[(1)] There exists $R>0$ and $C_F>0$ such that if $|x|\geq R$ then $-\nabla F(x)\cdot x \geq C_F|x|^2$.
	\item[(2)] The unperturbed equation \eqref{unperturbed} has finite number of equilibria, all of them being hyperbolic, and their stable and unstable manifolds intersect transversally.
	\item[(3)] The matrix $\int_0^\infty M(s)\,  ds$ is symmetric and $\|M(s)\| \leq E\lambda_{min}(A(s))$ for a constant $E>0 $ (cf. Assumption \ref{asm:M}) where  $A\in C^1([0,\infty);\R^{d\times d})$ with $A(s)$ being symmetric and positive definite matrix for every $s\geq 0$ such that
	\begin{itemize}
		\item[(A)] For almost every $s>0$ and for every $x\in \R^d$ we have
		$$
		\left(\frac{dA(s)}{ds}x,x\right)\leq -C_A(A(s)x,x),
		$$
		with a constant $C_A>0$ cf. Assumption \ref{asm:a1}.
		\item[(B)] For every $s\geq 0$ we have
		$$
		\frac{\lambda_{max}(A(s))}{\lambda_{min}(A(s))}\leq D,
		$$
		with a constant $D>0$, cf. Assumption \ref{asm:a2}.
	\end{itemize}
\end{itemize}

 The function $A$, which by (3)(A) must decay exponentially to zero as $t\to \infty$, defines the phase space for the memory term, it is the space
 $$
 L^2_A(\R^+)^d = \left\{ \eta:[0,\infty)\to \R^d\, :\ \int_0^\infty (A(s)\eta(s),\eta(s))\, ds < \infty \right\}.
 $$
 Now, the equation \eqref{unperturbed} defines the gradient dynamical system $S^0(t):\R^d\to \R^d$ for $t\geq 0$. This dynamical system has a  global attractor which
 consists of the finite number of equilibria and their connections. It structure is represented as a graph of partial order, the vertexes of this graph correspond to the equilibria of the system. An edge from $e_i$ to $e_j$ exists in this graph it there exists a bounded  solution of \eqref{unperturbed} which converges to $e_i$ as time tends to minus infinity and to $e_j$ as time tends to plus infinity.

 The problem governed by \eqref{perturbed} defines a dynamical system for $\varepsilon>0$ denoted as $S^\varepsilon(t):L^2_A(\R^+)^d\times \R^d \to L^2_A(\R^+)^d\times \R^d$, where the space $\R^d$ contains the current state of the system, and  $L^2_A(\R^+)^d$ its past. This system for $\varepsilon > 0$ is infinite dimensional.
The main result of the paper is the following theorem
\begin{thm}
	Assume (1)--(3) above. There exists $\varepsilon_0>0$ such that for every $\varepsilon \in [0,\varepsilon_0]$  the dynamical system governed by \eqref{perturbed} has a global attractor consisting of a finite number of equilibria and their connections. The graph that represents this system coincides with the graph for the unperturbed finite dimensional system $\{S^0(t)\}_{t\geq 0}$.
	\end{thm}

	The question is motivated by the results of \cite{borto}. There, the authors consider the infinite dimensional autonomous gradient dynamical system and they prove that upon small non-autonomous perturbation the structure of its attractor is preserved, that is, the phase portrait of the non-autonomous dynamics coincides with the autonomous one. Thus, the authors in \cite{borto} are able to fully characterize the non-autonomous dynamics for the problem which is small perturbation of the autonomous one. 
	%(see also \cite{LipBorto} for a similar result where the small perturbation is autonomous, but not $C^1$ - only Lipschitz).
	 Our result is of similar nature as \cite{borto}, but our main contribution stands in the fact that the unperturbed system is finite-dimensional and the perturbed one - infinite dimensional. This infinite dimensionality manifests itself in the presence of the memory term in the perturbed system.
	
	The proof that the structure of connections is exactly preserved upon perturbation consists of three ingredients:
	\begin{itemize}
		\item[(A)] the equilibria of the perturbed problem exist in the vicinity of the equilibria of the original one, and that these are all equilibria, cf. Theorem \ref{thm:equilibria},
		\item[(B)] no new connections arise when $\varepsilon>0$, i.e. the connections structure behaves upper-semi\-con\-ti\-nuously, cf. Theorem \ref{thm:upper_con},
		\item[(C)] the existing connections are preserved upon perturbation, i.e. the connections structure behaves lower-semicontinuously, cf. Theorem \ref{thm:lower}.
	\end{itemize}
	Fundamental ingredient in the proofs of these items is the fact, obtained in Sections \ref{sec:setup} and \ref{sec:LapCC}, that certain dynamical properties of the unperturbed problem can be continued for $\varepsilon>0$. In particular it is possible to construct the common Lyapunov function for $\varepsilon \in [0,\varepsilon_0]$. Moreover,  we  construct isolating blocks: i.e. the sets which isolate the equilibria of \eqref{unperturbed}, which after taking the Cartesian product with a certain ball in the memory space $L^2_A(R^+)^d$  also isolate the equilibria of \eqref{perturbed} with entry and exit behavior on the boundary being uniform with respect to $\varepsilon\in [0,\varepsilon_0]$. We prove that the new, infinite dimensional, variable $\eta \in L^2_A(R^+)^d$ can be bundled together with the stable variables belonging to the finite dimensional state $x\in \R^d$. Finally, we prove that the cone condition holds in these isolating sets with the same system of coordinates and the same quadratic form in the range of small $\varepsilon\in [0,\varepsilon_0]$.
	This opens the possibility of using the Hadamard's graph transform procedure to  construct the local stable and unstable manifolds of the equilibria as the Lipschitz graphs over the same systems of coordinates in the considered range of $\varepsilon$. The above assertion (A) follows from the construction of common isolating blocks with the cone condition, and (B) follows from the compactness argument (similar as in \cite{borto}), these results are contained in Section \ref{equi}.
	To get the most involved result (C), we need to prove that the local stable and unstable manifolds are actually $C^1$ close to each other in dependence on $\varepsilon$. We prove this by the differentiation of the graph transform. Moreover, we transport the smallness of $C^1$ distance between the local unstable manifolds along the flow in order to prove that the transversality of the intersection for $\varepsilon=0$ implies that this intersection is preserved upon the perturbation. This argument needs careful handling of the memory variable which appears in the system for $\varepsilon>0$. The result is contained in Section \ref{intersection}.

The fact that the norm of the memory term is weighted by the expression that decays exponentially to zero is a fundamental fact which allows us to treat the memory variable as the stable variable in the neighborhood of the equilibrium. The key result here is the dissipative estimate \eqref{eq:1} in Lemma \ref{lem:eta} on the time evolution of the norm of memory variable which is obtained in Section \ref{sec:memo_prop}. This estimate is derived using the concept from the seminal paper of Dafermos \cite{Dafermos}, who proved that in the linear problem of viscoelasticity the memory term is dissipative and has damping effect on the solution, which decays to zero due to this term's presence. Discoveries of Dafermos were later used in the context of global attractors for the nonlinear problem of viscoelasticity by Conti and Pata \cite{Pata}, who explored the dissipative nature of the memory term to obtain the existence of the global attractor. Dissipativity of the memory term in the context of global attractors has also been explored for the first order, reaction-diffusion type, problems in \cite{PataSquassina, Cara1, Cara2}, where the authors showed that it is possible to proceed directly using the variables of the system, without the Dafermos transformation. All these results, however, are of global nature. The novel contribution of this paper, is the exploitation of the dissipative nature of the memory term in the local argument realized in the neighbourhood of the equilibria and its application to recover the full intrinsic structure of the global attractor.

Another novelty of this paper is the development of the geometric methods to study the dynamics of the problem with memory/delay. Locally, in the neighbourhood of the equilibrium we apply the Hadamard procedure to obtain the existence of local stable and unstable manifolds. Then, we  show that local unstable manifolds for the perturbed problems are $C^1$-close to those of the original problem. We continue this proximity along the solutions, in order to prove that these manifolds are $C^1$ close  near the point of intersection, which lies in the neighborhood of the target point of the heteroclinic connection. Thus, we develop a geometric approach for problems with memory, which is useful for studying the dynamics of the problem. We stress that while we apply these methods to a particular class of problems, we expect that the proposed techniques can be extended to a wider class of problems with memory and/or delay. In particular, while we prove that the structure of connections in a gradient system is preserved under perturbation, we expect that our techniques will allow us to prove more general results, demonstrating that the dynamics of the problem remains unchanged under perturbation by the memory term. This is possible because, in our framework, the additional memory variable is, in some sense, damped. Moreover, as an interesting open question we ask, we ask, whether it is possible to study, in a similar way, the dynamics for the case where  the kernel in the memory term does not decay exponentially, but, for example, polynomially.

   One of the ingredients of our proof is the result on the existence of local stable and unstable manifolds. While we use the graph transform method, the alternative approach - the Lyapunov--Perron method - has been used to obtain the similar types of results for the problems with finite delay  in \cite{Krisztin}. We note that the results on the attractor structure for problems with delay or memory are not frequent:  we refer for example to \cite{Walther} for an overview of results and list of relevant references or to \cite{KV} for the study of fine structure of the unstable manifold of a periodic solution for finitely delayed ODE.

Due to the fact that some of the arguments of the article are  technical, they are moved to Appendices. Appendix 1 contains the results on asymptotic compactness of the dynamical system with memory. Appendix 2 is devoted to the Hadamard's graph transform procedure and the results on the existence of local stable and unstable manifolds. In Appendix 3 the results on $C^0$ dependence of those manifolds on parameter $\varepsilon$ are proved. Appendix 4 presents the approach to graph transform method using the Banach fixed point theorem, this argument allows us to differentiate the graph transform map and, using the fiber contraction argument,  obtain the $C^1$ smoothness of stable and unstable manifolds. In Appendix 5 the abstract cone conditions from Appendix 4 are demonstrated to hold  for the considered problem.   Finally, Appendix 6 is devoted to the proof that the derivatives of local stable and unstable manifolds depend continuously on the parameter $\varepsilon$.

\section{The weighted history space and its norm.}
Let $A: [0,\infty) \to \R^{d\times d}$ be a time dependent matrix function. This function will be used in the definition of the norm in the history space. Assumptions \ref{asm:a1}, \ref{asm:a2} on the function $A$ and Assumption \ref{asm:M} on related matrix function $M$ will be standing assumptions throughout the whole article.
\begin{asm}\label{asm1}\label{asm:a1}
	Assume that $A(s)$ is a symmetric and positive definite matrix for $s\geq 0$, $[0,\infty) \ni s\mapsto A(s)$ belongs to $C^1([0,\infty);\R^{d\times d})$ and that for almost every $s>0$ and every $u\in \R^d$
	\begin{equation}\label{eq:diffA}
	\left(\frac{dA(s)}{ds}u,u\right) \leq -C_A(A(s)u,u).
	\end{equation}
\end{asm}
\begin{lem}
	Under Assumption \ref{asm1} we have
	$$
	\int_0^\infty \|A(s)\|\, ds < \infty
	$$
\end{lem}
\begin{proof}
	We have
	$$
	 e^{C_As} \frac{d}{ds} (A(s)u,u) + C_Ae^{C_As}(A(s)u,u) \leq 0,
	$$
	for every $u\in \R^d$. Hence
	$$
	\frac{d}{ds} (e^{C_As}(A(s)u,u)) \leq 0,
	$$
	and
	$$
	e^{C_As}(A(s)u,u) \leq (A(0)u,u), \quad s \geq 0
	$$
	Finally
	$$
	(A(s)u,u) \leq e^{-C_As} (A(0)u,u),
	$$
	for every $s\geq 0$ and $u\in \R^d$.
	As $A(s)$ is symmetric and positively definite then for every $s$ we can find a vector $u(s)$ with norm one such that
	$$
	\|A(s)\| = (A(s)u(s),u(s)) \leq e^{-C_As}(A(0)u(s),u(s)) \leq \|A(0)\| e^{-C_As},
	$$
	and the assertion follows.
\end{proof}

%\begin{asm} \label{asm:a2}
%		Assume that for some constant $D>0$ and every $x\in \mathbb{R}^d$, $s\in \mathbb{R}^+$
%	\begin{equation}\label{eq:A_property}
%		\|A(s)\||x|^2 \leq D^2(A(s)x,x).
%	\end{equation}
%In other words
%$$
%\frac{\lambda_{max}(A(s))}{\lambda_{min}(A(s))} \leq D^2\ \ \textrm{for every}\ \ s\geq 0.
%$$
%\end{asm}
%\begin{lem}\label{lem:a_prop}
%
%	Under assumptions \ref{asm1} and \ref{asm:a2} for every $\eta \in L^2_A(\mathbb{R}^+)^d$
%	$$
%	\left|\int_0^\infty A(s)\eta(s)\, ds \right| \leq \left(D \sqrt{\int_0^\infty \|A(s)\|\, ds}\right)\|\eta\| := D'\|\eta\|.
%	$$
%\end{lem}
%\begin{proof}
%	\begin{align*}
%		& \left|\int_0^\infty A(s)\eta(s)\, ds \right|^2 \leq \left(\int_0^\infty \|A(s)\||\eta(s)|\, ds \right)^2.
%	\end{align*}
%	By the H\"older inequality
%	\begin{align*}
%		&
%	\left|\int_0^\infty A(s)\eta(s)\, ds \right|^2 \leq \left(\int_0^\infty \sqrt{\|A(s)\|} \sqrt{\|A(s)\|} |\eta(s)|\, ds \right)^2\\
%	& \ \ \leq \int_0^\infty \|A(s)\|\, ds \int_0^\infty \|A(s)\||\eta(s)|^2\, ds \leq D^2\int_0^\infty \|A(s)\|\, ds \int_0^\infty (A(s)\eta(s),\eta(s))\, ds.
%	\end{align*}
%	and the proof is complete.
%\end{proof}

We define the space $L^2_A(\R^+)^d$ with the norm  $\|\eta\|^2 = \int_0^\infty (A(s)\eta(s),\eta(s))\, ds.$
\begin{asm} \label{asm:a2}
	Assume that for some constant $\overline{D_A}>0$ and every $x\in \mathbb{R}^d$, $s\in \mathbb{R}^+$
	\begin{equation}\label{eq:A_property}
		\|A(s)\|\ |x|^2 \leq \overline{D_A}^2(A(s)x,x).
	\end{equation}
	In other words
	$$
	\frac{\lambda_{max}(A(s))}{\lambda_{min}(A(s))} \leq \overline{D_A}^2\ \ \textrm{for every}\ \ s\geq 0.
	$$
\end{asm}
\begin{lem}\label{lem:a_prop}
	
	Under Assumptions \ref{asm1} and \ref{asm:a2} for every $\eta \in L^2_A(\mathbb{R}^+)^d$
	$$
	\left|\int_0^\infty A(s)\eta(s)\, ds \right| \leq \left(\overline{D_A}\sqrt{ \int_0^\infty \|A(s)\|\, ds}\right)\|\eta\| := D_A\|\eta\|.
	$$
\end{lem}
\begin{proof} We begin with the estimate
	\begin{align*}
		& \left|\int_0^\infty A(s)\eta(s)\, ds \right|^2 \leq \left(\int_0^\infty \|A(s)\|\, |\eta(s)|\, ds \right)^2.
	\end{align*}
	By the H\"older inequality
	\begin{align*}
		&
		\left|\int_0^\infty A(s)\eta(s)\, ds \right|^2 \leq \left(\int_0^\infty \sqrt{\|A(s)\|} \sqrt{\|A(s)\|}\,  |\eta(s)|\, ds \right)^2\\
		& \ \ \leq \int_0^\infty \|A(s)\|\, ds \int_0^\infty \|A(s)\|\, |\eta(s)|^2\, ds \leq \overline{D_A}^2\int_0^\infty \|A(s)\|\, ds \int_0^\infty (A(s)\eta(s),\eta(s))\, ds.
	\end{align*}
	and the proof is complete.
\end{proof}

Now consider the function $M:[0,\infty) \to \R^{d\times d}$. We make the following assumption
\begin{asm}\label{asm:M}
	Assume that for every $s\geq 0$
	\begin{equation}\label{eq:boundA}
	{\|M(s)\|} \leq \overline{D_M}^2 \lambda_{min}(A(s)),
	\end{equation}
	with a constant $\overline{D_M} > 0$ and
	$$
	\int_0^\infty M(s)\, ds\ \ \mathrm{is\ symmetric}.
	$$
\end{asm}
\begin{rem} The fact that $\int_0^\infty M(s)\, ds$ is finite follows from \eqref{eq:boundA} and the assumptions on $A$. The symmetry of this integral is needed only  to construct the Lyapunov function in Lemma \ref{lem:LapFunc1}, and although it appears technical and having no natural explanation, it is unclear to us how to avoid this assumption in the argument. \end{rem}

The next result holds analogously to Lemma \ref{lem:a_prop}
\begin{lem}\label{lem:M}
	Under Assumptions \ref{asm:a1} and \ref{asm:M} we have
	$$
	\int_0^\infty \|M(s)\|\, ds < \infty
	$$
	and
	$$
	\left| \int_0^\infty M(s)\eta(s)\, ds \right|\leq D_M \|\eta\|\ \ \ \textrm{for every}\ \ \ \eta\in L^2_A(\R^+)^d,
	$$
	for every $\eta\in L^2_A(\R^+)^d$, where $D_M = \overline{D_M}\sqrt{\int_0^\infty \|M(s)\|\, ds}$.
\end{lem}
An example of $A(s)$ which satisfies the above assumptions is $A(s)=e^{-\kappa s} I$. Then $C_A=\kappa$, $\overline{D_A} = 1$, and we need, in addition to the symmetry of the integral of $M$ that
$$
\|M(s)\|\leq \overline{D_M}^2 e^{-\kappa s}.
$$

\section{Problem setup and its basic properties}\label{sec:setup}

\subsection{Unperturbed and perturbed problems. Dafermos transform and memory variable.} We consider the following ODE
\begin{equation}
	x'(t) = f(x(t)) \quad \textrm{where}\quad f\in C^2(\R^d;\R^d).  \label{eq:ODE}
\end{equation}
We assume that the ODE has a gradient form, i.e.
\begin{equation}\label{gradient}
f(x) = \nabla F(x) \ \ \textrm{where} \  \ F\in C^3(\R^d).
\end{equation}
Moreover we assume that 
\begin{equation}\label{diss2}
	\text{there exist constants}\  R>0, C_F>0\ \text{such that if}\ |x|\geq R\ \text{then}\ -\nabla F(x)\cdot x\geq C_F|x|^2. 
\end{equation}
We need a following simple property
\begin{lem}
	Assume \eqref{diss2}. Then there exist constants $\gamma>0$ and $\delta\in \R$ such that
	\begin{equation}\label{diss}
		F(x) \leq -\gamma |x|^2 + \delta.
	\end{equation}
\end{lem}
\begin{proof}
	If $|x|\leq R$ then for every constant $\gamma > 0$ we have 
	$$
	F(x) \leq \max_{|x|\leq R}F(x) \leq -\gamma |x|^2 + {\gamma} R^2 + \max_{|x|\leq R}F(x) = -\gamma |x|^2 + \delta_1(\gamma).$$
	If $|x|>R$ then $x=(1+c)x_0$ for some $|x_0|=\frac{R}{2}$ with $c>1$. Then
	\begin{align*}
	& F(x) = F(x_0) + \int_0^1 \nabla F(x_0(1+\theta c)) c x_0 \, d\theta = F(x_0) + \int_0^1 \nabla F(x_0(1+\theta c)) \frac{1+\theta c}{1+\theta c}c x_0 \, d\theta \\
	& \ \ \leq  F(x_0) - C_F \int_0^1 |x_0|^2c(1+\theta c) \, d\theta = F(x_0) - C_F \int_0^1 |x_0|^2c(1+\theta c) \, d\theta \\
	& \ \  = F(x_0) - C_F |x_0|^2c \left(1+\frac{c}{2}\right) = F(x_0) - C_F |x|^2\frac{c \left(2+c\right)}{2(1+c)^2}\\
	& \leq  \max_{|x|\leq R}F(x) - \frac{3}{8}C_F |x|^2,
	\end{align*}
	and the proof is complete. 
\end{proof}
 We perturb the above ODE with the additive linear distributed delay term that is multiplied by a small parameter $\varepsilon>0$.  This yields the equation
\begin{equation}
	x'(t) = f(x(t)) + \varepsilon\int_{-\infty}^t M(t-s) x(s)\, ds,  \label{eq:DDE}
\end{equation}
where $M(s) = \{M_{ij}(s)\}_{i,j=1}^d$ is a time dependent matrix.

Rearranging, we obtain
\begin{align*}
	& x'(t) = f(x(t)) + \varepsilon\int_{-\infty}^t M(t-s) (x(s)-x(t)+x(t))\, ds \\
	& \ \ \ \ \ \ \ \  = f(x(t)) + \varepsilon\left( \int_{-\infty}^t M(t-s)\, ds \right) x(t) +  \varepsilon\int_{-\infty}^t M(t-s)(x(s)-x(t))\, ds.
\end{align*}
After change of variables in time integrals we obtain
\begin{align*}
	& x'(t)  = f(x(t)) +  \varepsilon\left( \int_{0}^\infty M(s)\, ds \right) x(t)+\varepsilon\int_{-\infty}^t M(t-s)(x(s)-x(t))\, ds.
\end{align*}
 The Dafermos transform consists in using the correspondence between  functions $x:(-\infty,t] \to \R^d$ and pairs $(\eta,x_0)$ with $x_0\in \R^d$ and $\eta:[0,\infty) \to \R^d$ where $x_0=x(t)$ and $\eta(s) = x(t-s)-x(t)$. The inverse of this  correspondence maps  the pair $(\eta,x_0)$ to the function $x(s) = x_0 + \eta(t-s)$ for $s\leq t$. Thus, we  introduce the new variable
 $\eta^t:[0,\infty) \to \mathbb{R}^d$ as $\eta^t(s) = x(t-s)-x(t)$.
  Taking into account the influence of the initial data, we define
\begin{equation}\label{main2}
	\eta^t(s) = \begin{cases} x(t-s) - x(t)\ \ \textrm{for}\ \ s\leq t\\
		x(t-s) - x(t)=x_0+\eta^0(s-t)-x(t)\ \ \textrm{otherwise}.
	\end{cases}
\end{equation}
Using this variable, the equation of the problem takes the form
\begin{align}
	& x'(t)  = f^\varepsilon(x(t)) +  \varepsilon\int_{0}^\infty M(s)(x(t-s)-x(t))\, ds = f^\varepsilon(x(t)) +  \varepsilon\int_{0}^\infty M(s) \eta^t(s)\, ds ,\label{eq:main}
\end{align}
where $f^\varepsilon(x) = f(x) + \varepsilon\left( \int_{0}^\infty M(s)\, ds \right) x$, and equation \eqref{main2} governs the evolution of $\eta^t$.

\subsection{Energy inequality for the variable $\eta$.}\label{sec:memo_prop} The following lemma plays a crucial role in passing from ODE \eqref{eq:ODE} to \eqref{eq:DDE} as it shows that the tail $\eta$ can be treated
as a "contracting" direction from the point of view of geometric methods in dynamics. This is a crucial fact from the point of view of the tools used in this paper.

\begin{lem}\label{lem:eta}
	Let Assumption \ref{asm1} hold and let $x\in  C^1([0,\infty))^d$. Moreover let $\eta^t$ for $t\in [0,\infty)$ be given by  \eqref{main2} with the initial data  $x_0\in \mathbb{R}^d$ and $\eta^0\in L^2_A(\mathbb{R}^+)^d$. Then $\eta\in C^1([0,\infty);L^2_A(\R^+)^d)$,
	\begin{equation}\label{eq:1}
		\frac{d}{d t}\|\eta^t\|^2  + C_A\|\eta^t\|^2 \leq  - 2\left(\int_0^\infty A(s)\eta^t(s)\, ds,x'(t)\right),
	\end{equation}
	and
	\begin{equation}\label{eq:2}
		\|\eta^{t_2}\|^2 \leq e^{-C_A(t_2-t_1)}\|\eta^{t_1}\|^2 - 2 e^{-C_At_2} \int_{t_1}^{t_2}e^{C_At}\left(\int_0^\infty A(s)\eta^t(s)\, ds,x'(t)\right) \, dt\ \ \textrm{for} \ \ t_1< t_2.
	\end{equation}
\end{lem}
\begin{proof}
	Let $\eta^0 \in L^2_A(\R^+)^d$,  $x_0\in \R^d$ and $x \in C^1([0,\infty))^d$. Define $x(-s) = x_0+\eta^0(s)$ for $s> 0$. If $\eta^t$ is given by \eqref{main2} then the squared norm of $\eta^t$ is given by
	$$
	\|\eta^t\|^2 = \int_0^\infty (A(s)(x(t-s)-x(t)), (x(t-s)-x(t)))\, ds = \int_{-\infty}^t (A(t-s)(x(s)-x(t)), (x(s)-x(t)))\, ds.
	$$
	Let $t\geq 0$ and $h>0$. We calculate the right derivative of the above squared norm with respect to $t$.
	\begin{align*}
		& \frac{\|\eta^{t+h}\|^2 - \|\eta^t\| ^2}{h}  = \frac{1}{h}\int_{t}^{t+h} (A(t+h-s) (x(s)-x(t+h)),(x(s)-x(t+h)))\, ds \\
		& \ \ \ \  +\int_{-\infty}^t \left(\frac{A(t+h-s)-A(t-s)}{h}(x(s)-x(t)),(x(s)-x(t))\right)\, ds\\
		& \ \ \ \  + h\left(\int_{-\infty}^t A(t+h-s)\, ds \frac{x(t)-x(t+h)}{h},\frac{x(t)-x(t+h)}{h}\right)\\
		& \ \ \ \ + 2 \left(\int_{-\infty}^tA(t+h-s)(x(s)-x(t))\, ds, \frac{x(t)-x(t+h)}{h}\right).
	\end{align*}
	Passing to the limit with $h\to 0^+$, using the mean value theorem for integrals, the first term in the above sum tends to zero. Moreover, the limit of the third term is zero. In the second and fourth term we use the Lebesgue dominated convergence theorem to pass to the limit, whence
	\begin{align*}
		& \lim_{h\to 0^+}\frac{\|\eta^{t+h}\|^2 - \|\eta^t\| ^2}{h} \\
		& \ \ = \lim_{h\to 0^+}\int_{0}^\infty \left(\frac{A(s+h)-A(s)}{h}(x(t-s)-x(t)),(x(t-s)-x(t))\right)\, ds\\
		& \ \ \ \ \ \ \ - 2 \left(\lim_{h\to 0^+}\int_{0}^\infty A(s+h)(x(t-s)-x(t))\, ds, x'(t)\right) \\
		& \ \ = \int_{0}^\infty \left(\frac{d A(s)}{ds}\eta^t(s),\eta^t(s)\right)\, ds - 2 \left(\int_{0}^\infty A(s)\eta^t(s)\, ds, x'(t)\right).
	\end{align*}
	Similar calculation for $t>0$ and $h<0$ leads to the left derivative for $t>0$.
	Hence
	$$
	\frac{d}{dt}\|\eta^t\|^2 = \int_{0}^\infty \left(\frac{d A(s)}{ds}\eta^t(s),\eta^t(s)\right)\, ds - 2 \left(\int_{0}^\infty A(s)\eta^t(s)\, ds, x'(t)\right),
	$$
	and the assertion \eqref{eq:1} follows by Assumption \ref{asm1}.
	After multiplication by the integrating factor $e^{C_At}$ we deduce
	$$
	\frac{d}{d t} e^{C_At} \|\eta^t\|^2\, dt \leq - 2e^{C_At}\left(\int_0^\infty A(s)\eta^t(s)\, ds,x'(t)\right)
	$$
	Integrating from $t_1$ to $t_2$ we obtain \eqref{eq:2}.
\end{proof}
Similar argument leads to the following result
\begin{lem}\label{lem:diff_ogon}
	Let $x,y\in C^1([0,\infty)]$ and let $\eta,\xi$ be given by \eqref{main2} for $x$ and $y$, respectively.  Then
	\begin{equation}\label{eq:3}
		\frac{d}{d t}\|\eta^t-\xi^t\|^2  + C_A\|\eta^t-\xi^t\|^2 \leq  - 2\left(\int_0^\infty A(s)(\eta^t(s)-\xi^t(s))\, ds,(x(t)-y(t))'\right)
	\end{equation}
	and
	\begin{align*}
		%\label{eq:4}
		& \|\eta^{t_2}-\xi^{t_2}\|^2 \leq e^{-C_A(t_2-t_1)}\|\eta^{t_1}-\xi^{t_1}\|^2 \\
		& \qquad - 2 e^{-C_At_2} \int_{t_1}^{t_2}e^{C_At}\left(\int_0^\infty A(s)(\eta^t(s)-\xi^t(s))\, ds,x'(t)-y'(t)\right) \, dt\ \ \textrm{for} \ \ t_1< t_2.
	\end{align*}
\end{lem}

\subsection{Lyapunov function for the problem governed by \eqref{eq:DDE} and \eqref{main2}.}

The goal of this subsection is to construct a Lyapunov function which is uniform with respect to small $\varepsilon \geq 0$, i.e. which is valid not only for the ODE which is a gradient system (i.e. with $f=\nabla F$) but also  for \eqref{eq:DDE} with sufficiently small $\varepsilon$. To this end let $\mathcal{T}$ be a time interval either equal to $[0,T)$ or $[0,\infty)$ and assume that  $(\eta,x)\in  C^1(\mathcal{T};L^2_A(\R^+)^d)\times C^1(\mathcal{T})^d $ solve \eqref{eq:DDE}-- \eqref{main2}. We have the following lemma.

\begin{lem}
	\label{lem:LapFunc1}
	There exists $E_0 > 0$ such that for every $E\in (0,E_0)$ there exists $\varepsilon_0(E)>0$ such that for every $\varepsilon\in [0,\varepsilon_0)$  we have
	\begin{align}\label{eq:10c}
		&\frac{d}{d t}\left(E\|\eta^t\|^2-2F(x(t)) - \varepsilon\left(\int_0^\infty M(s)\, ds\,  x(t),x(t)\right)\right) + |x'(t)|^2 +
		%	\frac{C^2}{8\overline{D}^2\int_0^\infty\|A(s)\|\, ds}
		E\frac{C_A}{4}\|\eta^t\|^2  \leq   0 \ \ \text{for}\ t\in \mathcal{T}.
	\end{align}
\end{lem}
\begin{proof}
	Multiply \eqref{eq:main} by $2x'(t)$. Then
	$$
	2|x'(t)|^2 = 2(f^\varepsilon(x(t)),x'(t)) + 2\varepsilon\left(\int_0^\infty M(s)\eta^t(s)\, ds, x'(t)\right).
	$$
	Adding this equation to the inequality \eqref{eq:1} from Lemma~\ref{lem:eta} multiplied by $E>0$ we obtain
	\begin{align}\label{eq:9b}
		&\frac{d}{d t}E\|\eta^t\|^2 + 2|x'(t)|^2 + EC_A\|\eta^t\|^2 \leq  2(f(x(t)),x'(t)) \nonumber\\
		& \ \  + 2 \varepsilon\left(\int_0^\infty M(s)\, ds\  x(t),x'(t)\right) - 2\left(\int_0^\infty \left(E A(s) - \varepsilon M(s)\right)\eta^t(s)\, ds, x'(t)\right).
	\end{align}
	Using the fact that $f = \nabla F$ and $\int_0^\infty M(s)\, ds $ is symmetric we get
	\begin{align}\label{eq:10b}
		&\frac{d}{d t}\left(E\|\eta^t\|^2-2F(x(t)) - \varepsilon\left(\int_0^\infty M(s)\, ds\,  x(t),x(t)\right)\right) + 2|x'(t)|^2 + EC_A\|\eta^t\|^2\\
		& \ \  \leq   -2\left(\int_0^\infty \left(E A(s) - \varepsilon M(s)\right)\eta^t(s)\, ds, x'(t)\right). \nonumber
	\end{align}
	Choosing $\delta>0$ we estimate the term on the right-hand side as
	$$2\left(\int_0^\infty \left(E A(s) - \varepsilon M(s)\right)\eta^t(s)\, ds, x'(t)\right)\leq \delta |x'(t)|^2 + \frac{1}{\delta}\left|\int_0^\infty \left(E A(s)-\varepsilon M(s)\right)\eta^t(s)\, ds\right|^2$$.
	The last term can be estimated as
	\begin{align*}
		&\frac{1}{\delta}\left|\int_0^\infty \left(E A(s)-\varepsilon M(s)\right)\eta^t(s)\, ds\right|^2 \leq \frac{1}{\delta}\left(\int_0^\infty \sqrt{\|E A(s)-\varepsilon M(s)\|}\sqrt{\|E A(s)-\varepsilon M(s)\|}\, |\eta^t(s)|\, ds\right)^2\\
		& \ \leq \frac{1}{\delta}\int_0^\infty \|E A(s)- \varepsilon M(s)\|\, ds \int_0^\infty \|E  A(s)-\varepsilon M(s)\|\, |\eta^t(s)|^2\, ds.
	\end{align*}
	Moreover,
	$$
	\int_0^\infty \|E A(s)- \varepsilon M(s)\|\, ds \leq E \int_0^\infty\|A(s)\|\, ds + \varepsilon \int_0^\infty\|M(s)\|\, ds,
	$$
	and, using Assumptions \ref{asm:a2} and \ref{asm:M}
	\begin{align*}
		& \int_0^\infty \|E A(s)-\varepsilon M(s)\|\, |\eta^t(s)|^2\, ds \leq E\int_0^\infty \|A(s)\|\, |\eta^t(s)|^2\, ds + \varepsilon\int_0^\infty \|M(s)\|\, |\eta^t(s)|^2\, ds\\
		& \ \  \leq \left(E\overline{D_A}^2 + \varepsilon \overline{D_M}^2\right)\|\eta^t\|^2
	\end{align*}
	Choosing $\delta = 1$ we get
	\begin{align*}
		&\frac{d}{d t}\left(E\|\eta^t\|^2-2F(x(t)) - \varepsilon\left(\int_0^\infty M(s)\, ds\,  x(t),x(t)\right)\right) + |x'(t)|^2 + EC_A\|\eta^t\|^2\\
		& \ \  \leq   \left(E \int_0^\infty\|A(s)\|\, ds + \varepsilon \int_0^\infty\|M(s)\|\, ds\right)\left(E\overline{D_A}^2 + \varepsilon\overline{D_M}^2\right)\|\eta^t\|^2.
	\end{align*}
	Moving all terms to the left, the constant in front of $\|\eta^t\|^2$ is equal to
	$$
	- E^2\overline{D_A}^2 \int_0^\infty\|A(s)\|\, ds + E \left(C_A - \varepsilon\overline{D_M}^2 \int_0^\infty\|A(s)\|\, ds - \varepsilon\overline{D_A}^2\int_0^\infty\|M(s)\|\, ds \right) - \varepsilon^2 \overline{D_M}^2 \int_0^\infty\|M(s)\|\, ds.
	$$
  We want this expression to be non-negative. It can be rewritten as
	$$
	-E^2 G_1 + E (G_2 - \varepsilon G_3) - \varepsilon^2 G_4,
	$$
	where $G_1, G_2, G_3, G_4$ are positive constants. Take $E_0 = \frac{G_2}{2G_1}$. If $E\in (0,E_0)$, then
	$$
	-E^2 G_1 + E (G_2 - \varepsilon G_3) - \varepsilon^2 G_4 \geq EG_2 - E\frac{G_2}{2} - \varepsilon\frac{G_3G_2}{2G_1} - \varepsilon^2G_4 = E\frac{G_2}{2} - \varepsilon\frac{G_3G_2}{2G_1} - \varepsilon^2G_4.
	$$
	Now take $\varepsilon_0(E)$ such that
	$$
	\varepsilon\frac{G_3G_2}{2G_1} + \varepsilon^2G_4 \leq E\frac{G_2}{4},
	$$
	if only $\varepsilon\in (0,\varepsilon_0)$.
	This means that
	$$
	-E^2 G_1 + E (G_2 - \varepsilon G_3) - \varepsilon^2 G_4 \geq E\frac{G_2}{4} = E\frac{C_A}{4}.
	$$
	%If $\varepsilon = 0$ the maximum of this function is equal to $\frac{C^2}{4\overline{D}^2\int_0^\infty\|A(s)\|\, ds}$ and is attained for $E=E_0 = \frac{C}{2\overline{D}^2 \int_0^\infty\|A(s)\|\, ds}$. A straightforward calculation shows that there exists $\varepsilon_0 > 0$ and $E_0>0$ such that for every $\varepsilon\in [0,\varepsilon_0)$ the constant is greater than  $\frac{C^2}{8\overline{D}^2\int_0^\infty\|A(s)\|\, ds}$.
	The proof is complete.
\end{proof}
As a consequence of the above Lemma we obtained the following Lyapunov function  
$L_\varepsilon:L^2_A(\mathbb{R}^+)^d\times \mathbb{R}^d \to \R $ valid for every $E\in (0,E_0)$ and for every $\varepsilon \in [0,\varepsilon_0(E))$
\begin{equation}\label{eq:lyapunov}
	L_\varepsilon(\eta,x)=E\|\eta\|^2-2F(x) - \varepsilon\left(\int_0^\infty M(s)\, ds\,  x,x\right)
\end{equation}
Observe that  \eqref{diss} implies that
\begin{eqnarray} \label{eq:Lap-lbnd}
  L_\varepsilon(\eta,x)  \geq E\|\eta\|^2 + \left(2 \gamma  - \varepsilon\int_0^\infty \|M(s)\|\, ds\,\right)|x|^2 - 2 \delta.
\end{eqnarray}

%As the next Lemma shows, we can also obtain the following simpler family of Lyapunov functions dependent on $\varepsilon$
%\begin{equation}
%L(x(t),\eta^t)=\varepsilon\|\eta^t\|^2 - 2F(x(t)).
%\end{equation}
%\begin{lem}
%\label{lem:LapFunc2}
%    There holds the bound
%\begin{equation}\label{eq:9}
%\frac{d}{d t}\varepsilon\|\eta^t\|^2 - 2(f(x(t)),x'(t)) \leq -\varepsilon C\|\eta^t\|^2 - 2|x'(t)|^2.
%\end{equation}
%In particular if $f = \nabla F$ then
%\begin{equation}\label{eq:10}
%\frac{d}{d t}(\varepsilon\|\eta^t\|^2 - 2F(x(t))) \leq -\varepsilon C\|\eta^t\|^2 - 2|x'(t)|^2.
%\end{equation}
%	\end{lem}
%	\begin{proof}
%Multiply \eqref{eq:main} by $2x'(t)$. Then
%	$$
%	2|x'(t)|^2 = 2(f(x(t)),x'(t)) + 2\varepsilon\left(\left(\int_0^\infty M(s) ds \right)x(t),x'(t)\right) + 2\varepsilon\left(\int_0^\infty M(s)\eta^t(s)\, ds, x'(t)\right).
	%$$
%	Adding this equation to the inequality \eqref{eq:1} from Lemma~\ref{lem:eta} multiplied by $\varepsilon$ we immediately obtain .
%\begin{eqnarray*}
%  \frac{d}{d t}\varepsilon\|\eta^t\|^2  + C\varepsilon\|\eta^t\|^2 + 2|x'(t)|^2 \leq  - 2\varepsilon\left(\int_0^\infty A(s)\eta^t(s)\, ds,x'(t)\right) + 2(f(x(t)),x'(t)) \\
%  + 2\varepsilon\left(\left(\int_0^\infty M(s) ds \right)x(t),x'(t)\right) + 2\varepsilon\left(\int_0^\infty M(s)\eta^t(s)\, ds, x'(t)\right).
%\end{eqnarray*}
%\textbf{PZ: nie wydaje sie aby wyszlo co trzeba, nawet gdyby $A=M$}
%\end{proof}

\subsection{Solution boundedness and $C^0$ dependence on initial data.}

 We skip the proof of the existence and uniqueness of local solution for every initial data $(\eta^0,x_0)\in  L^2_A(\R^+)^d\times R^d$. The fact that every solution can be extended to a global one follows from the next lemma which is a consequence of the existence of a Lyapunov function obtained in Lemma \ref{lem:LapFunc1}.
We prove that every solution $(x,\eta)$ of \eqref{eq:DDE} and \eqref{main2} is globally bounded by a constant dependent on the initial data.
\begin{lem}\label{le:bound}
Assume \eqref{gradient} and \eqref{diss}.  Under assumptions of Lemma \ref{lem:LapFunc1},  if only
\begin{equation}\label{boundepsilon}
\varepsilon\int_0^\infty\|M(s)\|\, ds < 2\gamma,
\end{equation}
then every solution is bounded uniformly on bounded sets of initial data.
\end{lem}
\begin{proof}
From
the Lyapunov function \eqref{eq:lyapunov} and the bound \eqref{eq:Lap-lbnd} we obtain that
$$
 C(|x_0|, \|\eta^0\|) \geq L_\varepsilon(\eta^0,x_0) \geq L_\varepsilon(\eta^t,x(t)) \geq E\|\eta^t\|^2 + \left(2 \gamma  - \varepsilon\int_0^\infty \|M(s)\|\, ds\,\right)|x(t)|^2 - 2 \delta,
$$
where $C(\cdot,\cdot)$ is a continuous function independent of $\varepsilon$. This yields the assertion of the lemma.
\end{proof}
 In the next result, we prove the Lipschitz continuous dependence of the solution on the initial data.
	
\begin{lem}\label{lem:difference}
	Assume \eqref{gradient} and \eqref{diss}. There exists $\varepsilon_0$ such that for every $\varepsilon\in [0,\varepsilon_0]$ if $(\eta^t,x(t))$ and $(\xi^t,x(t))$ are two solutions with the initial data $(\eta^0,x_0)$ and $(\xi^0,y_0)$, respectively, then for every  $T>0$ there exists a constant $L(T)$ such that for every $t\in [0,T]$ we have
	$$
	|x(t)-y(t)| + \|\eta^t-\xi^t\|\leq L(T)(|x_0-y_0| + \|\eta^0-\xi^0\|).
	$$
\end{lem}
\begin{proof}
	In the proof by $D$ we will denote constants (which can vary from line to line) dependent on the initial data for both problems and by $C_i$ constants independent on these data. Subtracting \eqref{eq:main} for the two solutions we obtain
	\begin{equation}\label{eq:difference}
	(x(t)-y(t))'=f(x(t))-f(y(t)) + \varepsilon\int_0^\infty M(s)\, ds (x(t)-y(t)) + \varepsilon\int_0^\infty M(s)(\eta^t(s)-\xi^t(s))\, ds.
	\end{equation}
We multiply the above equation by $(x(t)-y(t))$, whence
	$$
	\frac{1}{2}\frac{d}{dt}|x(t)-y(t)|^2 \leq |f(x(t))-f(y(t))|\, |x(t)-y(t)| +  C_1 |x(t)-y(t)|^2 +  C_2 \|\eta^t-\xi^t\|\, |x(t)-y(t)|. 
	$$
	As from Lemma \ref{le:bound} the set  $\{\textrm{conv}\{x(t),y(t)\}\,:\ t\geq 0\}$ is bounded we deduce that
	\begin{equation}\label{eq:ineqone}
	\frac{d}{dt}|x(t)-y(t)|^2 \leq D|x(t)-y(t)|^2+C_3\|\eta^t-\xi^t\|^2.
	\end{equation}  
	From Lemma \ref{lem:diff_ogon} we obtain
	$$
	\frac{d}{dt}\|\eta^t-\xi^t\|^2 + C_A\|\eta^t-\xi^t\|^2\leq C_4\|\eta^t-\xi^t\|\ \ |(x(t)-y(t))'|.
	$$
	Substituting \eqref{eq:difference} and proceeding similarly as in the proof of \eqref{eq:ineqone}
	it follows that  $$
	\frac{d}{dt}\|\eta^t-\xi^t\|^2 + C_A\|\eta^t-\xi^t\|^2\leq D\ |x(t)-y(t)|^2 + C_5 \|\eta^t-\xi^t\|^2.
	$$
	Combining this inequality with \eqref{eq:ineqone}, we obtain
	$$
	\frac{d}{dt}(\|\eta^t-\xi^t\|^2+|x(t)-y(t)|^2) \leq D |x(t)-y(t)|^2 + C_6\|\eta(t)-\xi(t)\|^2,
	$$
	which, by the Gronwall lemma yields the required assertion.
\end{proof}

\subsection{Existence and uniform boundedness of global attractors.}The question which we address in the remaining part of the article is the following. Assume that $x'=f(x)$ is a Morse--Smale  system. The Morse--Smale property in our case means that the vector field $f$ has a finite number of hyperbolic equilibria such that the intersections of stable and unstable manifolds are transversal. If $\varepsilon>0$ is small, can we say that the problem with distributed memory has the same structure of the global attractor as the ODE?

The families of maps $\{S^\varepsilon(t)\}_{t\geq 0}:L^2_A(R^+)^d\times \R^d \to L^2_A(R^+)^d\times \R^d$ denote the semiflows the govern the solutions of the problem \eqref{eq:main}.

We prove that assumptions \eqref{gradient}, \eqref{diss}, and \eqref{diss2} imply that for $\varepsilon\in [0,\varepsilon_0]$ problems  have global attractors  $\mathcal{A}_\varepsilon \subset  L^2_A(\R^+)^d\times \R^d$ such that
\begin{equation}\label{eq:bound}
\bigcup_{\varepsilon\in[0,\varepsilon_0)} \mathcal{A}_\varepsilon\ \ \ \textrm{is bounded in}\ L^2_A(\R^+)^d \times \R^d.
\end{equation}
We denote the set of equilibria by
$$
\mathcal{E}_\varepsilon = \left\{ (0,x)\in L^2_A(\R^+)^d\times \R^d\,:\ \nabla F(x) + \varepsilon\int_0^\infty M(s)\, ds\, x = 0\right\}.
$$
Note that for all equilibria $\eta=0$, i.e. the memory variable must be equal to zero. 
%For the next lemma we assume that sets of equilibria are uniformly bounded. We will later prove that this bound indeed holds and follows from the fact that, for small $\varepsilon$, equilibria for $\varepsilon>0$ belong to isolating blocks for equilibria for $\varepsilon=0$ and there are no other equilibria.
\begin{lem} \label{lem:atr-bnd}
	Assume \eqref{gradient}, \eqref{diss}, and \eqref{diss2}.  Then there exists $\varepsilon_0$ such that  for every $\varepsilon\in [0,\varepsilon_0]$ the problems governed by \eqref{eq:main}-\eqref{main2} have global attractors $\mathcal{A}_\varepsilon$, that satisfy \eqref{eq:bound}.
\end{lem}
\begin{proof}
	We first show that if $\varepsilon$ is sufficiently small, then $\mathcal{E}_\varepsilon$ is bounded by a bound independent of $\varepsilon$. 
	We will show that if $(0,x)\in \mathcal{E}_\varepsilon$, then $|x|< R$, where $R$ is a constant from \eqref{diss2}. Indeed assume that $|x|\geq R$. Then
	$$
	C_F|x|^2 \leq -\nabla F(x)\cdot x = \varepsilon \left(\int_0^\infty M(s)\, ds\, x,x\right) \leq \varepsilon_0 \int_0^\infty \|M(s)\|\, ds |x|^2.    
	$$ 
	Hence, it suffices to take 
	$$\varepsilon_0 < \frac{\int_0^\infty\|M(s)\|ds}{C_F}$$ to arrive at a contradiction.  
We continue the argument by using the asymptotic compactness
result for the memory term whose proof is postponed to Appendix 1. We  use Lemma \ref{lemma:compact} from that appendix. The asymptotic compactness of the semiflow established in that lemma implies  that for every bounded set $\mathcal{B}\subset  L^2_A(\R^+)^d\times \R^d$ its $\omega$-limit set $\omega(\mathcal{B})$ (see Definition \ref{df:wB} in Appendix 1) is nonempty, compact and attracts $\mathcal{B}$ in the sense of Hausdorff semidistance in $ L^2_A(\R^+)^d\times \R^d$, cf. Lemma \ref{as:com}.

The argument now follows the lines of the proof of Theorem A.3 in \cite{Pata}.
%, while $x$ belongs to the isolating set around the equilibrium for $\varepsilon=0$.
Lemma \ref{lemma:compact}
as well as the existence of the Lyapunov function imply that for every initial data $(\eta^0,x_0)$ we can find an equilibrium $(0,x^*)$ such that
$S^\varepsilon(t)(\eta^0,x_0)\to (0,x^*)$ as $t\to \infty$.

Recall, that $L_\varepsilon(\eta,x)$  defined by \eqref{eq:lyapunov}  is Lyapunov function and we define   
$$
\mathcal{C}_\varepsilon = \left\{(\eta,x)\in L^2_A( \R^+)^d \times \R^d \,:\ L_\varepsilon(\eta,x) < \max_{(\xi,y)\in \mathcal{E}_\varepsilon} L_{\varepsilon}(\xi,y) + 1\right\}.
$$
As the set $\bigcup_{\varepsilon\in [0,\varepsilon_0]}\mathcal{E}_\varepsilon$ is bounded, so, from the bound \eqref{eq:Lap-lbnd}, is the set $\bigcup_{\varepsilon\in [0,\varepsilon_0]}\mathcal{C}_\varepsilon$.

If we fix $\mathcal{B}$, then there exists time $t^*(\mathcal{B})$ such that $S^\varepsilon(t)\omega(\mathcal{B})\subset \mathcal{C}_\varepsilon$ for $t\geq t^*$. Indeed, by continuity of $S^\varepsilon(t)$ for every $p\in \omega(\mathcal{B})$ there exists a neighborhood $\mathcal{U}_p$ and $t_p$ such that  $S^\varepsilon(t_p)\mathcal{U}_p \subset \mathcal{C}_{\varepsilon}$. As  $\mathcal{C}_{\varepsilon}$ is positively invariant, the inclusion $S^\varepsilon(t)\mathcal{U}_p \subset \mathcal{C}_{\varepsilon}$ holds for every $t\geq t_p$.
Sets $\{\mathcal{U}_p\}_{p\in \omega(\mathcal{B})}$ are open cover of  $\omega(\mathcal{B})$. We extract finite subcover,
	$\{\mathcal{U}_{p_n}\}_{n=1}^N$ whereas $t^* = \max\{t_{p_1},\ldots,t_{p_n}\}$. Since there exists a function $\psi(t) > 0$ such that
	$\lim_{t\to\infty}\textrm{dist}_{L^2_A(\R^+)^d\times \R^d}(S^\varepsilon(t)\mathcal{B},\omega(\mathcal{B})) \leq \lim_{t\to\infty}\psi(t) = 0$, for every $t$ and $(\eta,x)\in \mathcal{B}$ there exists $k(t) \in \omega(\mathcal{B})$ and $q(t)$ such that $S^\varepsilon(t)(\eta,x) = k(t) + q(t)$ and  $\|q(t)\|_{L^2_A(\R^+)^d\times \R^d} \leq 2\psi(t)$. Now
	$S^\varepsilon(t+t^*)(\eta,x) = S^\varepsilon(t^*)k(t) + S^\varepsilon(t^*)(k(t) + q(t))-S^\varepsilon(t^*)k(t)$ and  $S^\varepsilon(t^*)k(t) \in \mathcal{C}_\varepsilon$. Moreover $S^\varepsilon(t^*)$ is continuous and hence it is uniformly continuous in a neighbourhood of a compact set, therefore, for $t$ large enough
	$$
	\|S^\varepsilon(t^*)(k(t) + q(t))-S^\varepsilon(t^*)k(t)\|_{ L^2_A(\R^+)^d \times \R^d}\leq 1.
	$$
	This means that, for $t$ large enough
	$
	S_\varepsilon(t+t^*)(\eta,x)
	$ belongs to the ball centered at zero and with radius $\sup_{(\eta,x)\in \mathcal{C}_\varepsilon} \|(\eta,x)\|_{L^2_A(\R^+)^d\times \R^d} + 1$. Together with Lemma   \ref{lemma:compact} it is enough to guarantee the existence of the global attractor $\mathcal{A}_\varepsilon$ and the bound \eqref{eq:bound}.
\end{proof}
\subsection{Estimate for the difference of two solutions.}
In the next lemma we compare two solutions for different values of $\varepsilon$ and initial data. We assume that $\varepsilon\in [0,\varepsilon_0]$ with $\varepsilon_0$ being sufficiently small.
	\begin{lem}\label{lem_diff}
		Consider two solutions: one  $(\eta,x)$ of problem with $\varepsilon_1$  and the initial data $(\eta^{0},x_0)$ and another one $(\xi,y)$ of the problem with $\varepsilon_2$ and the initial data $(\xi^0,y_0)$.
		Then
		 	\begin{equation}\label{eq:close_1}
		 	|x(t)-y(t)|+\|\eta^t-\xi^t\|\leq Ce^{Ct}(\|\eta^0-\xi^0\| + |x_0-y_0|   + |\varepsilon_1-\varepsilon_2|),
		 \end{equation}
		 for every $t\geq 0$ where the constants $C$   depending on the initial data $(\eta^{0},x_0)$ and $(\xi^0,y_0)$ and are bounded on bounded sets of initial data.
	\end{lem}
	\begin{proof}
	Subtracting the equations for $x$ and $y$ we obtain
	\begin{align*}
	& (x(t)-y(t))' = f(x(t))-f(y(t)) + \varepsilon_1\left( \int_{0}^\infty M(s)\, ds \right) (x(t)-y(t))+(\varepsilon_1-\varepsilon_2)\left( \int_{0}^\infty M(s)\, ds \right) y(t) \\
	& \ + \varepsilon_1\int_{0}^\infty M(s) (\eta^t(s)-\xi^t(s))\, ds + (\varepsilon_1-\varepsilon_2)\int_{0}^\infty M(s) \xi^t(s)\, ds.
	\end{align*}
	Lyapunov function \eqref{eq:lyapunov} implies that sets $\{\textrm{conv}\{x(t),y(t)\}\, :\ t\geq 0\}$ and $\{\textrm{conv}\{\eta^t,\xi^t\}\, :\ t\geq 0\}$ are bounded by constants depending on the initial data of the problem. We denote the generic constant depending on the initial data by $C$. Multiplying the above equation by $x(t)-y(t)$ we obtain  \begin{align*}
	& \frac{d}{dt}|x(t)-y(t)|^2 \leq C
	\  |x(t)-y(t)|^2 + C |\varepsilon_1-\varepsilon_2|\,|x(t)-y(t)|+C\varepsilon_1\|\eta^t-\xi^t\|\, |x(t)-y(t)|\\
	& \ \ \ \leq C
	\  |x(t)-y(t)|^2 + C |\varepsilon_1-\varepsilon_2| ^2 +C\|\eta^t-\xi^t\|^2.
	\end{align*}
	Using \eqref{eq:3} it follows that
	\begin{align*}
	&\frac{d}{dt}\|\eta^t-\xi^t\|^2 + C_A\|\eta^t-\xi^t\|^2 \leq C\|\xi^t-\eta^t\|\, |x'(t)-y'(t)|\\
               & \ \ \leq C\|\xi^t-\eta^t\|\, |x(t)-y(t)| + C|\varepsilon_1-\varepsilon_2|\, \|\xi^t-\eta^t\|+C\varepsilon_1\|\eta^t-\xi^t\|^2.
	\end{align*}
	After straightforward calculations, and for sufficiently small $\varepsilon_0$,
	$$
	\frac{d}{dt}\|\eta^t-\xi^t\|^2\leq  C|x(t)-y(t)|^2 + C|\varepsilon_1-\varepsilon_2|^2,
	$$
	whence
		$$
	\frac{d}{dt}(|x(t)-y(t)|^2+\|\eta^t-\xi^t\|^2) \leq C
	\  (|x(t)-y(t)|^2 + \|\eta^t-\xi^t\|^2 +|\varepsilon_1-\varepsilon_2|^2),
	$$
	which yields the assertion of the lemma.
	\end{proof}
\subsection{Variational equation and $C^1$ dependence of the solution on initial data.} In the next lemma we characterize the derivative of the flow with respect to the initial data
	\begin{lem}\label{variational}
		Consider the mapping  $ L^2_A(\R^+)^d \times \R^d\ni (\eta^0,x_0)\mapsto (\eta(t),x_t)=S^\varepsilon(t)(\eta^0,x_0)$ defining the solutions of  \eqref{eq:main}--\eqref{main2}. The mapping $S^\varepsilon(t)$ is Fr\'{e}chet differentiable and its derivative at $(\eta^0,x_0)$ is defined as the linear mapping that assigns to $(\xi^0,w_0)$ the solution of the variational problem
		\begin{align}
			& w'(t) =  Df(x(t)) w(t)  + \varepsilon  \left(\int_0^\infty M(s)\, ds\right)  w(t) + \varepsilon \int_0^\infty M(s) \theta^t(s)\, ds.\label{weq_glob_2} \\
		&
		\theta^t(s) = \begin{cases} w(t-s) - w(t)\ \ \textrm{for}\ \ s\leq t\\
			w_0+\xi^0(s-t) - w(t)\ \ \textrm{otherwise},
		\end{cases}\label{variational:2}
		\end{align}
		with the initial data $w(0)=w_0$, $\theta^0=\xi^0$.
	\end{lem}
	\begin{proof}
	 We take two initial conditions $x_0, \overline{x}_0$ and $\eta^0, \overline{\eta}^0$ and call the corresponding solutions  $(\eta^t,x(t))$ and $( \overline{\eta}^t,\overline{x}(t))$. Their difference will be called
	$z(t) = \overline{x}(t)-x(t)$ and $\xi^t = \overline{\eta}^t-{\eta}^t$. They satisfy the equations
	\begin{align*}
		& z'(t) =  f(x(t)+z(t)) - f(x(t)) +  \varepsilon\left(\int_0^\infty M(s)\, ds\right)  z(t) + \varepsilon \int_0^\infty M(s) \xi^t(s)\, ds.
	\end{align*}
	and
	$$
	\xi^t(s) = \begin{cases} z(t-s) - z(t)\ \ \textrm{for}\ \ s\leq t\\
		 z_0+\xi^0(s-t) - z(t)\ \ \textrm{otherwise}.
	\end{cases}
	$$
	This motivates the  definition \eqref{weq_glob_2}--\eqref{variational:2} of a variational equation with unknowns $\theta^t$ and $w(t)$.
	Denote the difference $z-w = p$ and $\xi^t-\theta^t = \omega^t$. Then
	\begin{align*}
		& p'(t) =  f(x(t)+z(t)) - f(x(t)) - Df(x(t)) w(t)  + \varepsilon  \left(\int_0^\infty M(s)\, ds\right) p(t) + \varepsilon \int_0^\infty M(s) \omega^t(s)\, ds.
	\end{align*}

	$$
	\omega^t(s) = \begin{cases}p(t-s) - p(t)\ \ \textrm{for}\ \ s\leq t\\
		- p(t)\ \ \textrm{otherwise}.
	\end{cases}
	$$
	Rearranging the first equation and using the Taylor formula  with integral remainder we obtain  
	\begin{align*}
		& p'(t) =  Df(x(t))  p(t) +  \int_0^1 (1-\lambda)D^2f(x(t)+\lambda z(t))(z(t),z(t)) \, d\lambda + \varepsilon \int_0^t M(t-s)  p(s)\, ds.
	\end{align*}
	Integrating and using the fact that $x(t)$ and $z(t)$ are bounded, $D^2f$ is bounded on bounded sets and $p(0)=0$  we obtain
	$$
	|p(t)| \leq C\int_0^t |p(s)|ds +C \int_0^t |z(s)|^2\, ds + \varepsilon C \int_0^t \int_0^s |p(r)|\, dr\, ds \leq C \int_0^t |z(s)|^2\, ds +  C(1+\varepsilon t) \int_0^t |p(s)|\, ds.
	$$
	We need an estimate for $|z(s)|$. We have
	\begin{align*}
		& z'(t) =   \int_0^1 Df(x(t)+\lambda z(t)) \, d\lambda \, z(t)  +  \varepsilon \int_0^t M(s)  z(t-s)\, ds + \varepsilon \int_t^\infty M(s)  \, ds\,  z_0 + \varepsilon \int_t^\infty M(s) \xi^0(s-t)\, ds.
	\end{align*}
	Rewriting, we obtain
	\begin{align*}
		& z'(t) =  \int_0^1 Df(x(t)+\lambda z(t)) \, d\lambda \, z(t) +  \varepsilon \int_0^t M(t-s)  z(s)\, ds + \varepsilon \int_t^\infty M(s)  \, ds z_0 + \varepsilon \int_0^\infty M(s+t) \xi^0(s)\, ds.
	\end{align*}
	Now $x(t)$ is bounded and as is $z(t)$ because $\overline{x}(t)$ is attracted to the attractor and hence also bounded, and $Df$ is bounded on bounded sets.  We obtain
	$$
	|z(t)|\leq (1+\varepsilon Ct)|z_0| +\varepsilon Ct\|\xi^0\| + C (1+\varepsilon t)\int_0^t|z(s)|\, ds.
	$$
	By the Gronwall lemma
	$$
	|z(t)| \leq ((1+\varepsilon Ct)|z_0| +\varepsilon Ct\|\xi^0\|)e^{C t(1+\varepsilon t)}.
	$$
	This means that
	$$
	|z(t)|^2 \leq g(t)(|z_0|^2 + \|\xi^0\|^2),
	$$
	where by $g(t)$ we denote a generic increasing and continuous function of $t$.
	We deduce that
	$$
	|p(t)|\leq g(t) (|z_0|^2+\|\xi^0\|^2) + C(1+\varepsilon t)\int_0^t|p(s)|\ ds.
	$$
	By the Gronwall lemma
	$$
	|p(t)| \leq g(t) (|z_0|^2+\|\xi^0\|^2).
	$$
	Moreover,
	$$
	\|\omega^t\|^2 = \int_0^t (A(s)  p(t-s),p(t-s))\, ds - \int_0^\infty (A(s)p(t),p(t))\, ds,
	$$
	this means that
	$$
	\|\omega^t\|^2 \leq g(t)\left(|p(t)|^2+\int_0^t|p(s)|^2\, ds\right).
	$$
	We deduce that
	$$
	\|\omega^t\| \leq g(t)(|z_0|^2+\|\xi^0\|^2).
	$$
	We conclude that
	$$
	\lim_{|z_0|\to 0, \|\xi^0\|\to 0}\frac{|p(t)|+\|\omega^t\|}{|z_0|+\|\xi^0\|}\leq \lim_{|z_0|\to 0, \|\xi^0\|\to 0}\frac{g(t)(|z_0|^2+\|\xi^0\|^2)}{|z_0|+\|\xi^0\|}=0.
	$$
	This implies the Fr\'{e}chet differentiability of the flow and the fact that the derivative with respect to the initial data is the solution of the variational equation.
	\end{proof}
	
The following lemma implies the continuous dependence of the derivative with respect to the initial data on the parameter $\varepsilon$. 	
	
\begin{lem}\label{est:variational}
		Let $(\eta^{0,\varepsilon_1},x^{\varepsilon_1}_0)$ and $(\eta^{0,\varepsilon_2},x^{\varepsilon_2}_0)$ be the initial data for problems with $\varepsilon_1$ and $\varepsilon_2$, respectively. Moreover, let $(\theta^{0,\varepsilon_1},w^{\varepsilon_1}_0)$ and $(\theta^{0,\varepsilon_2},w^{\varepsilon_2}_0)$
	be the initial data for the variational problem \eqref{weq_glob_2}-\eqref{variational:2}. Then
	\begin{align*}
	& \left\|\frac{DS^{\varepsilon_2}(t)(\eta^{0,\varepsilon_2},x^{\varepsilon_2}_0)}{D(\eta,x)}(\theta^{0,\varepsilon_2},w^{\varepsilon_2}_0)
	-\frac{DS^{\varepsilon_1}(t)(\eta^{0,\varepsilon_1},x^{\varepsilon_1}_0)}{D(\eta,x)}(\theta^{0,\varepsilon_1},w^{\varepsilon_1}_0)\right\|_{ L^2_A(\R^+)^d\times \R^d}\\
	& \ \ \leq Ce^{Ct}(|\varepsilon_2-\varepsilon_1|+|w^{\varepsilon_2}_0-w^{\varepsilon_1}_0|+ \|\theta^{0,\varepsilon_2}-\theta^{0,\varepsilon_1}\|+|x^{\varepsilon_2}_0-x^{\varepsilon_1}_0| +\|\eta^{0,\varepsilon_2}-\eta^{0,\varepsilon_1}\|),
	\end{align*}
	for every $t\geq 0$, where $C$ depends on the initial data in a non-decreasing way and is bounded on bounded sets of initial data.
\end{lem}
\begin{proof}
	Denote
	$$
	\frac{DS^{\varepsilon_i}(t)(\eta^{0,\varepsilon_i},x^{\varepsilon_i}_0)}{D(\eta,x)}(\theta^{0,\varepsilon_i},w^{\varepsilon_i}_0) = (\theta^{\varepsilon_i,t},w^{\varepsilon_i}(t)).
	$$
		We have
	$$
	\frac{d}{dt}|w^{\varepsilon_2}(t)|^2 \leq C|w^{\varepsilon_2}(t)|^2 + C\|\theta^{\varepsilon_2,t}\|^2.
	$$
	Moreover, from \eqref{eq:1},
	$$
	\frac{d}{dt}\|\theta^{\varepsilon_2,t}\|^2 + C_A\|\theta^{\varepsilon_2,t}\|^2\leq C \|\theta^{\varepsilon_2,t}\|\,  |(w^{\varepsilon_2})'(t)|.
	$$
	It follows that
		$$
	\frac{d}{dt}\|\theta^{\varepsilon_2,t}\|^2\leq C |w^{\varepsilon_2}(t)|^2 + C\|\theta^{\varepsilon_2,t}\|^2,
	$$
	and
	$$
	\frac{d}{dt}(|w^{\varepsilon_2}(t)|^2+ \|\theta^{\varepsilon_2,t}\|^2)\leq C(|w^{\varepsilon_2}(t)|^2 + \|\theta^{\varepsilon_2,t}\|^2),
	$$
	whereas
	\begin{equation}\label{var_est1}
	|w^{\varepsilon_2}(t)|+ \|\theta^{\varepsilon_2,t}\| \leq Ce^{Ct}(|w^{\varepsilon_2}_0|+ \|\theta^{\varepsilon_2,0}\|).
	\end{equation}
	Now, we have the following equation for the difference between two solutions of variational equations along the equations on attractors	
	\begin{align*}
&	(w^{\varepsilon_2}(t)-w^{\varepsilon_1}(t))'=(Df(x^{\varepsilon_2}(t))-Df(x^{\varepsilon_1}(t)))w^{\varepsilon_2}(t) +  Df(x^{\varepsilon_1}(t))(w^{\varepsilon_2}(t)-w^{\varepsilon_1}(t))\\
& \ \  + (\varepsilon_2-\varepsilon_1) \int_0^\infty M(s) (s)\, ds w^{\varepsilon_2}(t) + \varepsilon_1 \int_0^\infty M(s) \, ds (w^{\varepsilon_2}(t)-w^{\varepsilon_1}(t)) \\
& \ \ + (\varepsilon_2-\varepsilon_1) \int_0^\infty M(s) \theta^{\varepsilon_2,t}(s)\, ds +
\varepsilon_1 \int_0^\infty M(s) (\theta^{\varepsilon_2,t}-\theta^{\varepsilon_1,t})\, ds.
\end{align*}
	Denote $w^{\varepsilon_2}(t)-w^{\varepsilon_1}(t) = z(t) $ and $\theta^{\varepsilon_2,t}-\theta^{\varepsilon_1,t}=\zeta^t$. We obtain
	$$
	\frac{d}{dt}|z(t)|^2 \leq C|x^{\varepsilon_2}(t)-x^{\varepsilon_1}(t)|\,|w^{\varepsilon_2}(t)|\,  |z(t)| + C|z(t)|^2+C|\varepsilon_2-\varepsilon_1|(|w^{\varepsilon_2}(t)|+\|\theta^{\varepsilon_2,t}\|)|z(t)|+C\varepsilon_1\|\zeta^t\|\, |z(t)|.
	$$
Using \eqref{var_est1} and Lemma \ref{lem_diff} we obtain
	$$
\frac{d}{dt}|z(t)|^2 \leq Ce^{Ct}|x^{\varepsilon_2}_0-x^{\varepsilon_1}_0|^2+Ce^{Ct}\|\eta^{0,\varepsilon_2}-\eta^{0,\varepsilon_1}\|^2+ C|z(t)|^2+Ce^{Ct}|\varepsilon_2-\varepsilon_1|^2+C\varepsilon_1\|\zeta^t\|^2,
$$
where the constants $C$ depend on the initial data for original problems and variational problems.
		We need to derive the estimate on the difference of the norms  $\|\theta^{\varepsilon_1,t}-\theta^{\varepsilon_2,t}\| = \|\zeta^t\|$. To this end, we use  \eqref{eq:3}, whence
		$$
		\frac{d}{dt}\|\zeta^t\|^2+C_A\|\zeta^t\|^2\leq C\|\zeta^t\|\, |z'(t)|.
		$$
		It follows that
		$$
		\frac{d}{dt}\|\zeta^t\|^2\leq Ce^{Ct}|x^{\varepsilon_2}_0-x^{\varepsilon_1}_0|^2+Ce^{Ct}\|\eta^{0,\varepsilon_2}-\eta^{0,\varepsilon_1}\|^2+ C|z(t)|^2+Ce^{Ct}|\varepsilon_2-\varepsilon_1|^2.
		$$
		We deduce the estimate
			$$
		\frac{d}{dt}(|z(t)|^2+\|\zeta^t\|^2) \leq C(|z(t)|^2+\|\zeta^t\|^2) + Ce^{Ct}(|x^{\varepsilon_2}_0-x^{\varepsilon_1}_0|^2+\|\eta^{0,\varepsilon_2}-\eta^{0,\varepsilon_1}\|^2+ |\varepsilon_2-\varepsilon_1|^2),
		$$
		and the Gronwall lemma yields the desired assertion.
		\end{proof}
		
%\section{Further properties of the weighted history norm.} \label{memo_prop}
%	We show some flexibility in the estimate. Now let $B$ be arbitrary matrix.  Then, the analogous argument leads to the following results
%	\begin{lem}
%		Let $(x,\eta)$ be a solution. Then,
%		\begin{equation}\label{eq:5}
%			\frac{d}{d t}\|B \eta^t\|^2  + C\|B\eta^t\|^2 \leq  - 2\left(\int_0^\infty A(s)B \eta^t(s)\, ds,B x'(t)\right)
%		\end{equation}
%		and
%		\begin{equation}\label{eq:6}
%			\|B \eta^{t_2}\|^2 \leq e^{-C(t_2-t_1)}\|B \eta^{t_1}\|^2 - 2 e^{-Ct_2} \int_{t_1}^{t_2}e^{Ct}\left(\int_0^\infty A(s)B \eta^t(s)\, ds,B x'(t)\right) \, dt\ \ \textrm{for} \ \ t_1< t_2.
%		\end{equation}
%		Moreover if $(x,\eta)$ and $(y,\xi)$ are two solutions.
%		\begin{equation}\label{eq:7}
%			\frac{d}{d t}\|B(\eta^t-\xi^t)\|^2  + C\|B(\eta^t-\xi^t)\|^2 \leq  - 2\left(\int_0^\infty A(s)B(\eta^t(s)-\xi^t(s))\, ds,B (x(t)-y(t))'\right)
%		\end{equation}
%		and
%		\begin{align*}
%			& \|B(\eta^{t_2}-\xi^{t_2})\|^2 \leq e^{-C(t_2-t_1)}\|B(\eta^{t_1}-\xi^{t_1})\|^2 \\
%			& \qquad - 2 e^{-Ct_2} \int_{t_1}^{t_2}e^{Ct}\left(\int_0^\infty A(s)B(\eta^t(s)-\xi^t(s))\, ds,B(x'(t)-y'(t))\right) \, dt\ \ \textrm{for} \ \ t_1< t_2.
%		\end{align*}
%	\end{lem}

\section{Continuation of  isolating blocks with cone conditions.}
\label{sec:LapCC}

The goal of this section is to show that local  dynamical properties of \eqref{eq:ODE} "survive" as we pass to \eqref{eq:DDE}.
The dynamical objects whose "survival" we prove  are
 the isolating blocks satisfying cone conditions that are continued from \eqref{eq:ODE}. These blocks are certain closed "box-like" neighbourhoods of the equilibria for which the isolation property means that the entry and exit behavior of the dynamical system on the boundary is well defined. The precise definition of isolating blocks as well as the argument that their existence guarantees the existence of local stable and unstable manifolds are moved to Appendix 2. We show here that after extending them by Cartesian product with appropriate ball in the memory variable $\eta$ included in the stable variables, the same blocks are valid   for sufficiently small $\varepsilon$ in \eqref{eq:DDE}. The continuation of isolating block is established in Section~\ref{subsec:cont-iso-block}. The cone conditions are discussed in Section~\ref{subsec:coneCond}.

\subsection{Isolating block for $\varepsilon=0$.}
\label{subsec:cont-iso-block}
The next result follows from \cite[Theorem 26]{Zgliczynski}. We give a short proof for the completeness of exposition. In this section we use the notation
$$
B_u(\delta) = \prod_{k=1}^{u_1}[-\delta, \delta] \times \prod_{k=u_1+1}^{u_1+u_2}\{(x,y)\,:\, x^2+y^2 \leq \delta^2\},
$$
and
$$
B_s(\delta) = \prod_{k=u_1+u_2+1}^{u_1+u_2+s_1}[-\delta, \delta]\times \prod_{k=u_1+u_2+s_1+1}^{u_1+u_2+s_1+s_2}\{(x,y)\,:\; x^2+y^2 \leq \delta^2\}.
$$
\begin{lem}\label{isolation:lem5}
Let $x_0$ be such that $f(x_0) = 0$. Assume that this equilibrium is hyperbolic, that is, that the spectrum of the matrix $D f(x_0)$ does not intersect the imaginary axis, with $s$ equal to the dimension of its stable space, and $u=d-s$ the dimension of its unstable space. Let $s=s_1+2s_2$, where $s_1$ is the dimension of the generalized eigenspace related with real stable eigenvalues, and $2s_2$ is the dimension of the generalized eigenspace related with complex stable eigenvalues. Analogously, $u=u_1+2u_2$. For every sufficiently small constant $\kappa>0$ there exists the nonsingular matrix $T_\kappa$ and a number $\delta_0 >0$ such that for every $\delta\in (0,\delta_0)$ the set
 $$
N_\kappa(\delta) = T_\kappa \left(B_u(\delta) \times B_s(\delta)\right)  + x_0.
$$
is an isolating block with cones for $\varepsilon=0$, i.e. for equation \eqref{eq:ODE}.
\end{lem}
\begin{proof}
We begin with the discussion of the role of the constant $\kappa$. We fix $\kappa >0$,  this will be some small number.
	Let $T_\kappa$ be an invertible matrix such that $T_\kappa^{-1} D f(x_0) T_\kappa$ is the Jordan form, that on the diagonal has either real eigenvalues $\lambda$  of $D f(x_0)$, or blocks $\begin{pmatrix}\alpha & \beta\\ -\beta & \alpha \end{pmatrix}$ in case of complex eigenvalues, and all off-diagonal terms have absolute values   not greater than $\kappa$. Assume that the eigenvalues in $T_\kappa$ are sorted such that: first there are real positive eigenvalues, then complex eigenvalues with positive real part, then negative real eigenvalues, and finally complex eigenvalues with negative real part. Moreover assume that
 $T_\kappa^{-1} D f(x_0) T_\kappa = \begin{pmatrix} A& 0\\ 0 & B\end{pmatrix}$, where $A+A^T \in \R^{u\times u}$ is negative definite and $B+B^T \in \R^{s\times s}$ is positive definite.  For every $\kappa > 0$ such change of coordinates $T_\kappa$ exists.  We first prove that the set $N_\kappa(\delta)$
 is an isolating block for sufficiently small $\kappa$.

 If we denote $x = x_0+T_\kappa y$, we obtain the system
 $$
 y' = T_\kappa^{-1}D f(x_0) T_\kappa y + T_\kappa^{-1} f(x_0 + T_\kappa y) - T_\kappa^{-1}D f(x_0) T_\kappa y = h(y).
 $$
 Now for
 $y\in  B_u(2\delta) \times B_s(2\delta),$ we deduce, by the Taylor theorem, as $f\in C^2(\R^d,\R^d)$ that
 $$
 |T_\kappa^{-1} f(x_0 + T_\kappa y) - T_\kappa^{-1}D f(x_0) T_\kappa y | \leq C_\kappa \delta^2
 $$\color{black}
 where $C_\kappa$ depends on $\kappa$ but not on $\delta \in (0,\delta_0)$.
We need to prove that:
\begin{itemize}
	\item if
	 $y\in  B_u(2\delta) \times \partial B_s(\delta),$
  then
	\begin{equation}\label{eq:1is}
	h_i(y) y_i < 0 \ \ \textrm{for} \ \ i\in \{u_1+2u_2+1,\ldots, u_1+2u_2+s_1\}
	\end{equation}
	and
		\begin{equation}\label{eq:2is}
	h_i(y) y_i + h_{i+1}(y)y_{i+1}< 0 \ \ \textrm{for} \ \ i\in \{u_1+2u_2+s_1+1,\ldots,u_1+2u_2+s_1+j,\ldots u_1+2u_2+s_1+(2s_2-1)\},
	\end{equation}
	where $j$ are odd numbers,
	\item if
	$y\in \left(B_u(2\delta) \setminus (\textrm{int} B_u(\delta))\right) \times B_s(\delta),$
	then
		\begin{equation}\label{eq:3is}
	h_i(y)y_i > 0 \ \ \textrm{for} \ \ i\in \{1,\ldots, u_1\}
	\end{equation}
	and
		\begin{equation}\label{eq:4is}
	h_i(y) y_i + h_{i+1}(y)y_{i+1}> 0 \ \ \textrm{for} \ \ i\in \{ u_1+1,\ldots,u_1+j,\ldots u_1+2u_2-1\},
	\end{equation}
	where $j$ are odd numbers,
\end{itemize}
By the Lipschitz continuous dependence on the initial condition, on bounded sets of initial data and compact time intervals these conditions imply the isolation given in Definition \ref{isolation}.

To prove the first assertion observe that for $i\in u_1+u_2+1,\ldots, u_1+u_2+s_1$.
$$
h_i(y)y_i = \lambda_i y_i^2 + G(y),
$$
where
$$
|G(y)| \leq (d-1)\kappa 4 |\delta|^2 + 2C_\kappa
|\delta|^3,
$$
the first term coming from off diagonal values (at most $d-1$) in $T_\kappa^{-1}D f(x_0) T_\kappa$, and the second one from the remainder which is a product of number which is dominated by the euclidean norm of a vector bounded by $C_\kappa |\delta|^2$ and a number bounded by $2\delta$. This means that we can choose $\kappa$ small enough (related to the lowest eigenvalue $\lambda_i$) and $\delta_0$ (that is chosen according to $C_\kappa$) and  get \eqref{eq:1is}.

For the complex pairs of eigenvalues, the off diagonal terms in blocks $\begin{pmatrix}
	\alpha & \beta \\ -\beta & \alpha
\end{pmatrix}$ cancel and we obtain
$$
h_i(y)y_i + h_{i+1}(y)y_{i+1} = \textrm{Re}\, \lambda_i (y_i^2+y_{i+1}^2) + G(y),
$$
with
$$
|G(y)| \leq 8 (d-2)\kappa |\delta|^2 + 2C_k|\delta|^3,
$$
and \eqref{eq:2is} holds analogously as \eqref{eq:1is}. Verification of \eqref{eq:3is} and \eqref{eq:4is} follows analogously.

To see that the cone condition holds it is enough to take the matrix $Q$ such that $q_{ij} = 0$ for $i\neq j$, $q_{ii} = -1$ for $i=1,\ldots, u$ and $q_{ii}=1$ for $i=u+1,\ldots, d$ and see that $QT_\kappa^{-1}D f(x_0) T_\kappa+T_\kappa^{-1}D f(x_0) T_\kappa Q$ is positive definite, which must be preserved on a small neighborhood of $y=0$.
\end{proof}

\subsection{Continuation of isolation property for $\varepsilon>0$}\label{sec:isoeps}
In the subsequent part of this section we will show that it is possible to choose $\delta$ and $\kappa$ as well as $R>0$ such that if $N_\kappa(\delta)$ is an isolating block with cones for \eqref{eq:ODE} then the set $\overline{B}_{L_A^2(\R^+)^d}(0,R)\times N_\kappa(\delta)$ is an isolating block with cones for \eqref{eq:DDE}. We start from an estimate. Substitute \eqref{eq:main} in \eqref{eq:1}.
Then we obtain
\begin{align*}
& \frac{d}{d t}\|\eta^t\|^2  + C_A\|\eta^t\|^2 \leq  - 2\left(\int_0^\infty A(s)\eta^t(s)\, ds,f(x(t))\right) -2  \varepsilon\left(\int_0^\infty A(s)\eta^t(s)\, ds,\left( \int_{0}^\infty M(s)\, ds \right) x(t)\right)\\
& \ \ \  - 2\varepsilon \left(\int_0^\infty A(s)\eta^t(s)\, ds,\int_{0}^\infty M(s) \eta^t(s)\, ds \right) .
\end{align*}
After computations which use Lemmas \ref{lem:a_prop} and \ref{lem:M} it follows that
\begin{align*}
	& \frac{d}{d t}\|\eta^t\|^2  + C_A\|\eta^t\|^2 \leq  2 D_A \|\eta^t\| |f(x(t))| + 2  \varepsilon D_A\| \eta^t\| \int_{0}^\infty \|M(s)\|\, ds |x(t)| + 2\varepsilon D_AD_M \|\eta^t\|^2.
\end{align*}

\begin{align*}
	& \frac{d}{d t}\|\eta^t\|^2  \leq \|\eta^t\| \left(  2 D_A  |f(x(t))| + 2  \varepsilon D_A \int_{0}^\infty \|M(s)\|\, ds |x(t)| + 2\varepsilon D_A D_M \|\eta^t\|- C_A\|\eta^t\| \right).
\end{align*}

The above computation leads is a straightforward way to the following lemma.
\begin{lem}
\label{lem:iso-tail-pz}

Suppose that
$f(x_0) = 0$ and that $N \subset \mathbb{R}^d$ is a compact set containing $x_0$. Moreover let
 $\varepsilon < \frac{C_A}{2DAD_M}$ and
$$
R > \frac{2 D_A  \left(\sup_{z\in N}|f(z)|+  \varepsilon  \int_{0}^\infty \|M(s)\|\, ds  \cdot \sup_{z\in N}|z|  \right) }{C_A-2\varepsilon D_AD_M}.
$$
Then for $\eta \in \partial B_{L^2_A(\R^+)^d}(0,R)$, and $y\in N$ there holds
$$\frac{d}{dt} \|\eta^t\|^2_{L^2_A(\R^+)^d} < 0.
$$
\end{lem}

In the next theorem we demonstrate that the entry and exit behavior established in Lemma \ref{isolation:lem5} are preserved for $\varepsilon>0$. The new isolating block will have the form $$(0,x_0)+\overline{B}_{L^2_A(\R^+)^d}(0,R)\times T_\kappa(B_u(\delta)\times B_s(\delta)).$$
The entry set will be given by
$$(0,x_0)+\partial {B}_{L^2_A(\R^+)^d}(0,R)\times T_\kappa(B_u(\delta)\times \partial B_s(\delta))$$
and the exit set by
$$(0,x_0)+\overline{B}_{L^2_A(\R^+)^d}(0,R)\times T_\kappa(\partial B_u(\delta)\times B_s(\delta)).$$
To this end let us first rewrite the equation \eqref{eq:DDE} in the changed variables $y$.
$$
y'(t) = h(y(t)) + \varepsilon T_{\kappa}^{-1} \left(\int_0^\infty M(s)\, ds\right) (x_0+T_\kappa y(t)) + \varepsilon T_\kappa^{-1}\int_0^\infty M(s) \eta^t(s)\, ds.
$$

Now choose $\kappa$ and $\delta_0$ such that Lemma \ref{isolation:lem5} holds and assume that $\delta < \frac{\delta_0}{2}$. For such $\delta$ let $r(\delta)$ be a smallest possible number such that $N_\kappa(2\delta) \subset B(x_0,r)$. Note that $r\to 0$ as $\delta \to 0$. Take $\eta^t \in \overline{B}_{L^2_A(\R^+)^n}(0,R)$ and $x=T_\kappa y+x_0 \in N_\kappa(2\delta)$. We rewrite the $i$-th equation of the above system as
$$
y_i'(t) = h_i(y(t)) + g_i(\eta^t,y(t)),
$$
where
$$
g_i(\eta^t,y(t)) = \varepsilon \left(T_{\kappa}^{-1} \left(\int_0^\infty M(s)\, ds\right) (x_0+ T_\kappa y(t))\right)_i + \varepsilon \left(T_\kappa^{-1}\int_0^\infty M(s) \eta^t(s)\, ds\right)_i
$$
hence
\begin{equation}\label{eq:gi}
|g_i(\eta,y)|\leq |g(\eta,y)| \leq \varepsilon \|T_\kappa^{-1}\|D_M\|\eta^t\| + \varepsilon  \|T_\kappa^{-1}\|\int_0^\infty \|M(s)\|\, ds (|x_0| + r) \leq  \varepsilon C_1 (R+|x_0|+r),
\end{equation}
for a constant $C_1>0$. Now we are ready to state our result about the continuation of the isolating block for $\varepsilon >0$.
\begin{thm}\label{thm:isol}
There exists $\kappa>0$, $\varepsilon_0 >0$, $\delta>0$, and $R>0$ such that for every  fixed point $x_0$, every $\varepsilon\in (0,\varepsilon_0)$ the set $\overline{B}_{L^2_A(\R^+)^d}(0,R)\times N_\kappa(\delta)$ is an isolating block for \eqref{eq:DDE}, i.e.
\begin{itemize}
	\item[(I)] we have the entry behavior on $\partial B_s(\delta)$, that is, if
	$(\eta,y)\in  \overline{B}_{L^2_A(\mathbb{R}_+)^d}(0,R) \times B_u(2\delta) \times \partial B_s(\delta),$
	then
	\begin{equation}\label{eq:1is_e}
		(h_i(y)+g_i(\eta,y)) y_i < 0 \ \ \textrm{for} \ \ i\in \{u_1+2u_2+1,\ldots, u_1+2u_2+s_1\},
	\end{equation}
	and
	\begin{align}\label{eq:2is_e}
		& (h_i(y)+g_i(\eta,y)) y_i + (h_{i+1}(y)+g_{i+1}(\eta,y))y_{i+1}< 0 \\
		&\nonumber \ \ \textrm{for} \ \ i\in \{u_1+2u_2+s_1+1,\ldots,u_1+2u_2+s_1+j,\ldots, u_1+2u_2+s_1+(2s_2-1)\},
	\end{align}
	where $j$ are odd numbers,
	\item[(II)] we have the exit behavior on $\partial B_u(\delta)$, that is, if
	$(\eta,y)\in  \overline{B}_{L^2_A(\mathbb{R}_+)^d}(0,R)\times \left(B_u(2\delta) \setminus (\textrm{int} B_u(\delta))\right) \times B_s(\delta),$
	then
	\begin{equation}\label{eq:3is_e}
		(h_i(y)+g_i(\eta,y))y_i > 0 \ \ \textrm{for} \ \ i\in \{1,\ldots, u_1\},
	\end{equation}
	and
	\begin{equation}\label{eq:4is_e}
		(h_i(y)+g_i(\eta,y)) y_i + (h_{i+1}(y)+g_{i+1}(\eta,y))y_{i+1}> 0 \ \ \textrm{for} \ \ i\in \{ u_1+1,\ldots,u_1+j,\ldots u_1+2u_2-1\},
	\end{equation}
	where $j$ are odd numbers,
	\item[(III)] we have the entry behavior on $\partial B_{L^2_A(\mathbb{R}_+)^d}(0,R)$, that is, if $(\eta,y) \in  \partial B_{L^2_A(\mathbb{R}_+)^d}(0,R) \times B_u(2\delta)\times B_s(\delta) $, then
	$$
	\frac{d}{dt}\|\eta_t\|^2_{L^2_A(\mathbb{R}_+)^d} < 0 \ \ \textrm{at}\ \ t=0. \,
	$$
\end{itemize}
\end{thm}
\begin{proof}
	We first provide the condition needed for  (III) to hold: this is the entry condition for the variable $\eta$. Following Lemma \ref{lem:iso-tail-pz} we need that $R>0$ and $\delta_0>0$ should satisfy
			\begin{equation}\label{eq:con1}
	R > \frac{2 D_A  \left(\sup_{z\in N_\kappa(2\delta)}\|D f(z)\|r(\delta) +  \varepsilon  \int_0^\infty\|M(s)\|\,  ds  r(\delta)  \right) }{C_A-2\varepsilon D_AD_M}.
	\end{equation}
	We switch to the conditions needed for (I) and (II), that is, for \eqref{eq:1is_e}--\eqref{eq:4is_e}.

	We first note  that there exists $\delta_0$ and a constant $C_2 > 0$ such that if only $\delta\in (0,\delta_0)$ and $y\in  B_u(2\delta)\times \partial B_s(\delta) \cup (B_u(2\delta)\setminus \textrm{int} B_u(\delta)) \times B_s(\delta)$ then
	\begin{equation}\label{eq:real}
	C_2 \delta \leq |h_i(y)|\ \ \textrm{for indexes corresponding to real eigenvalues},
	\end{equation}
	and
	\begin{equation}\label{eq:complex}
	C_2 \delta^2 \leq |h_i(y)y_i + h_{i+1}(y)y_{i+1}|\ \ \textrm{for indexes corresponding to complex eigenvalues}.
	\end{equation}
	In order to guarantee \eqref{eq:1is_e} and \eqref{eq:3is_e} we need $h_i(y)$ to have the same sign as and $h_i(y)+g_i(\eta,y)$ for $y\in  B_u(2\delta)\times \partial B_s(\delta) \cup (B_u(2\delta)\setminus \textrm{int} B_u(\delta)) \times B_s(\delta)$ and $\eta\in \overline{B}_{L^2_A(\mathbb{R}_+)^d}(0,R)$.
	For the complex eigenvalues we need, on the other hand,
	that $(h_i(y)+g_i(\eta,y)) y_i + (h_{i+1}(y)+g_{i+1}(\eta,y))y_{i+1}$ and $h_i(y)y_i + h_{i+1}(y)y_{i+1}$ have the same signs. Therefore, in view of
	 \eqref{eq:real} it is sufficient to prove that $|g_i(\eta,y)| < C_2 \delta$ and $|g_i(\eta,y)y_i+g_{i+1}(\eta,y)y_{i+1}| < C_2\delta^2$. Using \eqref{eq:gi} it is enough if the following inequality holds
	 \begin{equation}\label{eq:con2}
	 \varepsilon C_1(R+|x_0|+r(\delta)) < C_2\delta.
	 \end{equation}

First choose $\varepsilon_0$ such that $C_A-2\varepsilon_0 D_AD_M > \frac{C_A}{2}$. We need to guarantee \eqref{eq:con1}. Note that if it holds for $\varepsilon=\varepsilon_0$ then the same inequality holds for $\varepsilon \in (0,\varepsilon_0)$.
Now pick $R>0$. We need to choose $\delta$ small enough such that
$$\frac{4 D_A  \sup_{z\in N_\kappa(2\delta)}\|D f(z)\|r(\delta) }{C_A} < \frac{R}{2}.
$$ This is possible, because by decreasing $\delta$ we can make $r(\delta)$ arbitrarily small. Now it is possible to choose sufficiently small $\varepsilon_0$ such that $\varepsilon_0 C_1 (R+|x_0|+r(\delta)) < C_2\delta$ and $\frac{\varepsilon_0 4D_A \int_0^\infty\|M(s)\|\, ds\, r(\delta)}{C_A} < \frac{R}{2}$. Thus both inequalities are satisfied and the proof is complete.
\end{proof}

\subsection{Continuation of cone condition for $\varepsilon>0$.}
\label{subsec:coneCond}
The goal of this section is to show that  cone conditions from (\ref{eq:ODE}) "survive" for sufficiently small $\varepsilon$ for (\ref{eq:DDE}) on the same isolating block. We refer to Appendix 2 for the definition of the cone condition. We assume that it holds for  the equation without memory, that is
$$
x'(t) = f^\varepsilon (x(t)) = f(x(t)) + \varepsilon \int_0^\infty M(s)\, ds x(t).
$$
For this equation we have a quadratic form $Q$  (a symmetric matrix) and a set $N$ (isolating block) such that
on $N$ we have for any $x \in \mathbb{R}^d $ and $|\varepsilon| \leq \Delta$ for some $G>0$
\begin{equation}
	x^t(Df^\varepsilon(N)^TQ + Q Df^\varepsilon(N))x \geq  G |x|^2.
\end{equation}
Note that
$$
Df^\varepsilon(x) = Df(x) + \varepsilon \int_0^\infty M(s)\, ds.
$$
 Hence, if $\varepsilon>0$ is small enough, then the same form $Q$ is valid both for $f$ and for $f^\varepsilon$ (with possibly smaller constant $G$).
Let $(\eta,x)$ and $(\xi,y)$ be two solutions of \eqref{eq:main}--\eqref{main2} such $x,y \in N$ and $\|\eta\|, \|\xi\|\leq R$.
Let $E>0$ be any positive constant. We hope that for the quadratic form
\begin{equation}
	\widetilde{Q}(\eta,x)=Q(x) -  E \|\eta\|^2  \label{eq:tildeQ-new}
\end{equation}
we will have cone-conditions on the set (isolated block)
\begin{equation}
	\widetilde{N}= \{\|\eta \| \leq R\}\times N.
\end{equation}
We have
\begin{eqnarray*}
	\frac{d}{dt} \left( Q(x(t)-y(t),x(t)-y(t))  \right) =
	(x'-y')^\top Q (x-y) +  (x-y)^\top Q (x'-y').
\end{eqnarray*}
Since
\begin{eqnarray*}
	x'-y'&=&\left(f^\varepsilon(x) + \varepsilon \int_0^\infty M(s)\eta(s)ds \right)
	-\left(f^\varepsilon(y)+ \varepsilon \int_0^\infty M(s)\xi(s)ds\right) \\
	&=&(f^\varepsilon(x) - f^\varepsilon(y))  + \varepsilon \int_0^\infty M(s) (\eta(s) - \xi(s))ds \\
	&=&  \overline{Df^\varepsilon [x,y]}(x - y)  + \varepsilon \int_0^\infty M(s) (\eta(s) - \xi(s))ds,
\end{eqnarray*}
we obtain using Lemma \ref{lem:M}
\begin{align*}
	&(x'-y')^\top Q (x-y) +  (x-y)^\top Q (x'-y') = (x-y)^\top\left( \overline{Df^\varepsilon[x,y]}^\top Q + Q \overline{Df^\varepsilon[x,y]} \right)(x-y) \\
&\ \ \ \ 	+ \varepsilon \left( \int_0^\infty M(s) (\eta(s) - \xi(s))ds \right)^\top Q (x-y) + \varepsilon (x-y)^T Q \int_0^\infty M(s) (\eta(s) - \xi(s))ds    \\
&\ 	\geq G |x-y|^2 -2 \varepsilon D_M \|Q\| \cdot |x-y| \cdot \|\eta - \xi\|.
\end{align*}
From Lemma~\ref{lem:diff_ogon} we have
\begin{align*}
	&\frac{d}{dt} \|\eta^t - \xi^t\|^2 \leq - C_A \|\eta - \xi\|^2 - 2\left(\int_0^\infty A(s) (\eta(s) - \xi(s))ds, (x(t)-y(t))' \right) \\
&\ \ \ 	\leq - C_A \|\eta - \xi\|^2  - 2 \left(\int_0^\infty A(s) (\eta(s) - \xi(s))ds, \overline{Df^\varepsilon[x,y]}(x-y) \right) + \\
&\ \ \ \ \ \ 	-2\varepsilon \left(\int_0^\infty A(s) (\eta(s) - \xi(s))ds, \int_0^\infty M(s) (\eta(s) - \xi(s))ds \right) \\
&\ \  \ 	\leq - C_A \|\eta - \xi\|^2  + 2 D_A \|\eta - \xi\| \cdot \|Df^\varepsilon (N)\| \cdot |x-y | + \varepsilon 2 D_A D_M \|\eta - \xi\|^2
\end{align*}
Now we are ready to demonstrate that the cone condition holds. From previous derivations we obtain

\begin{align*}
&	\frac{d}{dt} \left( Q(x(t)-y(t),x(t)-y(t)) - E  \|\eta^t - \xi^t\|^2 \right) \\
&\ \ \ \geq 	G |x-y|^2 -2\varepsilon D_M \|Q\| \cdot |x-y| \cdot \|\eta - \xi\|   \\
&\ \ \ \ \ \ 	-  E \left(- C_A \|\eta - \xi\|^2  + 2 D_A\|\eta - \xi\| \cdot \|Df^\varepsilon (N)\| \cdot |x-y | + \varepsilon 2 D_M D_A \|\eta - \xi\|^2 \right) \\
&\ \ \ 	= G |x-y|^2  +  2\left(- \varepsilon D_M \|Q\|  -  E D_A \|Df^\varepsilon (N)\| \right) |x-y| \cdot \|\eta - \xi\| \\
&\ \ \ \ \ \ 	+ \left(C_AE - 2\varepsilon D_AD_M \right) \|\eta - \xi\|^2
\end{align*}
The expression on the right-hand side is a quadratic form in terms of $(|x-y|,\|\eta-\xi\|)$  with the matrix
\begin{equation}
	B= \left[
	\begin{array}{cc}
		G &   \left(- \varepsilon D_M \|Q\|  -  E D_A \|Df^\varepsilon (N)\| \right) \\
		\left(- \varepsilon D_M \|Q\|  -  E D_A \|Df^\varepsilon (N)\| \right) &  \left(C_AE - 2\varepsilon D_A D_M  \right) \\
	\end{array}
	\right]
\end{equation}
Consider first the case with $\varepsilon=0$. Matrix $B$ becomes
\begin{equation}
	B_0= \left[
	\begin{array}{cc}
		G &   -  E D_A \|Df (N)\|  \\
		-  E D_A \|Df (N)\| &  C_AE  \\
	\end{array}
	\right]
\end{equation}
It is positive definite provided the following condition holds
\begin{equation}
	0< \det(B_0) = C_AE G - E^2 D_A^2  \|Df (N)\|^2,
\end{equation}
which is satisfied if
\begin{equation}
	E < \frac{C_A G}{D_A^2  \|Df (N)\|^2}.
\end{equation}
Since $\det B$ depends continuously on $\varepsilon$ we obtain the following theorem.

\begin{lem}\label{lem:E1}
	For any $E>0$ satisfying
	\begin{equation}
		E < \frac{C_A G}{D_A^2  \|Df (N)\|^2}.
	\end{equation}
	there exists $\varepsilon_1=\varepsilon_1(E)$, such that cone condition holds for quadratic form (\ref{eq:tildeQ-new}) for any $\varepsilon \leq \varepsilon_1$.
\end{lem}

The next result will be needed for the continuous dependence of local stable and unstable manifolds on $\varepsilon$. We obtain this continuity using the abstract results from Appendix 3. To use them, we must verify  the cone conditions with parameter given in Definition \ref{df:hconespar}. Consider the problem governed by \eqref{eq:main}--\eqref{main2} and
denote the solution with $\varepsilon \in [0,\Delta]$ and initial data $ (\eta^0,x_0) \in \widetilde{N}$ by $S^\varepsilon(t)(\eta^0,x_0)=(\eta^t_{\varepsilon,x_0,\eta_0},x_{\varepsilon,x_0,\eta_0})$.
Define
\begin{align*}
& \widehat{Q}(x,\eta,\varepsilon) = \widetilde{Q}(x,\eta) + L|\varepsilon|^2 = {Q}(x) - E\|\eta\|^2 + L|\varepsilon|^2,
\end{align*}
where $L$ can be either positive or negative constant. In the next lemma we prove that it is possible to choose this constant in such a way, that $\widehat{Q}$ satisfies the conditions given in Definition \ref{df:hconespar}. If $L$ is a large positive number then we obtain continuous dependence of unstable manifold on $\varepsilon$ and if $L$ is a large negative number then we obtain the continuous dependence of stable manifold on $\varepsilon$.
\begin{lem}\label{lem:cone:par}
	There exists $L_0>0$ and $E_{max}>0$ such that for every $|L|\geq L_0$ and $E\in (0,E_{max})$ there exists $\Delta(E) > 0$ such that the cone conditions with parameter given in (i) and (ii) of Definition \ref{df:hconespar} are satisfied on the h-set $\widetilde{N}$ with $\varepsilon_1,\varepsilon_2 \in [0,\Delta(E)]$
\end{lem}
\begin{proof}
Assume that $(\eta^0_1,x^0_1,\varepsilon_1)$ and $(\eta^0_2,x^0_2,\varepsilon_2)$ are such that
$$
\widehat{Q}(x^0_1-x^0_2,\eta^0_1-\eta^0_2,\varepsilon_1-\varepsilon_2) = 0.
$$
We must prove that
$$
\frac{d}{dt}\widehat{Q}(x_{\varepsilon_1,x^0_1\eta^0_1}(t)-x_{\varepsilon_2,x^0_2\eta^0_2}(t),\eta_{\varepsilon_1,x_1^0,\eta_1^0}(t)-\eta_{\varepsilon_2,x^0_2\eta^0_2}(t),\varepsilon_1-\varepsilon_2) \geq 0\ \ \textrm{for}\ \ t=0.
$$
Denoting, for simplicity, $(\eta_1(t),x_1(t)) = (\eta_{\varepsilon_1,x_1^0,\eta_1^0}(t),x_{\varepsilon_1,x^0_1\eta^0_1}(t))$ and \\
 $(\eta_2(t),x_2(t)) = (\eta_{\varepsilon_2,x_2^0,\eta_2^0}(t),x_{\varepsilon_2,x^0_2\eta^0_2}(t))$ we should prove that
\begin{equation}\label{extended_cone}
\frac{d}{dt}\left((x_1(t)-x_2(t))^\top Q (x_1(t)-x_2(t)) - E\|\eta_1(t)-\eta_2(t)\|^2\right) \geq 0\ \ \textrm{for}\ \ t=0.
\end{equation}
We estimate both terms from below separately
\begin{align*}
& \frac{d}{dt}((x_1(t)-x_2(t))^\top Q (x_1(t)-x_2(t))
 = (x_1'(t)-x_2'(t))^\top Q (x_1(t)-x_2(t)) + (x_1(t)-x_2(t))^\top Q (x_1'(t)-x_2'(t))\\
 & \ \ \ = (f(x_1(t))-f(x_2(t)))^\top Q (x_1(t)-x_2(t)) +  (\varepsilon_1 x_1(t) - \varepsilon_2 x_2(t))^\top\int_0^\infty M^\top(s)\, dsQ (x_1(t)-x_2(t))\\
 & \ \ \ \qquad  + \left( \varepsilon_1 \int_0^\infty M(s) \eta_1^t(s)\, ds - \varepsilon_2 \int_0^\infty M(s) \eta_2^t(s) \ ds\right)^\top Q (x_1(t)-x_2(t))\\
 & \ \ \ \qquad +  (x_1(t)-x_2(t))^\top Q (f(x_1(t))-f(x_2(t))) +  (x_1(t)-x_2(t))^\top Q \int_0^\infty M(s)\, ds (\varepsilon_1 x_1(t) - \varepsilon_2 x_2(t)) \\
 & \ \ \ \qquad + (x_1(t)-x_2(t))^\top Q \left( \varepsilon_1 \int_0^\infty M(s) \eta_1^t(s)\, ds - \varepsilon_2 \int_0^\infty M(s) \eta_2^t(s) \ ds \right)\\
  & \ \ \ = (x_1(t)-x_2(t))^\top (Df(N)^\top Q + QDf(N)) (x_1(t)-x_2(t)) \\
  & \ \ \ \qquad +   \varepsilon_2 ( x_1(t) -  x_2(t))^\top\left(\int_0^\infty M^\top(s)\, dsQ+  Q \int_0^\infty M(s)\, ds\right) (x_1(t)-x_2(t))\\
 & \ \ \ \qquad + 2\varepsilon_2(x_1(t)-x_2(t))^\top Q \left(   \int_0^\infty M(s) (\eta_1^t(s) -  \eta_2^t(s))\, ds \right)\\
 & \ \ \ \qquad + 2(\varepsilon_1-\varepsilon_2) (x_1(t)-x_2(t))^\top Q \int_0^\infty M(s)\, ds  x_1(t) + 2(\varepsilon_1-\varepsilon_2)(x_1(t)-x_2(t))^\top Q \left( \int_0^\infty M(s) \eta_1^t(s) \, ds \right) \\
 & \ \ \ = I_1 + I_2 + I_3 + I_4.
\end{align*}
Now
$$
I_1 + I_2 \geq G |x_1(t)-x_2(t)|^2,
$$
where $G$ can be chosen uniformly for $\varepsilon\in [0,\Delta].
$
Moreover
$$
I_3 \geq - 2 \Delta |x_1(t)-x_2(t)|\  \|Q\| D_M \|\eta_1^t - \eta_2^t\|,
$$
and
$$
I_4 \geq - 2 |\varepsilon_1-\varepsilon_2| \ |x_1(t)-x_2(t)|\  \|Q\| \left(\int_0^\infty \|M(s)\|\, ds \sup_{x\in N}|x| + D_M R\right).
$$
For simplicity we use the following notation for the constant which will appear several times in the subsequent computations  $\overline{R} = \int_0^\infty \|M(s)\|\, ds \sup_{x\in N}|x| + D_M R $.
Summarizing, we obtain
\begin{align*}
& \frac{d}{dt}((x_1(t)-x_2(t))^\top Q (x_1(t)-x_2(t)) \geq G |x_1(t)-x_2(t)|^2 - 2 \Delta |x_1(t)-x_2(t)|\, \|Q\| D_M \|\eta_1^t - \eta_2^t\| \\
& \ \ \ - 2 |\varepsilon_1-\varepsilon_2|\, |x_1(t)-x_2(t)|\, \|Q\| \overline{R}.
\end{align*}
We estimate the second term in \eqref{extended_cone} from Lemma~\ref{lem:diff_ogon}
\begin{align*}
&\frac{d}{dt} \|\eta_1^t - \eta_2^t\|^2 \leq - C_A \|\eta_1^t - \eta_2^t\|^2 - 2\left(\int_0^\infty A(s) (\eta_1^t(s) - \eta_2^t(s))ds, (x_1(t)-x_2(t))' \right) \\
&\ \ \ 	= - C_A \|\eta_1^t- \eta_2^t\|^2  - 2 \left(\int_0^\infty A(s) (\eta_1^t(s) - \eta_2^t(s))ds, f(x_1(t))-f(x_2(t))\right)  \\
& \
 \ \ \ - 2 \left(\int_0^\infty A(s) (\eta_1^t(s) - \eta_2^t(s))ds, \varepsilon_1 \int_0^\infty M(s)\, ds x_1(t)-\varepsilon_2 \int_0^\infty M(s)\, ds  x_2(t)\right)  \\
&\ \ \ \ \ \ 	-2 \left(\int_0^\infty A(s) (\eta_1^t(s) - \eta_2^t(s))ds, \varepsilon_1\int_0^\infty M(s) \eta_1^t(s) \, ds - \varepsilon_2\int_0^\infty M(s) \eta_2^t(s)ds \right) \\
&\ \ \ 	= - C_A \|\eta_1^t- \eta_2^t\|^2  - 2 \left(\int_0^\infty A(s) (\eta_1^t(s) - \eta_2^t(s))ds, f(x_1(t))-f(x_2(t))\right)  \\
& \
\ \ \ - 2 (\varepsilon_1-\varepsilon_2)\left(\int_0^\infty A(s) (\eta_1^t(s) - \eta_2^t(s))ds,  \int_0^\infty M(s)\, ds x_1(t)+ \int_0^\infty M(s) \eta_1^t(s) \, ds  \right)  \\
& \
\ \ \ - 2 \varepsilon_2 \left(\int_0^\infty A(s) (\eta_1^t(s) - \eta_2^t(s))ds,    \int_0^\infty M(s)\, ds (x_1(t) -  x_2(t))+\int_0^\infty M(s) (\eta_1^t(s)-\eta_2^t(s)) \, ds \right).
\end{align*}
It follows that
\begin{align*}&\frac{d}{dt} \|\eta_1^t - \eta_2^t\|^2 \leq - C_A \|\eta^t_1 - \eta^t_2\|^2  + 2 D_A \|\eta_1^t - \eta_2^t\| \cdot \|Df (N)\| \cdot |x_1(t)-x_2(t) | \\
& \ \ \ +2 |\varepsilon_1-\varepsilon_2| D_A\|\eta_1^t-\eta_2^t\| \overline{R}+ 2\Delta  D_M D_A \|\eta_1^t - \eta_2^t\|^2 + 2\Delta  D_A \int_0^\infty\|M(s)\|\, ds \|\eta_1^t - \eta_2^t\|\  |x_1(t)-x_2(t)|.
\end{align*}
Putting together the two estimates we obtain
\begin{align*}
&\frac{d}{dt}\left((x_1(t)-x_2(t))^\top Q (x_1(t)-x_2(t)) - E\|\eta_1(t)-\eta_2(t)\|^2\right) \geq \\
& \ \ \ \geq G |x_1(t)-x_2(t)|^2 + E(C_A- 2\Delta  D_MD_A) \|\eta^t_1 - \eta^t_2\|^2 \\
& \ \ \  - 2  |x_1(t)-x_2(t)| \left(\Delta \|Q\| D_M + ED_A \|Df (N)\|+E\Delta  D_A\int_0^\infty\|M(s)\|\, ds\right) \|\eta_1^t - \eta_2^t\| \\
& \ \ \ - 2 |\varepsilon_1-\varepsilon_2|\, |x_1(t)-x_2(t)|\, \|Q\| \overline{R}   -2 E |\varepsilon_1-\varepsilon_2| D_A\|\eta_1^t-\eta_2^t\| \overline{R}.
\end{align*}
We need the right-hand side of the last estimate to be bounded from below by $0$ at $t=0$, on the boundary of the cone, i.e. for ${Q}(x^0_1-x^0_2) - E\|\eta^0_1-\eta_2^0\|^2 + L|\varepsilon_1-\varepsilon_2|^2 = 0$, whereas we can estimate from above as follows
$$
|\varepsilon_1-\varepsilon_2| \leq \frac{1}{\sqrt{|L|}}\left(\sqrt{\|Q\|} \cdot |x^0_1-x^0_2| + \sqrt{E}\|\eta^0_1-\eta_2^0\|\right)
$$
We deduce that, at $t=0$ we have
\begin{align*}
&\frac{d}{dt}\left((x_1(t)-x_2(t))^\top Q (x_1(t)-x_2(t)) - E\|\eta_1(t)-\eta_2(t)\|^2\right)|_{t=0}\geq  \left(G - \frac{2\|Q\|^\frac{3}{2}\overline{R}}{\sqrt{|L|}} \right) |x_1^0-x_2^0|^2 \\
& \ \ \ \ + E\left(C_A- 2\Delta  D_AD_M-\frac{2 \sqrt{E} D_A \overline{R}}{\sqrt{|L|}} \right) \|\eta^t_1 - \eta^t_2\|^2 \\
& \ \ \  - 2  |x_1^0-x_2^0| \left(\Delta \|Q\| D_M + ED_A \|Df (N)\|+E\Delta  D_A\int_0^\infty\|M(s)\|\, ds+\frac{E \sqrt{\|Q\|} D_A  \overline{R}}{\sqrt{|L|}} + \frac{\|Q\| \overline{R} \sqrt{E} }{\sqrt{|L|}}\right) \|\eta_1^t - \eta_2^t\|.
\end{align*}
We are free to choose sufficiently large (positive or negative) $L$, sufficiently small $\Delta$ and sufficiently small $E$. We already have the upper bound on $E$ in Lemma  \ref{lem:E1} given by $E\leq E_{max}$.

Now suppose that $|L|$ is large enough and $\Delta$ is small enough such that
$$
\sqrt{|L|} \geq \max\left\{ \frac{4 \|Q\|^{\frac{3}{2}}\overline{R}}{G},\frac{8\sqrt{E_{max}}D_A\overline{R}}{C_A}  \right\} \qquad \textrm{and}\qquad \Delta \leq \frac{C_A}{8 D_A D_M}.
$$

With these assumption the above estimate takes the form
\begin{align*}
&\frac{d}{dt}\left((x_1(t)-x_2(t))^\top Q (x_1(t)-x_2(t)) - E\|\eta_1(t)-\eta_2(t)\|^2\right)|_{t=0}\geq \frac{G}{2}  |x_1^0-x_2^0|^2 + \frac{EC_A}{2} \|\eta^t_1 - \eta^t_2\|^2 \\
& \ \ \ \   - 2  |x_1^0-x_2^0| \left(\Delta \|Q\| D_M + ED_A \|Df (N)\|+E\Delta  D_A\int_0^\infty\|M(s)\|\, ds+\frac{E \sqrt{\|Q\|} D_A  \overline{R}}{\sqrt{|L|}} + \frac{\|Q\| \overline{R} \sqrt{E} }{\sqrt{|L|}}\right) \|\eta_1^t - \eta_2^t\|.
\end{align*}
The above quadratic form on $|x_1^0-x_2^0|$ and $\|\eta^t_1 - \eta^t_2\|$ is nonnegatively defined provided
$$
\frac{EC_AG}{4} \geq \left(\Delta \|Q\| D_M + ED_A \|Df (N)\|+E\Delta  D_A\int_0^\infty\|M(s)\|\, ds+\frac{E \sqrt{\|Q\|} D_A  \overline{R}}{\sqrt{|L|}} + \frac{\|Q\| \overline{R} \sqrt{E} }{\sqrt{|L|}}\right)^2.
$$
But we know that
\begin{align*}
&
\left(\Delta \|Q\| D_M + ED_A \|Df (N)\|+E\Delta  D_A\int_0^\infty\|M(s)\|\, ds+\frac{E \sqrt{\|Q\|} D_A  \overline{R}}{\sqrt{|L|}} + \frac{\|Q\| \overline{R} \sqrt{E} }{\sqrt{|L|}}\right)^2 \\
& \ \ \leq 5 \Delta^2 \|Q\|^2 D_M^2 + 5E^2D_A^2 \|Df (N)\|^2+5E^2\Delta^2  D_A^2\left(\int_0^\infty\|M(s)\|\, ds\right)^2+\frac{5 E^2 \|Q\| D_A^2  \overline{R}^2}{|L|} + \frac{5\|Q\|^2 \overline{R}^2 {E} }{|L|}.
\end{align*}
Hence we need the following five inequalities
\begin{align*}
& 5 \Delta^2 \|Q\|^2 D_M^2 \leq \frac{EC_AG}{20},\ \  5E^2D_A^2 \|Df (N)\|^2 \leq \frac{EC_AG}{20},\ \  5E^2\Delta^2  D_A^2\left(\int_0^\infty\|M(s)\|\, ds\right)^2 \leq \frac{EC_AG}{20},\\
& \ \ \  \frac{5 E^2 \|Q\| D_A^2  \overline{R}^2}{|L|} \leq \frac{EC_AG}{20},\ \  \frac{5\|Q\|^2 \overline{R}^2 {E} }{|L|}  \leq \frac{EC_AG}{20},
\end{align*}
or, after the simplification,
\begin{align*}
& 100 \Delta^2 \|Q\|^2 D_M^2 \leq EC_AG,\ \  100 ED_A^2 \|Df (N)\|^2 \leq C_AG,\ \  100 E\Delta^2  D_A^2\left(\int_0^\infty\|M(s)\|\, ds\right)^2 \leq C_AG,\\
& \ \ \  {100 E \|Q\| D_A^2  \overline{R}^2} \leq {C_AG |L|},\ \  {100 \|Q\|^2 \overline{R}^2  } \leq {C_AG}{|L|} .
\end{align*}
We see that it is enough to choose
$$
|L| \geq \max\left\{ \frac{100 E_{max} \|Q\| D_A^2\overline{R}^2}{C_AG},\frac{100 \|Q\|^2\overline{R}^2}{C_AG} \right\}
$$
and the last two inequalities hold. We are now free to pick $E$ which satisfies
$$
0 < E \leq \frac{C_AG}{100 D_A^2 \|Df(N)\|^2}\ \ \textrm{and}\ \ E< E_{max},
$$
and we finally need to pick $\Delta$ such that
$$
\Delta^2 \leq \min \left\{  \frac{EC_AG}{100 \|Q\|^2 D_M^2}, \frac{C_AG}{100 E D_A^2\left(\int_0^\infty\|M(s)\|\, ds\right)^2} \right\}.
$$
The proof is complete.
\end{proof}

\section{Equilibria, persistence of connections in the limit $\varepsilon\to 0$.}\label{equi}
In this section we relate the equilibria of the  ODE \eqref{eq:ODE} with the equilibria of the problem \eqref{eq:main}--\eqref{main2}. We show that if $\varepsilon$ is small, then for every equilibrium $e$ of \eqref{eq:ODE} there exists an equilibrium $(0,e^{\varepsilon})$ of \eqref{eq:main}--\eqref{main2} in its vicinity and the perturbed system has no other equilibria.  Moreover, we show the upper semi-continuity result on the connections between the equilibria, that is, if the two equilibria of \eqref{eq:main}--\eqref{main2} are connected for a sequence of parameters $\varepsilon\to 0$, then the connection also exists for $\varepsilon=0$. We remind that the limit equation \eqref{eq:ODE}
has only a finite number of isolated and hyperbolic equilibria and the system is Morse--Smale, i.e. the intersections of their stable and unstable manifolds are always transversal.
\subsection{Continuation of equilibria} We begin from the result on the continuation of the equilibria of the system from $\varepsilon=0$ to the positive values.
\begin{thm}\label{thm:equilibria}
	There exists $\varepsilon_0>0$ and $R>0$ such that for every $\varepsilon\in [0,\varepsilon_0)$ if $e\in \R^d$ is an equilibrium for \eqref{eq:ODE} with an isolating block with cones $N^x$ then the problem governed by \eqref{eq:main}--\eqref{main2} has an equilibrium $(0,e^\varepsilon)$ which is unique in the set $\widetilde{N}^x =  \{ \|\eta\|\leq R \} \times N^x$.  Moreover $(0,e^\varepsilon)$ are the only equilibria for \eqref{eq:main}--\eqref{main2}.
	\end{thm}
\begin{proof}
	Denote by $\mathcal{E}$ the set of equilibria of \eqref{eq:ODE} and let $e\in \mathcal{E}$. Take $R$ and $N^x$ from Theorem \ref{thm:isol} (isolating set) and take $E$ satisfying the constraints from Lemma \ref{lem:E1} (the cone condition), Lemma \ref{lem:cone:par} (the cone condition with parameters) and Lemma \ref{lem:LapFunc1} (the Lyapunov function). Now take $\varepsilon \in [0,\varepsilon_0]$, where $\varepsilon_0$ satisfies all the constraints of the previous results: the constraint of  Lemma \ref{lem:E1} (the cone condition), Lemma \ref{lem:cone:par} (the cone condition with parameters), Lemma \ref{lem:LapFunc1} (the Lyapunov function) and Theorem \ref{thm:isol} (isolating set) and the constraints of Section \ref{sec:setup}. From Lemma \ref{lem:LapFunc1} we deduce that the equilibria of \eqref{eq:main} must have the $\eta$ component equal to zero. Theorem \ref{graph_t} together with Lemma \ref{lem:E1} imply that the  problem governed by \eqref{eq:main}--\eqref{main2} has a unique equilibrium in the set $\widetilde{N}^x$. We denote this equilibrium by $(0,e^\varepsilon)$. It must be $f^\varepsilon(e^\varepsilon) = 0$. We must show that problem \eqref{eq:main} does not have other equilibria than the ones which lie in $\widetilde{N}^x$.
	
	From Lemma \ref{lem:atr-bnd}  we deduce that if $f^\varepsilon(e^\varepsilon) = 0$ then $|e^\varepsilon| \leq R$. Assume that $e^\varepsilon \notin \bigcup_{y\in \mathcal{E}}\textrm{int}\, N^y$. Now let
	$$\beta = \min\left\{|f(x)|\,:\ |x|\leq {R}, x\notin \bigcup_{y\in \mathcal{E}}\textrm{int}\, N^y\right\}.
	$$
	This is a positive constant. We have
	$$
	|f^\varepsilon(e^\varepsilon)| \geq |f(e^\varepsilon)| - \varepsilon|e^\varepsilon| \int_0^\infty \|M(s)\|\, ds \geq  \beta - \varepsilon  R \int_0^\infty \|M(s)\|\, ds.
	$$
	Decreasing $\varepsilon$ is necessary we note that we must have $|f^\varepsilon(e^\varepsilon)| > 0$, a contradiction.
\end{proof}
\subsection{Connections are preserved in the limit.} We pass to the proof that if the heteroclinic connections exists for $\varepsilon>0$, they also must exist for $\varepsilon=0$. In the sequel we always assume that $E, \varepsilon_0, R$ satisfy the constraints which are needed for all results of Sections \ref{sec:setup} and \ref{sec:LapCC} to hold. 
\begin{df}
	Let $\varepsilon\geq 0$. The function $(\eta,x): \R \to L^2_A(\R^+)^d\times \R^d$ is a bounded complete (eternal) solution for \eqref{eq:main}--\eqref{main2} if for every $t\in \R$ the function $(\eta^{t+\cdot},x(t+\cdot)): \R \to L^2_A(\R^+)^d\times \R^d$ is a solution for \eqref{eq:main}--\eqref{main2}  and moreover $\sup_{t\in \mathbb \R}\left(|x(t)|+ \|\eta^t\|_{L^2_A(\R^+)^d}\right)$ is bounded.
\end{df}
The existence of the Lyapunov function in  Lemma \ref{lem:LapFunc1} and the  compactness of global attractors obtained in Lemma \ref{lem:atr-bnd} directly imply the following result.  \begin{lem}
	The pair $t\mapsto (\eta^t,x(t))$ is a bounded complete solution for \eqref{eq:main}--\eqref{main2} if and only if there exists two points $e^\varepsilon_1, e^\varepsilon_2 \in \mathbb{R}^d$ satisfying $f^\varepsilon(e^\varepsilon_1) = f^\varepsilon(e^\varepsilon_2) = 0$ such that
	$$
	\lim_{t\to -\infty} (\eta^t,x(t)) = (0,e^\varepsilon_1), \ \ \ 	\lim_{t\to \infty} (\eta^t,x(t)) = (0,e^\varepsilon_2)
	$$
	In such case we say that the exists a connection between the equilibria $e^\varepsilon_1$ and $e^\varepsilon_2$.
\end{lem}
In the next lemma we prove that the existing connections are preserved in the limit, see \cite[Proposition 4]{borto}. It is interesting to notice that we do not need uniform asymptotic compactness with respect to the memory variable because the limit system does not contain the memory term and it is enough to obtain the connection in the phase space $\R^d$.
\begin{lem}\label{lem:connection}
	If for a sequence $\varepsilon^n\to 0^+$ there exist connections between equilibria $e^{\varepsilon^n}_1$ and $e^{\varepsilon^n}_2$ through the system \eqref{eq:main}-\eqref{main2} where $\lim_{n\to\infty} e^{\varepsilon^n}_1 = e_1$ and $\lim_{n\to\infty} e^{\varepsilon^n}_2 = e_2$ then $e_1$ and $e_2$ are equilibria of \eqref{eq:ODE} and there exists a sequence of equilibria $e_1=g_1,\ldots, g_N=e_2$ such that there exist complete trajectories of \eqref{eq:ODE} which connect $e_i\to e_{i+1}$ for $i\in \{1,\ldots,N-1\}$.
\end{lem}
\begin{proof}
	The fact that $e_1$ and $e_2$ are equilibria of \eqref{eq:ODE} follows from the definition of $f^\varepsilon$.
	Denote by $(\eta_n^{(\cdot)},x_n(\cdot)):\R\to L^2_A(\R^+)^d\times \R^d$ the bounded complete solutions for $\epsilon_n$ such that for each $n$
	$$
	\lim_{t\to -\infty} (\eta_n^t,x_n(t)) = (0,e_1^{\epsilon^n}), \ \ \ 	\lim_{t\to \infty} (\eta_n^t,x_n(t)) = (0,e_2^{\epsilon^n}).
	$$
	Now $x_n$ and $\eta_n$ are bounded uniformly with respect to $t$. Estimate \eqref{eq:Lap-lbnd} implies that they are also bounded uniformly with respect to $n$. Hence $x_n'(t)$ is also uniformly bounded with respect to both $n$ and $t$.
	For every $\delta > 0$ there exists $t_n^1$ such that if only $t\leq t_n^1$ then $|x_n(t) - e_1^{\varepsilon_n}| \leq \delta$. From the Arzela--Ascoli lemma, using the diagonal argument we can construct a function $x^1:\R\to \R^n$ such that $x_n(t+t_n^1)\to u^1(t)$ uniformly for $t$ on every bounded time interval. Since
	$$
	x_n(t_n^1+t) = x_n(t_n^1) + \int_{t_n^1}^{t_n^1+t} \left(f(x_n(s)) + \varepsilon_n \int_{0}^\infty M(r)\, dr x_n(s) + \varepsilon_n \int_0^\infty M(r)\eta_n^s(r)\, dr\right) \, ds,
	$$
	we can pass to the limit with $n$ to infinity whence
	$$
	u_1(t) = u_1(0) + \int_0^t f(u_1(s))\, ds,
	$$
	i.e. $x_1$ solves \eqref{eq:ODE}. Now let $t \leq 0$ be fixed. We have
	$$
	|u_1(t) - e_1| \leq |u_1(t) - x_n(t+t_n^1)| + |x_n(t+t_n^1) - e_1^{\varepsilon_n}| + |e_1^{\varepsilon_n}- e_1| \leq |u_1(t) - x_n(t+t_n^1)| + \delta + |e_1^{\varepsilon_n}- e_1|.
	$$
	Passing with $n$ to infinity we deduce that
	$$
	|u_1(t) - e_1| \leq \delta.
	$$
	Since $\lim_{t\to -\infty} u_1(t) = e$, an equilibrium of \eqref{eq:ODE}, we deduce by taking $\delta$ small enough related to minimal distance between the equilibria of the system, that it must be $e = e_1$.
		Now $\lim_{t\to \infty}u_1(t) = g_2$, an equilibrium of \eqref{eq:ODE}. If $g_2=e_2$ the proof is complete. Otherwise for every $n$ there exists $k(n) \to \infty$ as  $n \to \infty$ and $\tau^2_{n}$ such that $|x_{k(n)}(\tau^2_n) - g_2|\leq \frac{1}{n}$.  Hence $x_{k(n)}(\tau^2_n+t)$ converges to $e_2$ \textbf{to $g_2$} uniformly on bounded time intervals. This means that for every sufficiently small $\delta$ and every $n$ there exists a maximal $t^2_n>\tau^2_n$ such that for $t\in [\tau^2_n,t^2_n]$ we have $|x_n(t)-g_2|\leq \delta$ and it must be $t^2_n-\tau^2_n \to \infty$ as $n\to \infty$.  Solutions $x_n(t_n^2+t)$, again from the Arzela--Ascoli lemma converge to $u_2(t)$, the solution of \eqref{eq:ODE}, uniformly on bounded time intervals. Moreover for every $t\leq 0$ we are able to find $n_0$ such that for every $n\geq n_0$ we have $\tau_n^2 < t+t_n^2$. Then
		$$
		|u_2(t)-g_2| \leq |u_2(t) - x_n(t+t_n^2)| + |x_n(t+t_n^2) - g_2| \leq |u_2(t) - x_n(t+t_n^2)| + \delta.
		$$
		Passing with $n\to \infty$ we deduce that $|u^2(t) - g_2| \leq \delta$ for every $t\leq 0$ and it is enough to choose $\delta$ sufficiently small so that $\lim_{t\to -\infty} u_2(t) = g_2$. Now, $\lim_{t\to \infty}u_2(t) = g_3$. If $g_3 = e_2$ the proof is complete. If not, we continue the procedure, which is always possible if the equilibrium is not $e_2$. Since the number of equilibria of \eqref{eq:ODE} is finite and the system is gradient, the procedure must end after finite number of steps, which concludes the proof.
\end{proof}
% TODO: \usepackage{graphicx} required
\begin{figure}
	\centering
	\includegraphics[width=0.7\linewidth]{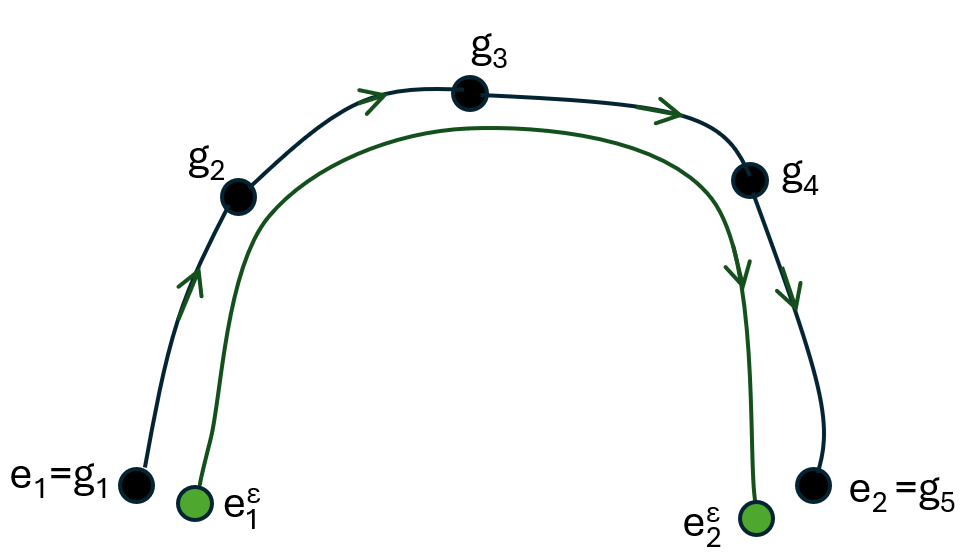}
	\caption{Illustration of the proof of Lemma \ref{lem:connection}. We pass to the limit in the $x$ component of the connection between equilibria $(0,e^\varepsilon_1)$ and $(0,e^\varepsilon_2)$. In the limit obtained via the Arzela--Ascoli lemma we obtain possibly several connections with intermediate points $g_2, g_3, \ldots, g_{N-1}$.}
	\label{fig:limitconnection}
\end{figure}

As the limit system \eqref{eq:ODE} is Morse--Smale, the existence of the sequence of connections $e_1=g_1\to g_2\to \ldots \to g_N=e_{2}$ implies the existence of connection $e_1\to e_2$, whence we can formulate the following Theorem
\begin{thm}\label{thm:upper_con}
	If for a sequence $\varepsilon^n\to 0^+$ there exist connections between equilibria $(0,e^{\varepsilon^n}_1)$ and $(0,e^{\varepsilon^n}_2)$ through the system \eqref{eq:main}--\eqref{main2} where $\lim_{n\to\infty} e^{\varepsilon^n}_1 = e_1$ and $\lim_{n\to\infty} e^{\varepsilon^n}_2 = e_2$ then $e_1$ and $e_2$ are equilibria of \eqref{eq:ODE} and there exists a connection $e_1\to e_2$ through the system \eqref{eq:ODE}.
\end{thm}

\section{Continuation of the intersection of manifolds.}\label{intersection}

\subsection{Disks and their reparameterization over tangent space.}

We begin with a definition of a embedded disk. If $V$ is a Banach space and  $B\subset W$ is an open set in a finite dimensional Banach space $W$, the $C^1$ mapping $g:B\to \textrm{range}(g)\subset V$ is an embedded disk if the linear map $Dg(p):W\to \textrm{range}(Dg(p))\subset V$ has the rank equal to the dimension of $W$.  We will work with the disks which have the coordinates, in which they are graphs of Lipschitz functions. to this end we need the following definition
\begin{df}
	Let $V,W$ be two disjoint Banach spaces, such that $V\subset Z$ and $W\subset Z$ for a Banach space $Z$. Let $B_1 = B(0,\delta_1)$ be a ball in $V$ and $B_2 = B(0,\delta_2)$ be a ball in $W$. The function $g:B_1\to B_2$ is a disk if $f$ is Lipschitz, where the domain is equipped with the norm of $V$ and the range, with the norm of $W$. We will say that the disk is $C^1$ if $g\in C^1(B_1)$. Moreover we will say that $\text{graph}(g) = \{v+g(v)\, :\ v\in B_1\}$.  
\end{df}
Let $g:B_1\to B_2$ be a $C^1$ disk, and let $v_0 \in B_1$ and $w_0 = g(v_0)$. Denote by $z_0\in \text{graph}(g)$ the point $z_0=v_0+g(v_0)$. Then $T_{z_0} g$, the tangent space of the disk $g$ at point $z_0$ is given by
$$
T_{z_0} g =\left\{ v+D g(v_0)v\,:\ v\in V\right\}.
$$    
We prove the following result that says that if $g_\varepsilon$ is a disk that is $C^1$ close to $g$ then locally $g_\varepsilon$ is a disk over the tangent space of $g$ with an arbitrarily small Lipschitz constant. 
\begin{lem}\label{lemma_disks}
	Let $\R^d = V+W$ where $V,W$ are subspaces such that $\textrm{dim}\, V = m$ and $\textrm{dim}\, W = d-m$ and let $B_1\subset V$ and $B_2\subset W$ be balls centered at zero. Let $P$ be a nonempty set of parameters and let  $g_\varepsilon:P\times B_1\to B_2$ for $\varepsilon>0$ and  $g:B_1\to B_2$ be $C^1$  disks such that
	$$
	\lim_{\varepsilon\to 0} \sup_{\eta\in P}\sup_{v\in B_1}\left(|g_\varepsilon(\eta,v)-g(v)|+\left|D_vg_\varepsilon (\eta,v)-Dg (v)\right|\right) =0
	$$
	and let $z_0\in \text{graph} (g)$. Let moreover  $M$ be a $d\times d$ nonsingular matrix such that $M\cdot (\R^m\times (0)_{d-m}) = T_{z_0} g$. Then for every $L>0$ there exists $\varepsilon_0>0$ and  $\delta_1>0$ and $\delta_2>0$, and balls $B(0,\delta_1)\subset \R^m \times (0)_{d-m}$, $B(0,\delta_2)\subset (0)_m\times \R^{d-m}$  such that for every $\varepsilon\in [0,\varepsilon_0]$ and for every $\eta\in P$ there exists a disk $h_{\varepsilon,\eta}: B(0,\delta_1)\to B(0,\delta_2)$ with a Lipschitz constant $L$ such that
	$$
	M\cdot \text{graph}(h_{\varepsilon,\eta}) \subset  \text{graph}(g_\varepsilon(\eta,\cdot))-z_0. 
	$$ 
	Moreover $\delta_2\leq c(\varepsilon_0) + c(\delta_1)\delta_1$, where $\lim_{\varepsilon\to 0}c(\varepsilon) = 0$ and $\lim_{\delta_1\to 0}c(\delta_1) = 0.$
	\end{lem}
	 \begin{proof}
	 	We will use the notation $\overline{G([a,b])} = \mathrm{conv}\{ g(x)\, :\  x = \lambda a + (1-\lambda b), \lambda\in [0,1] \}$. Translate the graph of $g_\varepsilon(\eta,\cdot)$ by $z_0=v_0+g(v_0)$. Such translation defines a function $b_{\varepsilon,\eta}(v)=g_{\varepsilon}(\eta,v+v_0)-g(v_0)$ on a ball $B$ in $V$ centered at zero with a radius independent of $\varepsilon$. This function is of class $C^1$, therefore if $z=v+b_{\varepsilon,\eta}(v)\in \text{graph}(b_{\varepsilon,\eta})$, then 
	 		$$
	 	z = v+g_\varepsilon(\eta,v+v_0)-g(v_0)\in v+g_\varepsilon(\eta,v_0)+ \overline{D_v g_\varepsilon(\eta,[v_0,v+v_0])}v-g(v_0).
	 	$$
	 	This point can be written as 
	 	\begin{align*}
	 	z & \in v+D g(v_0)v+g_\varepsilon(\eta,v_0)-g(v_0)+ \overline{\left(D_v g_\varepsilon(\eta,\cdot)-D g(\cdot)\right)([v_0,v+v_0])}v\\
	 	& +\overline{D g([v_0,v+v_0])-D g(v_0)}v.
	 	\end{align*}
	 	We will write
	 	$$z = v+D g(v_0)v + \Delta_1 + \Delta_2 + \Delta_3.
	 	$$
	 	Observe that  $v+D g(v_0)v \in T_{z_0}g$. Denoting $\Pi_{M\cdot (\overline{v},0)}$ the projection on the tangent space $T_{z_0}f$ and by $\Pi_{M\cdot (0,\overline{w})}$ the complementary projection, we can represent the considered point $z\in \text{graph}(b_{\varepsilon,\eta})$ as
	 	$$
	 	z = v+D g(v_0)v + \Pi_{M\cdot (\overline{v},0)}(\Delta_1+\Delta_2+\Delta_3) + \Pi_{M\cdot (0,\overline{w})}(\Delta_1+\Delta_2+\Delta_3).
	 	$$
	 	We first prove that there exists $\delta_1 > 0$ such that for every $(\overline{v},0) \in B(0,\delta_1)$ there exists $v\in B$ such that
	 	$$
	 	(\overline{v},0) = M^{-1} \left(v+D g(v_0)v + \Pi_{M\cdot (\overline{v},0)}(\Delta_1+\Delta_2+\Delta_3)\right).
	 	$$
	 	The mapping
	 	$$
	 	V\ni v\mapsto M^{-1} \left(v+ Dg(v_0)v  \right)\in \R^m \times (0)_{d-m}
	 	$$
	 	is a linear invertible mapping from $m$ dimensional space into $m$ dimensional space. Therefore  for every $\delta$ with $B(0,\delta) \subset B$ we can find a ball $B(0,{\delta}_1(\delta)) \subset \R^m \times (0)_{d-m}$ such that
	 	$$
	 	\{ (\overline{v},0) \in   B(0,{\delta}_1)  \} \subset \left\{ M^{-1} \left(v+D g(v_0)v \right)\, : \ v\in B(0, \delta)\right\}.
	 	$$
	 	Consider the homotopy
	 	$$
	 	B(0, \delta) \times[0,1] \ni (v,\theta ) \mapsto H(v,\theta )  =  M^{-1} \left(v+D g(v_0)v + \theta \Pi_{M\cdot (\overline{v},0)}(\Delta_1+\Delta_2+\Delta_3) \right).
	 	$$
	 	Now for constants $C_1, C_2 > 0 $ we have
	 	$$
	 	|H(v,\theta)| \geq C_1|v| - C_2|\Delta_1| - C_2|\Delta_2| - C_2|\Delta_3|\  \ \textrm{for}\ \  (v,\theta )\in B(0, \delta) \times[0,1].
	 	$$
	 	Assuming that $|v| = \delta$, by taking $\delta$ small enough we  obtain $C_2|\Delta_3| \leq \frac{C_1}{6}\delta$ and then, by taking  $\varepsilon$ small enough we get 
	 	$C_2|\Delta_1| \leq \frac{C_1}{6}\delta$ and $C_2|\Delta_2| \leq \frac{C_1}{6}\delta$. Hence for $|v| = \delta$ we have  $|H(v,\theta)| >  \frac{C_1}{2}\delta$, and, by the homotopy invariance of the Brouwer degree if only $\delta$ and $\varepsilon$ are small enough  we obtain the existence of ${\delta}_1(\delta)$ such that for every  $(\overline{v},0) \in B(0,{\delta}_1)$ there exists $(0,\overline{w})\in (0)_m\times 
	 	\R^{d-m}$ for which the point $M\cdot((\overline{v},\overline{w}))$ belongs to the graph of $b_{\varepsilon,\eta}$. The estimate on $|M^{-1}\Pi_{M\cdot(0,\overline{w})}(\Delta_1+\Delta_2+\Delta_3)|$ gives us the bound on the radius $\delta_2$. 
	 	
	 	In the next step we  show that for a given  $\overline{v}$ the point $\overline{w}$ is unique and that the dependence $(0,\overline{w}) = h_{\varepsilon,\eta}(\overline{v},0)$  is Lipschitz with a constant that can be made arbitrarily small by decreasing, if necessary, the radius ${\delta}_1$ and $\varepsilon$.   Consider two points in graph of $b_{\varepsilon,\eta}$, and denote them by $z_1=v_1+g_\varepsilon(\eta,v_1+v_0)-g(v_0)=M\cdot (\overline{v}_1,\overline{w}_1)$ and $z_2=v_2+g_\varepsilon(\eta,v_2+v_0)-g(v_0)=M\cdot (\overline{v}_2,\overline{w}_2)$.
	 	Now
	 	\begin{align*}		
	 	& z_1-z_2=M\cdot (\overline{v}_1-\overline{v}_2,\overline{w}_1-\overline{w}_2)  = g_\varepsilon(\eta,v_1+v_0)-g_\varepsilon(\eta,v_2+v_0)+v_1-v_2\\
	 	&\ \ \ \in v_1-v_2 + Dg(v_0)(v_1-v_2) \\
	 	& \ \ \ + \overline{\left(D_v g_\varepsilon(\eta,\cdot)-Dg(\cdot)\right)([v_0+v_1,v_0+v_2])}(v_1-v_2)
	 	+\overline{\left(Dg(v_0+B(0,\delta))-Dg(v_0)\right)}(v_1-v_2).
	 	\end{align*}
	 	We will write
	 	$$
	 	z_1-z_2 = v_1-v_2 + Dg(v_0)(v_1-v_2) + \Delta_1 + \Delta_2.
	 	$$
	 	It follows that
	 	\begin{align*}
	 		& M\cdot  (\overline{v}_1-\overline{v}_2,0) =  v_1-v_2 + Dg(v_0)(v_1-v_2)  + \Pi_{M\cdot  (\overline{v},0)}(\Delta_1+\Delta_2),
	 	\end{align*}
	 	and, 
	 	$$
	 	M\cdot  (0,\overline{w}_1-\overline{w}_2) =  \Pi_{M\cdot  (0,\overline{w})}(\Delta_1+\Delta_2).
	 	$$
	 	The first of the above two equations implies that there exist constants $C_1, C_2(\varepsilon), C_3(\delta_1) > 0$ such that
	 	$$
	 	|(\overline{v_1}-\overline{v}_2,0) | \geq C_1 |v_1-v_2| - (C_2(\varepsilon)+C_3(\delta_1))|v_1-v_2|,
	 	$$
	 	where $C_2$ and $C_3$ can be made as small as we need by taking sufficiently small $\delta_1$ and $\varepsilon$. Moreover
	 	$$
	 	| (0,\overline{w}_1-\overline{w}_2) | \leq (C_4(\varepsilon)+C_5(\delta_1)) |v_1-v_2|,
	 	$$
	 	where $C_4, C_5$ again can be made as small as necessary by taking sufficiently small $\delta_1$ and $\varepsilon>0$. Both above inequalities imply that 
	 	$$| (0,\overline{w}_1-\overline{w}_2) | \leq  \frac{C_4(\varepsilon)+C_5(\delta_1)}{C_1-C_2(\varepsilon)-C_3(\delta_1)}|(\overline{v}_1-\overline{v}_2,0) | $$
	 	 with the Lipschitz constant being as small as we need, which can be obtained by taking small $\delta_1$ and small $\varepsilon>0$. This gives the restriction on the radius $\delta_1$ and $\varepsilon$,  in order to get the desired Lipschitz constant $L$. \end{proof}
	 	We have locally reparametrized the graphs of $g_\varepsilon(\eta,\cdot)$  over the tangent space of the disk $g$ to get the family of disks with arbitrarily small Lipschitz constant. Note that the radius  $\delta_1$ can be arbitrarily decreased in the above lemma. In the next result we study the dependence on  $\eta$, which was treated in the previous result as a parameter.
	 	\begin{lem}\label{lem:stable}
	 		Under assumptions of the previous lemma let $P$ be an open and  bounded set in the Banach space $X$. Assume that $g_\varepsilon:P\times B_1 \to B_2$ (with $g_0(\eta,v)=g(v)$) is Fr\'{e}chet differentiable with respect to $\eta\in P$ and $\|D_\eta g_\varepsilon(\eta,v)\|_{\mathcal{L}(X;W)}\leq E$. Then, denoting $h_\varepsilon(\eta,\cdot) = h_{\varepsilon,\eta}(\cdot)$, we obtain
	 		$$	
	 		|h_\varepsilon(\eta_1,(\overline{v}_1,0))-h_\varepsilon(\eta_2,(\overline{v}_2,0))| \leq L|(\overline{v}_1-\overline{v}_2,0)| + K E\|\eta_1-\eta_2\|_X \ \textrm{for} \ \eta_1,\eta_2\in P, (\overline{v}_1,0), (\overline{v}_2,0) \in B(0,\delta_3),
	 		$$
	 		for every $\varepsilon\in [0,\varepsilon_0]$
	 		where, as in Lemma \ref{lemma_disks}, $L$ can be made arbitrarily small by decreasing, if necessary, constants $\varepsilon_0$ and $\delta_1$, and $K$ is a constant independent of $\varepsilon$.  
	 		\end{lem}
	 		 \begin{proof}From Lemma \ref{lemma_disks}, the function $h_\varepsilon:P\times B(0,\delta_1)\to B(0,\delta_2)$ is well defined. We need to show the Lipschitz condition. Choose $\eta_1,\eta_2\in P$, and consider the two points, one in the graph of $b_{\varepsilon,\eta_1}$, and the second one in the graph of in the graph of $b_{\varepsilon,\eta_2}$. Denote them by $z_1=v_1+g_\varepsilon(\eta_1,v_1+v_0)-g(v_0)=M\cdot(\overline{v}_1,\overline{w}_1)$ and $z_2=v_2+g_\varepsilon(\eta_2,v_2+v_0)-g(v_0)=M\cdot(\overline{v}_2,\overline{w}_2)$. Then, for $\overline{v}_1,\overline{v}_2 \in B(0,\delta_3)$ we have $\overline{w}_1 = h_{\varepsilon}(\eta_1,\overline{v}_1)$ and $\overline{w}_2 = h_{\varepsilon}(\eta_2,\overline{v}_2)$.  We have
	 		 	\begin{align*}
	 		 		&   z_1-z_2= M\cdot (\overline{v}_1-\overline{v}_2,\overline{w}_1-\overline{w}_2) = v_1-v_2 + g_\varepsilon(\eta_1,v_1+v_0)-g_\varepsilon(\eta_2,v_2+v_0) \\
	 		 		& \ \in \overline{D_vg_\varepsilon([(\eta_2,v_2+v_0),(\eta_1,v_1+v_0)])}(v_1-v_2)+\overline{D_\eta g_\varepsilon([(\eta_2,v_2+v_0),(\eta_1,v_1+v_0)])}(\eta_1-\eta_2)  +v_1-v_2 \\
	 		 		& \ = \overline{(D_vg_\varepsilon(\cdot)-Dg(\cdot))([(\eta_2,v_2+v_0),(\eta_1,v_1+v_0)])}(v_1-v_2)\\
	 		 		& \qquad \qquad +\left(\overline{Dg([v_2+v_0,v_1+v_0])}-D g(v_0)\right)(v_1-v_2) \\
	 		 		&\qquad\qquad  + \overline{D_\eta g_\varepsilon([(\eta_2,v_2+v_0),(\eta_1,v_1+v_0)])}(\eta_1-\eta_2)  +v_1-v_2 + D g(v_0)(v_1-v_2)\\
	 		 		&  = \Delta_1+\Delta_2+\Delta_3 + v_1-v_2 + D g(v_0)(v_1-v_2).
	 		 	\end{align*}
	 		 	
	 		 	We project this formula on the tangent space $T_{z_0}g$ and its complement. We obtain
	 		 	\begin{align*}
	 		 		& M\cdot  (\overline{v}_1-\overline{v}_2,0) \in  D g(v_0)(v_1-v_2)+v_1-v_2  + \Pi_{M\cdot  (\overline{v},0)}(\Delta_1+\Delta_2+\Delta_3).
	 		 	\end{align*}
	 		 	This means that
	 		 	$$
	 		 	|(\overline{v}_1-\overline{v}_2,0) | \geq C_1 |v_1-v_2| - (C_2(\varepsilon)+C_3(\delta_1))|v_1-v_2| - C_4E\|\eta_1-\eta_2\|_X,
	 		 	$$
	 		 	where $C_2(\varepsilon) \to 0$  as $\varepsilon\to 0$, $C_3(\delta_1)\to 0$ as $\delta_1 \to 0$ and $C_1, C_4$ are constants. Now
	 		 	\begin{align*}
	 		 		& M\cdot  (0,\overline{w}_1-\overline{w}_2) \in  \Pi_{M\cdot  (0,\overline{w})}(\Delta_1+\Delta_2+\Delta_3).
	 		 	\end{align*}
	 		 	It follows that
	 		 	$$
	 		 	|(0,\overline{w}_1-\overline{w}_2) |\leq (C_5(\varepsilon)+C_6(\delta_1))|v_1-v_2|+C_7E\|\eta_1-\eta_2\|_X,
	 		 	$$
	 		 	where, again, $\lim_{\varepsilon\to 0}C_5(\varepsilon) = 0$ and $\lim_{\delta_1\to 0}C_6(\delta_1) = 0$. Summarizing, we obtain
	 		 	$$
	 		 	|(0,\overline{w}_1-\overline{w}_2) |\leq \frac{C_5(\varepsilon)+C_6(\delta_1)}{C_1-C_2(\varepsilon)-C_3(\delta_1)}|(\overline{v}_1-\overline{v}_2,0) | + \left(C_7E+\frac{C_5(\varepsilon)+C_6(\delta_1)}{C_1-C_2(\varepsilon)-C_3(\delta_1)}C_4E\right)\|\eta_1-\eta_2\|_X.
	 		 	$$
	 		 	Note that decreasing $\varepsilon$ and the radius $\delta_1$ of the box can make the constant $\frac{C_5(\varepsilon)+C_6(\delta_1)}{C_1-C_2(\varepsilon)-C_3(\delta_1)}$ arbitrarily small, which ends the proof of the lemma.
	 		 	
	 		 	\end{proof} 
	 	\subsection{Transport of disks.} The following lemma says that a finite dimensional $C^1$ disk stays a disk after composition with a mapping. This new disk can be reparameterized as a disk over the tangent space, with arbitrarily small Lipschitz constant. 
	 	
	 		 	\begin{lem}\label{lem:transport_disks} Let $\R^d=V+W$ where $V,W$ are two subspaces such that $\text{dim}\, V = n$ and $\text{dim}\, W = d-n$.  Let $g:B_1\to B_2$ be a $C^1$  disk, where $B_1\subset V$ and $B_2 \subset W$ are balls.  
	 		Let moreover $G:\R^d\to \R^d$ be a $C^1$ mapping.  Assume that $z_0\in \text{graph}(g)$ be such that $DG(z_0)$ is invertible and let $M$ be a $d\times d$  invertible matrix such that $M\cdot((0)_{d-n}\times \R^n) = T_{G(z_0)}G\circ (\textrm{Id}_{B_1}+g)$. Then for every $L>0$ there exists  $\delta_3, \delta_4>0$, balls $B(0,\delta_3)\subset  (0)_{d-n} \times \R^n$, $B(0,\delta_4)\subset \R^{d-n} \times (0)_n$ and the  disk $h:B(0,\delta_3)\to B(0,\delta_4)$ with the Lipschitz constant  $L$ such that 
	 		$$
	 		M \cdot \text{graph}(h) \subset G(\text{graph}(g))-G(z_0).
	 		$$
	 		Moreover $\delta_4\leq c(\delta_3)\delta_3$, where $c(\delta_3)\to 0$ as $\delta_3\to 0$. 
	 	\end{lem}
	 	\begin{proof}
	 		Denote $z_0=v_0+g(v_0)$. Projection onto the tangent space $T_{G(z_0)}G\circ (Id_{B_1}+g)$ will be denoted by $\Pi_{M\cdot (0,\overline{v})}$ and $\Pi_{M\cdot (\overline{w},0)}$ will be the complementary projection. There exists a ball $B\subset V$ such that $v_0 + B \subset B_1$. Define the function $b:B\to W$ as $b(v)=g(v_0+v)-g(v_0)$ so that $Db(v) = Dg(v_0+v)$. For $v\in B$ we have
	 		\begin{align*}
	 			G(v_0+v+g(v_0+v)) = G(z_0+v+b(v)) = G\left(z_0+v+Dg(v_0)v+\Delta(v)\right),
	 		\end{align*}
	 		where $\Delta(v) \in o(|v|)$. Denote
	 		$$
	 		\Delta_1 (v) = v+Dg(v_0)v+\Delta(v).
	 		$$
	 		 Then
	 		\begin{align*}
	 			&  G(v_0+v+g(v_0+v)) = G(z_0)+DG(z_0)\Delta_1(v) + \Delta_2(\Delta_1(v))\\
	 			& \qquad  = G(z_0)+DG(z_0)\left(v+Dg(v_0)v\right)+DG(z_0)\Delta(v) + \Delta_2(\Delta_1(v)),
	 		\end{align*}
	 		with $\Delta_2(\Delta_1(v))\in o(|v|)$.
	 		Denote $DG(z_0)\Delta(v) + \Delta_2(\Delta_1(v)) = \Delta_3(v)$. This quantity belongs to $o(|v|)$. Then we have
	 		$$
	 		G(v_0+v+g(v_0+v)) - G(z_0) =DG(z_0)\left(v+Dg(v_0)v\right)+\Pi_{M\cdot (\overline{w},0)}\Delta_3(x) + \Pi_{M\cdot (0,\overline{v})}\Delta_3(v),
	 		$$
	 		and the expression $DG(z_0)\left(v+Dg(v_0)v\right)+\Pi_{M\cdot (0,\overline{v})}\Delta_3(v)$ belongs to the tangent space $T_{G(z_0)}G\circ (Id_{B_1}+g)$. We need to show that for every ball $B(0,\delta)\subset B$
	 		there exists $\delta_3(\delta)$ such that for every $(0,\overline{v})\in B(0,\delta_3)$ there exists $v\in B(0,\delta)$ such that
	 		$$
	 		M\cdot (0,\overline{v}) = DG(z_0)\left(v+Dg(v_0)v\right)+\Pi_{M\cdot (0,\overline{v})}\Delta_3(v).
	 		$$
	 		The argument follows by homotopy. Indeed, the invertibility of the linear mapping $v\mapsto DG(z_0)(v+Dg(v_0)v)$ implies that  for every $\delta>0$ with $B(0,\delta)\subset B$ there exists ${\delta_3}(\delta)$ such that
	 		$$
	 		\{ (0,\overline{v}) \in B(0,\delta_3)  \} \subset \left\{ M^{-1} DG(z_0)\left(v+Dg(v_0)v\right)\, : \ v\in B(0, \delta)\right\},
	 		$$
	 		and, decreasing $\delta$ if necessary, the result follows similarly as in Lemma \ref{lemma_disks} by considering the homotopy
	 		$$
	 		B(0,\delta)\times [0,1]\ni (v,\theta) \mapsto M^{-1}\left(DG(z_0)\left(v+Dg(v_0)v\right)+\theta \Pi_{M\cdot (0,\overline{v})}\Delta_3(v)\right).
	 		$$ Now, we define $h(0,\overline{v}) = M^{-1}(\Pi_{M\cdot (\overline{w},0)}\Delta_3(v))$. The estimate on $|M^{-1}(\Pi_{M\cdot (\overline{w},0)}\Delta_3(v))|$ gives the bound on the radius $\delta_4$ which behaves like $c(\delta_3)\delta_3$ with $c(\delta_3)\to 0$ as $\delta_3\to 0$. 
	 		
	 		To demonstrate that the point $(\overline{w},0)$ is uniquely determined for a given $(0,\overline{v})$ and the Lipschitz condition holds with arbitrarily small constant, consider the two points in the graph of $b$ denoting them by $v_1+b(v_1)$ and $v_2+b(v_2)$. We also denote $z_1=z_0+v_1+b(v_1)$ and $z_2=z_0+v_2+b(v_2)$.  The function $G:\R^d\to \R^d$ can be treated as the function of two variables $G:V\times W \to \R^d$ and then its partial derivatives are denoted by $D_vG$ and $D_wG$. Consider the difference
	 		\begin{align*}
	 			& G(z_1)-G(z_2)=G(z_0+v_1+b(v_1)) - G(z_0+v_2+b(v_2))\\
	 			& \qquad  \in \overline{D_vG([z_1,z_2])} (v_1-v_2) + \overline{D_wG([z_1,z_2])}(b(v_1)-b(v_2)) \\
	 			& \qquad \subset \left(\overline{D_vG([z_1,z_2])}   + \overline{D_wG([z_1,z_2])}\ \overline{D b([v_1,v_2])}\right)(v_1-v_2) \\
	 			& \qquad \subset  D_vG(z_0)(v_1-v_2)+ D_wG(z_0)Dg(v_0)(v_1-v_2)\\
	 			& \qquad  + D_wG(z_0)\left(\overline{D g([v_0+v_1,v_0+v_2])} - Dg(v_0) \right)(v_1-v_2) \\
	 			& \qquad \qquad +   \left(\overline{D_vG([z_1,z_2])} - D_vG(z_0) + \left(\overline{D_wG([z_1,z_2])}-D_wG(z_0)\right)\overline{D g([v_0+v_1,v_0+v_2])}\right)(v_1-v_2)\\
	 			& \qquad\qquad = D_vG(z_0)(v_1-v_2)+ D_wG(z_0)Dg(v_0)(v_1-v_2) + \Delta.
	 		\end{align*}
	 		Now denote $G(z_1) = M\cdot (\overline{w}_1,\overline{v}_1)$ and  $G(z_2) = M\cdot (\overline{w}_2,\overline{v}_2)$. Because 
	 		$$D_vG(z_0)(v_1-v_2)+ D_wG(z_0)Dg(v_0)(v_1-v_2) = DG(z_0)(I_{B_1}+Dg(v_0))(v_1-v_2) \in T_{G(z_0)}G\circ (I_{B_1}+g),$$ we have
	 		\begin{align*}
	 			& M\cdot  (0,\overline{v}_1-\overline{v}_2) \in DG(z_0)(I_{B_1}+Dg(v_0))(v_1-v_2) +  \Pi_{M\cdot (0,\overline{v})}\Delta,
	 		\end{align*}
	 		and
	 		\begin{align*}
	 			& M\cdot  (\overline{w}_1-\overline{w}_2,0) \in  \Pi_{M\cdot (\overline{w},0)}\Delta.
	 		\end{align*}
	 		Now, as $M^{-1} DG(z_0)(I_{B_1}+Dg(v_0))$ is invertible  we deduce that
	 		$$
	 		|(0,\overline{v}_1-\overline{v}_2) | \geq C_1 |v_1-v_2| - C_2(\delta_3)|v_1-v_2|,
	 		$$
	 		and
	 		$$
	 		| (\overline{w}_1-\overline{w}_2,0) | \leq C_3(\delta_3) |v_1-z_2|,
	 		$$
	 		where $C_1$ is a fixed constant and $C_2, C_3$ can be made as small as we need by taking sufficiently small radius $\delta_3$. This implies the required Lipschitz condition 
	 		$$
	 		| (\overline{w}_1-\overline{w}_2,0) | \leq \frac{C_3(\delta_3)}{C_1-C_2(\delta_3)}|(0,\overline{v}_1-\overline{v}_2) |.
	 		$$
	 		\end{proof}
	 	
	In the next lemma we consider $G_\varepsilon$, a perturbation of the mapping $G$ defined in the extended space $X\times \R^d$ and disks $(g_{\varepsilon,X},g_{\varepsilon,W})$, where $g_{\varepsilon,W}$ are perturbations of $g$ and $g_{\varepsilon,X}$ are disks with values in $X$. We provide conditions under which the disks $g_{\varepsilon,W}$ transported in the extended space $X\times \R^d$ stay close  to the disk $g$ transported in $\R^d$.    	 	
	 	
	 	\begin{lem} \label{lem:transport_eps} Let $\R^d=V+W$ where $V,W$ are two subspaces such that $\text{dim}\, V = n$ and $\text{dim}\, W = d-n$. Let moreover $X$ be a Banach space. Let $g:B_1\to B_2$ be a $C^1$  disk, where $B_1\subset V$ and $B_2 \subset W$ are balls.  Let moreover $g_\varepsilon=(g_{\varepsilon,X},g_{\varepsilon,W}):B_1\to B_2\times P$ be a family of $C^1$  disks given for $\varepsilon>0$ with $P\subset X$ being an open and bounded set such that $\sup_{v\in B_1}\left\|Dg_{\epsilon,X}(v)\right\|_{\mathcal{L}(V;X)} \leq E(\varepsilon)$, with $E(\varepsilon)\to 0$ as $\varepsilon\to 0$. Let moreover 
	 		$$
	 		\lim_{\varepsilon\to 0}\sup_{v\in B_1}\left(|g_{\varepsilon,W}(v)-g(v)|+\left|D g_{\varepsilon,W}(v)-D g(v)\right|\right) = 0.
	 		$$
	 		Next, let $G_\varepsilon=(F_{\varepsilon,X},G_{\varepsilon,d}):X \times \R^d\to X\times \R^d$ be a family of $C^1 $mappings and let $G:\R^d\to \R^d$ be a $C^1$ mapping such that
	 		$$
	 		\lim_{\varepsilon\to 0}\sup_{\eta\in P}\sup_{u\in B_1+B_2}\left(|G_{\varepsilon,d}(\eta,u)-G(u)|+\left|D_uG_{\varepsilon,d}(\eta,u)-DG(u)\right|\right)=0.
	 		$$ Moreover let $\|D G_{\varepsilon}(\eta,u)\|_{\mathcal{L}(X\times \R^d;X\times \R^d)}$ be bounded uniformly in $\varepsilon$ for $(\eta,u)$ in bounded sets in $X\times \R^d$. Assume that $z_0 = v_0+g(v_0) \in \text{graph}(g)$ is such that $DG(z_0)$  is invertible and let $M$ be a $d\times d$  invertible matrix such that $M\cdot((0)_{d-n}\times \R^n) = T_{G(z_0)}G\circ (\textrm{Id}_{B_1}+g)$. Pick $L>0$. There exist $\varepsilon_0>0$, $\delta_3, \delta_4>0$, with $\delta_4\leq c(\varepsilon_0)+c(\delta_3)\delta_3$ where $c(r)\to 0$ as $r\to 0$,  balls $B(0,\delta_3)\subset  (0)_{d-n} \times \R^n$, $B(0,\delta_4)\subset  \R^{d-n} \times (0)_n$ and a bounded set  $Q \subset X$ such that for every $\varepsilon\in (0,\varepsilon_0]$ there exists the  disk $h_\varepsilon=(h_{\varepsilon,W},h_{\varepsilon,X}):B(0,\delta_3)\to B(0,\delta_4)\times Q$ with the Lipschitz constant for $h_{\varepsilon,W}$ equal to $L$ and 
	 		$$
	 		 M \cdot \text{graph}(h_{\varepsilon,W}) \subset G_{\varepsilon,d}(\text{graph}(g_\varepsilon))-G(z_0).
$$
Moreover $h_{\varepsilon,X}(0,\overline{v}) = G_{\varepsilon,X}(g_{\varepsilon,X}(v),v+g_{\varepsilon,W}(v))$, where $v$ is such that $	 		 M \cdot ((0,\overline{v})+h_{\varepsilon,W}(0,\overline{v})) =  G_{\varepsilon,d}(g_{\varepsilon,X}(v),v+g_{\varepsilon,W}(v)))-G(z_0)
$.
	 		\end{lem}
	 		\begin{proof}
	 			In the lemma statement we denote $\textrm{graph}(g_\varepsilon) = \{ (g_{\varepsilon,X}(p),v+g_{\varepsilon,W}(p))\,:\ p\in B_1 \}$. Let $B\subset V$ be a ball such that $v_0+B\subset B_1$.
 If $v\in B$, then the point in the graph of $g_\varepsilon$ is denoted by
	 			\begin{align*}
	 				&  (g_{\varepsilon,X}(v+v_0),v+v_0+g_{\varepsilon,W}(v+v_0)) \\
	 				& \ \ =    \left(g_{\varepsilon,X}(v_0)+Dg_{\varepsilon,X}(v_0)v+\Delta_1(v),v_0+v+g(v_0)+g_{\varepsilon,W}(v_0)-g(v_0)+
	 				Dg_{\varepsilon,W}(v_0)v+\Delta_2(v)\right)\\
	 				& \ \ = (g_{\varepsilon,X}(v_0),z_0) + \left(Dg_{\varepsilon,X}(v_0)v,v+Dg_{\varepsilon,W}(v_0)v\right) + (0,g_{\varepsilon,W}(v_0)-g(v_0)) +  \left(\Delta_1(v),\Delta_2(v)\right)\\
	 				& \ \  = (g_{\varepsilon,X}(v_0),z_0) + \left(Dg_{\varepsilon,X}(v_0)v,v+Dg(v_0)v\right) +   \left(0,\left(Dg_{\varepsilon,W}(v_0)-Dg(v_0)\right)v\right) \\
	 				& \qquad \qquad + (0,g_{\varepsilon,W}(v_0)-g(v_0)) +  \left(\Delta_1(v),\Delta_2(v)\right)\\
	 				& \ \  = (g_{\varepsilon,X}(v_0),z_0) + I + II + III + IV.
	 			\end{align*}
	 			Now, consider the image of this point by the map $G_\varepsilon=(G_{\varepsilon,X},G_{\varepsilon,d})$.
	 			\begin{align*}
	 				& G_{\varepsilon,d}(g_{\varepsilon,X}(v+v_0),v+v_0+g_{\varepsilon,W}(v+v_0)) \\
	 				& \  = G(z_0) + (G_{\varepsilon,d}(g_{\varepsilon,X}(v_0),z_0)-G(z_0)) + D G_{\varepsilon,d}(g_{\varepsilon,X}(v_0),z_0)I \\
	 				& \quad + DG_{\varepsilon,d}(g_{\varepsilon,X}(v_0),z_0)(II+III+IV)\\
	 				& \quad +(\overline{DG_{\varepsilon,d}([(g_{\varepsilon,X}(v_0),z_0),(g_{\varepsilon,X}(v+v_0),v+v_0+g_{\varepsilon,W}(v+v_0))]}-DG_{\varepsilon,d}(g_{\varepsilon,X}(v_0),z_0))\cdot\\
	 				& \qquad \qquad \cdot(I+II+III+IV)\\
	 				& \  = G(z_0) + DG(z_0)\left(v+Dg(v_0)v\right)  + \left(D G_{\varepsilon,d}(g_{\varepsilon,X}(v_0),z_0)-DG(z_0)\right)I + \\
	 				& \quad  +  (G_{\varepsilon,d}(g_{\varepsilon,X}(v_0),z_0)-G(z_0)) + D G_{\varepsilon,d}(g_{\varepsilon,X}(v_0),z_0)(II+III+IV)\\
	 				& \quad +(\overline{DG_{\varepsilon,d}([(g_{\varepsilon,X}(v_0),z_0),(g_{\varepsilon,X}(v+v_0),v+v_0+g_{\varepsilon,W}(v+v_0))]}-DG_{\varepsilon,d}(g_{\varepsilon,X}(v_0),z_0))\cdot\\
	 				& \qquad \qquad \cdot(I+II+III+IV).
	 			\end{align*}
	 			We can rewrite the above equation as 
	 			$$
	 			G_{\varepsilon,d}(g_{\varepsilon,X}(v+v_0),v+v_0+g_{\varepsilon,W}(v+v_0)) - G(z_0) = DG(z_0)\left(v+Dg(v_0)v\right)  + \Delta_3 + \Delta_4+\Delta_5 + \Delta_6. 
	 			$$
	 			As in Lemma \ref{lem:transport_disks}  we must prove that for every ball $B(0,\delta)\subset B$ there exists $\delta_3(\delta)$ such that for every $(0,\overline{v}) \in B(0,\delta_3)$ we can find $v\in B(0,\delta)$ such that
	 			$$
	 			M\cdot(0,\overline{v}) = DG(z_0)\left(v+Dg(v_0)v\right)  + \Pi_{M\cdot(0,\overline{v})}(\Delta_3 + \Delta_4+\Delta_5 + \Delta_6).
	 			$$
	 			The proof follows, again, by the homotopy invariance of the Brouwer degree as in the previous Lemma. Indeed, there exist constants $c(\varepsilon)\to 0$ as $\varepsilon\to 0$ and $c(|v|)\to 0$ as $|v|\to 0$, for which we obtain   
	 			\begin{align*}	 				
	 			& |\Delta_3| = |\left(D G_{\varepsilon,d}(g_{\varepsilon,X}(v_0),z_0)-DG(z_0)\right)I| \leq c(\varepsilon) |I| \leq c(\varepsilon)|v|, \\
	 			& |\Delta_4| = |(G_{\varepsilon,d}(g_{\varepsilon,X}(v_0),z_0)-G(z_0))| \leq c(\varepsilon),\\
	 			& |\Delta_5| = |D G_{\varepsilon,d}(g_{\varepsilon,X}(v_0),z_0)(II+III+IV)| \leq C|II+III+IV| \leq c(\varepsilon)+(c(\varepsilon)+c(|v|))|v|,\\
	 			& |\Delta_6| \leq c(|v|)(|I+II+III+IV|)\leq c(|v|)( |v| + c(\varepsilon)).
	 			\end{align*}
	 			Thus 
	 			$$
	 			|\Delta_3+\Delta_4+\Delta_5+\Delta_6| \leq c(\varepsilon) + c(|v|)|v|,
	 			$$ 
	 			 and the homotopy to be used is. 
	 			$$
	 			B(0,\delta)\times [0,1]\ni (v,\theta) \mapsto M^{-1}(DG(z_0)\left(v+Dg(v_0)v\right)  + \theta \Pi_{M\cdot(0,\overline{v})}(\Delta_3 + \Delta_4+\Delta_5 + \Delta_6)).
	 			$$
	 			Now, consider two points in the graph of $g_\varepsilon$. Denote them by
	 			$(\eta_1,z_1) = (g_{\varepsilon,X}(v_1+v_0),v_1+v_0 + g_{\varepsilon,W}(v_1+v_0))$ and $(\eta_2,z_2) = (g_{\varepsilon,X}(v_2+v_0),v_2+v_0 + g_{\varepsilon,W}(v_2+v_0))$. Now
	 			$
	 			G_{\varepsilon,d}(\eta_1,z_1) = M\cdot (\overline{w}_1,\overline{v}_1)
	 			$ and $
	 			G_{\varepsilon,d}(\eta_2,z_2) = M\cdot (\overline{w}_2,\overline{v}_2).
	 			$
	 			The function $G_{\varepsilon,d}:X\times \R^d \to \R^d$ will be treated as the function $G_{\varepsilon,d}:X\times V\times W \to \R^d$ and the corresponding derivatives will be denoted by $D_\eta G_{\varepsilon,d}$, $D_v G_{\varepsilon,d}$, and $D_w G_{\varepsilon,d}$. We calculate
	 			\begin{align*}
	 				& M\cdot (\overline{w}_2-\overline{w}_1,\overline{v}_2-\overline{v}_1) =  G_{\varepsilon,d}(\eta_2,z_2) - G_{\varepsilon,d}(\eta_1,z_1) \\
	 				& \ \in
	 				\overline{D_v G_{\varepsilon,d}([(\eta_1,z_1),(\eta_2,z_2)])}(v_2-v_1)\\
	 				& \qquad \qquad + \overline{D_w G_{\varepsilon,d}([(\eta_1,z_1),(\eta_2,z_2)])}(g_{\varepsilon,W}(v_2+v_0)-g_{\varepsilon,W}(v_1+v_0))\\
	 				& \qquad\qquad  + \overline{D_\eta G_{\varepsilon,d}([(\eta_1,z_1),(\eta_2,z_2)])}(g_{\varepsilon,X}(v_2+v_0)-g_{\varepsilon,X}(v_1+v_0))\\
	 				& \ \subset
\overline{D_v G_{\varepsilon,d}([(\eta_1,z_1),(\eta_2,z_2)])}(v_2-v_1) + \overline{D_w G_{\varepsilon,d}([(\eta_1,z_1),(\eta_2,z_2)])}\ \overline{D g_{\varepsilon,W}([v_2+v_0,v_1+v_0])} (v_2-v_1)\\
	 				& \qquad \qquad  + \overline{D_\eta G_{\varepsilon,d}([(\eta_1,z_1),(\eta_2,z_2)])}\ \overline{D g_{\varepsilon,X}([v_2+v_0,v_1+v_0])} (v_2-v_1).
	 			\end{align*}
	 			Furthermore, treating $G:\R^d\times \R^d$ as the function $G:X\times \R^d\times \R^d$, independent on the variable in $X$, we obtain  
	 			\begin{align*}
	 				& M\cdot (\overline{w}_2-\overline{w}_1,\overline{v}_2-\overline{v}_1) \in D_v G(z_0)\left(Dg(v_0)(v_2-v_1)+v_2-v_1\right)\\
	 				& \qquad +
	 				\overline{(D_vG_{\varepsilon,d}(\cdot)-D_vG(\cdot))([(\eta_1,z_1),(\eta_2,z_2)])}(v_2-v_1) + \overline{D_vG([z_1,z_2])-D_vG(z_0)}(v_2-v_1)\\
	 				& \qquad + \overline{(D_w G_{\varepsilon,d}(\cdot)-D_w G (\cdot))([(\eta_1,z_1),(\eta_2,z_2)])}\ \overline{D g_{\varepsilon,W}([v_2+v_0,v_1+v_0])}(v_2-v_1)\\
	 				& \qquad + \left(\overline{D_w G ([z_1,z_2])}-D_wG(z_0)\right)\overline{D g_{\varepsilon,W}([v_2+v_0,v_1+v_0])}(v_2-v_1)\\
	 				& \qquad + D_wG(z_0)\overline{(D g_{\varepsilon,W}(\cdot)-Dg(\cdot))([v_2+v_0,v_1+v_0])}(v_2-v_1)\\
	 				& \qquad +D_wG(z_0)\left(\overline{Dg([v_2+v_0,v_1+v_0])}-Dg(v_0)\right)(v_2-v_1)\\
	 				& \qquad  + \overline{D_\eta G_{\varepsilon,d}([(\eta_1,z_1),(\eta_2,z_2)])}\ \overline{D g_{\varepsilon,X}([v_2+v_0,v_1+v_0])} (v_2-v_1) \\
	 				& \ \ = I + II+III+IV+V+VI+VII+VIII.
	 			\end{align*}
	 			The  term $I$ in the last sum belongs to the tangent space $T_{G(z_0)}G\circ(Id_{B_1}+g)$. Terms $II, IV, VI$ satisfy the estimate $|II+IV+VI|\leq C_2(\varepsilon)|v_2-v_1|$ with $C_2(\varepsilon)\to 0$ as $\varepsilon\to 0$. Terms $III,V,VII$ satisfy the bound $|III+V+VII|\leq C_3(\delta_3)|v_2-v_1|$, where $C_3(\delta_3)$ also tends to zero as $\delta_3\to 0$. As for the term $VIII$ we have the bound $|VIII|\leq C_4 E(\varepsilon)|v_2-v_1|$, where $C_4$ is a bound on $\|D_\eta G_{\varepsilon,d}(\eta,u)\|_{\mathcal{L}(X;\R^d)}$ for $\eta,u\in P\times (B_1+B_2)$.
	 			It follows that
	 			$$
	 			|(0,\overline{v}_2-\overline{v}_1) | \geq C_1 |v_2-v_1| - (C_2(\varepsilon)+C_3(\delta_3)+C_4E(\varepsilon))|v_2-v_1|,
	 			$$
	 			and
	 			$$
	 			|(\overline{w}_2-\overline{w}_1,0) | \leq  (C_2(\varepsilon)+C_3(\delta_3)+C_4E(\varepsilon))|v_2-v_1|.
	 			$$
 Hence
	 			$$
	 			|(\overline{w}_2-\overline{w}_1,0) | \leq  \frac{C_2(\varepsilon)+C_3(\delta_3)+C_4E(\varepsilon)}{C_1-C_2(\varepsilon)-C_3(\delta_3)-C_4E(\varepsilon)}|(0,\overline{v}_2-\overline{v}_1) |,
	 			$$
	 			and we have the required Lipschitz condition with the arbitrarily small constant obtained by decreasing $\varepsilon$ and the radius $\delta_3$.
	 			Finally, let us estimate
	 			\begin{align*}
	 				& G_{\varepsilon,X}(\eta_2,z_2) - G_{\varepsilon,X}(\eta_1,z_1) \\
	 				& \ \in
	 				\overline{D_v G_{\varepsilon,X} ([(\eta_2,z_2),(\eta_1,z_1)])}(v_2-v_1)\\
	 				& \qquad \qquad + \overline{D_w G_{\varepsilon,X} ([(\eta_2,z_2),(\eta_1,z_1)])}(g_{\varepsilon,W}(v_2+v_0)-g_{\varepsilon,W}(v_1+v_0)) \\
	 				& \qquad \qquad  + \overline{D_\eta G_{\varepsilon,X} ([(\eta_2,z_2),(\eta_1,z_1)])}(g_{\varepsilon,X}(v_2+v_0)-g_{\varepsilon,X}(v_1+v_0))\\
	 				& \ \subset 
	 								\overline{D_v G_{\varepsilon,X} ([(\eta_2,z_2),(\eta_1,z_1)])}(v_2-v_1)\\
	 								& \qquad \qquad + \overline{D_w G_{\varepsilon,X} ([(\eta_2,z_2),(\eta_1,z_1)])}\ \overline{D g_{\varepsilon,W}([v_2+v_0,v_1+v_0])}(v_2-v_1)\\
	 				& \qquad \qquad  + \overline{D_\eta G_{\varepsilon,X} ([(\eta_2,z_2),(\eta_1,z_1)])}\ \overline{D g_{\varepsilon,X}([v_2+v_0,v_1+v_0])}(x_2-x_1).
	 			\end{align*}
	 			This means that
	 			$$
	 			\|G_{\varepsilon,X}(\eta_2,z_2) - G_{\varepsilon,X}(\eta_1,z_1)\|_X \leq L_1|v_2-v_1| \leq  L_2|(0,\overline{v}_2-\overline{v}_1) |,
	 			$$
	 			which completes the proof.
	 			
	 			\end{proof}
	 	
\subsection{Preparation of good coordinate system and representation of manifolds for $\varepsilon=0$.}

In this subsection we come back to the considered problem, and we work with system \eqref{eq:ODE}. We recall that in that system 
there exists a finite number of equilibria  $\{ e_1,\ldots, e_N\}$ which are all hyperbolic, that is $Df(e_i) = D^2F(e_i)$ is a nonsingular matrix for every $e_i$.

Using the Hadamard--Perron theorem, cf. for example \cite[Theorem 3.2.1]{Guckenheimer}, and the fact that $f\in C^2(\mathbb{R}^d;\mathbb{R}^d)$ we deduce that each equilibrium has the stable and unstable manifold $W^u(e_i)$ and $W^s(e_i)$ which is of class $C^2$.  If, for the two equilibria $e_i$ and $e_j$, there exists the solution $\gamma$ that connects $e_i$ to $e_j$ then this solution
belongs to both the unstable manifold of $e_i$ and the stable manifold of $e_j$.

Assume  that for system \eqref{eq:ODE} equilibrium $e_i$ is connected to $e_j$. The goal of this subsection is to prepare the coordinates, in which  $W^u(e_i)$ and $W^s(e_j)$
have an appropriate local representation near the point of intersection, such that it will be possible to establish later the transversal intersection of these manifolds for $\varepsilon >0$.

In Section \ref{sec:LapCC} we have verified the cone and isolation conditions which, by results of Appendix 2 guarantee the existence of the local stable and unstable manifolds of all equilibria in
 their neighborhoods.  These neighborhoods are isolating h-sets with cones. So, if the equilibrium $e_i$ is connected to $e_j$, that is,  if $W^u(e_i) \cap W^s(e_j) \neq \emptyset$, we can take $z \in W^s_{loc,N(e_j)}(e_j)$ as an intersection point,  where $N(e_j)$ is an h-set with cones for the equilibrium $e_j$. Then  $z \in W^u(e_i)$.

We assume that the intersection of $W^u(e_i)$ and $W^s_{loc,N(e_j)}(e_j)$ is transversal, that is the algebraic sum of the both tangent spaces constitutes the whole $\R^d$ 
$$T_z W^u(e_i) + T_z W^s_{loc,N(e_j)}(e_j) = \mathbb{R}^d.$$

Let $\textrm{dim} W^u(e_i) = u_i$ and $\textrm{dim} W^s_{loc,N(e_j)}(e_j) = s_j$. Then $u_i + s_j \geq d+1$.
Tangent space $T_z W^u(e_i)$ is the $u_i$ dimensional subspace of $\mathbb{R}^d$ and $T_z W^s_{loc,N(e_j)}(e_j)$ is its $s_j$ dimensional subspace. The intersection of both spaces is $c$ dimensional subspace of $\mathbb{R}^d$ where $c = u_i+s_j - d$. Denote $u_i = k_1+c$ and $s_j = k_2+c$. There exists an invertible $d\times d$ matrix $M$ such that $M\cdot (\mathbb{R}^{k_1}\times \mathbb{R}^{c}\times (0)_{k_2}) = T_z W^u(e_i)$ and $M\cdot ((0)_{k_1} \times \mathbb{R}^{c}\times \mathbb{R}^{k_2}) = T_z W^s_{loc,N(e_j)}(e_j)$.
This matrix defines the linear change of coordinates in $\mathbb{R}^d$. We will denote the new coordinates by $(\overline{x},\overline{a},\overline{y}) \in \mathbb{R}^{k_1}\times \mathbb{R}^c \times \mathbb{R}^{k_2}$.
\begin{figure}
	\centering
	\includegraphics[width=0.7\linewidth]{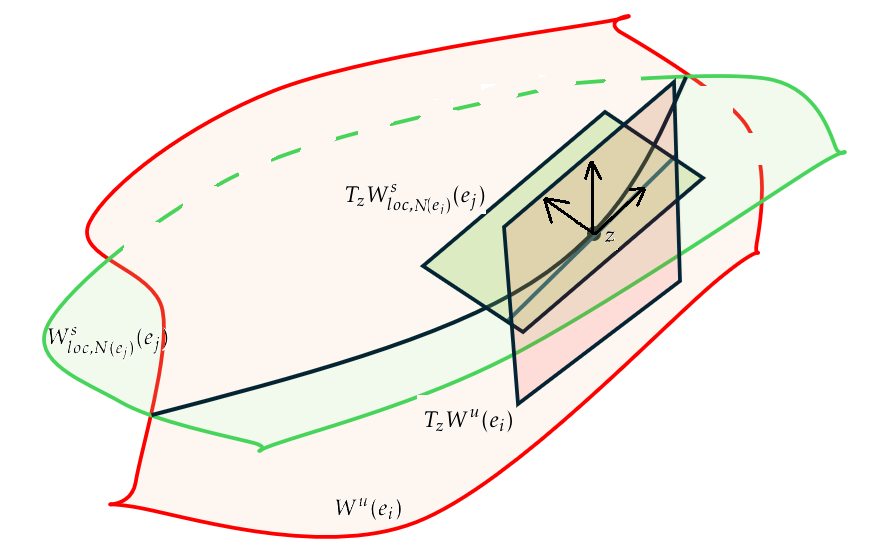}
	\caption{Illustration of the linear transformation of the system of coordinates. In the picture both $W^s(e_j)$ and $W^u(e_i)$ are two dimensional and the intersection is one dimensional. Then, three arrows represent the new coordinate system with $a=1$, $k_1=1$, and $k_2=1$.}
	\label{fig:coordinates}
\end{figure}

\begin{figure}
	\centering
	\includegraphics[width=0.7\linewidth]{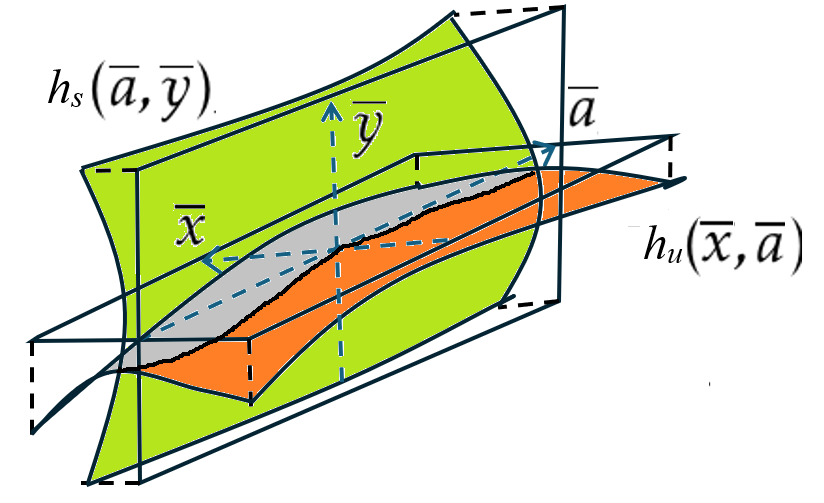}
	\caption{Illustration of Lemma \ref{lem:box}. Local stable and unstable manifold are, respectively, vertical and horizontal disks over the tangent spaces contained in a box centered in the point of intersection.}
	\label{fig:local}
\end{figure}

Both manifolds $W^u(e_i)$ and $W^s(e_j)$ are $C^1$ disks, which are tangent to $M\cdot   (  \mathbb{R}^{k_1}\times \mathbb{R}^c \times (0)_{k_2})$ and $M\cdot   ( (0)_{k_2}\times \mathbb{R}^c \times \mathbb{R}^{k_2})$, respectively. We obtain the following lemma 
\begin{lem}\label{lem:box}
	There exist constants $\delta_{k_1}, \delta_{k_2}, \delta_c > 0$, balls $B(0,\delta_c)\subset \R^c$, $ B(0,\delta_{k_1})\subset \R^{k_1}$, $ B(0,\delta_{k_2})\subset \R^{k_2}$  and Lipschitz functions $h^s:B(0,\delta_c)\times B(0,\delta_{k_2}) \to B(0,\delta_{k_1})$ and $h^u: B(0,\delta_{k_1})\times B(0,c) \to B(0,\delta_{k_2})$ such that
	\begin{equation}\label{eq:disk_1}
	z+ M \cdot \textrm{graph}(h^s)\subset  W^s(e_j),
	\end{equation}
	and
	\begin{equation}\label{eq:disk_2}
	z+M \cdot \textrm{graph}(h^u) \subset W^u(e_i).
	\end{equation}
	Moreover, the radii $\delta_{k_1}, \delta_{k_2}, \delta_c$ can be decreased to make the Lipschitz constants of both disks $h^s, h^u$ arbitrarily small.
\end{lem}
\begin{proof}
	Since we know that local stable manifold of $e_j$ and local unstable manifold of $e_i$ are $C^1$ disks, the lemma is a consequence of Lemma \ref{lemma_disks} (with $\varepsilon = 0$ and $m=c+k_2$) and Lemma \ref{lem:transport_disks}  (with $n=k_1+c$ and $G=S^0(t)$ where $t$ is such that $z=S^0(t)z_0$ for a point $z_0$ in a local unstable manifold of $e_i$).  We only need to guarantee the choice of $\delta_{k_1}, \delta_{k_2}, \delta_c > 0$ for which $\text{range}(h^u)\subset  B(0,\delta_{k_2})$ and $\text{range}(h^s) \subset B(0,\delta_{k_1})$. 
	In Lemmas \ref{lemma_disks} and \ref{lem:transport_disks} we construct the disks contained in the graphs of these manifolds: $h_s:B(0,\delta_1)\to B(0,\delta_2)$ with $\delta_2 \leq c_1(\delta_1)\delta_1$ and $h_u:B(0,\delta_3)\to B(0,\delta_3)$ with $\delta_4 \leq c_2(\delta_3)\delta_3$, where $c_1, c_2$ are some functions which tend to zero as the argument tends to zero. We decrease $\delta_1$ and $\delta_3$  to have $\delta_1=\delta_3=\delta$,  $c_1(\delta)\leq \frac{1}{4}$ and $c_2(\delta)\leq \frac{1}{4}$. Now we can pick $\delta_{k_1} = \delta_{k_2} = \delta_c = \frac{\delta}{2}$. 
	For $(\overline{x},\overline{a}) \in B(0,\delta_{k_1})\times B(0,c)$ we have $|\overline{x}|^2+|\overline{a}|^2\leq \frac{\delta^2}{2} < \delta^2$, so 
	$|h^u(\overline{x},\overline{a})|\leq c_1(\delta)\delta \leq \frac{\delta}{4} < \frac{\delta}{2} = \delta_{k_1}$, and, likewise, for $(\overline{a},\overline{y}) \in B(0,c)\times B(0,\delta_{k_2})$ we have $|\overline{a}|^2+|\overline{y}|^2\leq \frac{\delta^2}{2} < \delta^2$, so 
	$|h^s(\overline{a},\overline{y})|\leq c_2(\delta)\delta \leq \frac{\delta}{4} < \frac{\delta}{2} = \delta_{k_2}$. The proof is complete. 
	\end{proof}
Note that above argument wstays valid if we scale down the size of the box by decreasing $\delta$.

\subsection{Intersection of invariant manifolds for $\varepsilon >0$}

Now we consider the problem with $\varepsilon>0$. Results of Appendix 4 and Appendix 5 guarantee that local stable and unstable manifolds of all equilibria for $\varepsilon>0$ are $C^1$ disks, and dependence of their derivatives on parameter $\varepsilon$ is established in Appendix 6. These results are needed to guarantee that assumptions of Lemmas \ref{lem:stable} and \ref{lem:transport_eps} hold. The local stable manifold of the equilibrium $(0,e_j^\varepsilon)$ is a $C^1$ disk over variables $(\eta,v) \in \overline{B}_{L^2_A( \R^+)^d}(0,R) \times B_1$, where $B_1$ is a $k_2+c$ dimensional ball, and the values are in $\R^{k_1}$. Likewise, the local unstable manifold of $(0,e_i^\varepsilon)$ is the disk over the variable $y \in B_2$, where $B_2$ is a $k_1+c$ dimensional ball, and the values are in $L^2_A( \R^+)^d \times \mathbb{R}^{k_2}$. In the next lemma we prove that for every $\eta \in B_{L^2_A( \R^+)^d}(0,R)$ the section of the stable manifold of $(0,e_j^\varepsilon)$   and the image by $S^\varepsilon(t)$ of the local unstable manifold of $(0,e_i^\varepsilon)$ are both disks in the box constructed in Lemma   \ref{lem:box}. We also  calculate the Lipschitz constants of these disks.
\begin{lem}\label{lem:boxeps}
Consider the box $B(0,\delta_{k_1}) \times B(0,\delta_c)\times B(0,\delta_{k_2})$ from Lemma \ref{lem:box}. There exists $\varepsilon_0>0$ such that for every $\varepsilon \in (0,\varepsilon_0]$, if  $(0,e_j^\varepsilon)$ is an equilibrium for $\varepsilon$ that corresponds to $e^j$, and  $\overline{B}_{L^2_A( \R^+)^d}(0,R)\subset L^2_A( \R^+)^d$ is a ball such that local stable manifold of $(0,e_j^\varepsilon)$ is defined as a disk with $\eta \in \overline{B}_{L^2_A( \R^+)^d}(0,R)$ then there exists a disk 
$$
h^s_{\varepsilon}:\overline{B}_{L^2_A( \R^+)^d}(0,R) \times B(0,\delta_c)\times B(0,\delta_{k_2}) \to  B(0,\delta_{k_1}) 
$$
in the stable manifold $W^s((0,e_j^\varepsilon))$ translated to $z$ with finite dimensional variables transformed by $M$. Moreover the following Lipschitz condition holds
$$
|h_{s,\varepsilon}(\eta_2,\overline{a}_2,\overline{y}_2) - h_{s,\varepsilon}(\eta_1,\overline{a}_1,\overline{y}_1)|\leq D_1 |(0,\overline{a}_1-\overline{a}_2,\overline{y}_1-\overline{y}_2) | + D_2E\|\eta_1-\eta_2\|_{L^2_A( \R^+)^d},
$$
  for every $\eta_1, \eta_2 \in B_{L^2_A( \R^+)^d}(0,R)$, $\overline{a}_1, \overline{a}_2 \in B(0,\delta_c)$, $\overline{y}_1, \overline{y}_2 \in B(0,\delta_{k_2})$ where the constant $D_1$ can be made arbitrarily small by decreasing the radii $\delta_1, c$  and $\varepsilon$,  $E$ is the Lipschitz constant of the $\eta$ variable in the local unstable manifold of  $(0,e_j^\varepsilon)$, and $D_2>0$ is a constant.

  Moreover, there exists a disk $h^u_{\varepsilon}=(h^u_{\varepsilon,L^2_A(\R^+)^d},h^u_{\varepsilon,k_2}):B(0,\delta_{k_1}) \times B(0,\delta_c) \to L^2_A(\R^+)^d \times B(0,\delta_{k_2}) $ in the unstable manifold $W^u((0,e_i^\varepsilon))$ translated to the point $z$ with finite dimensional variables transformed by $M$,  satisfying
  $$
  |h^u_{\varepsilon,k_2}(\overline{x}_2,\overline{a}_2) - h^u_{\varepsilon,k_2}(\overline{x}_1,\overline{a}_1)|\leq D_3 |(\overline{x}_2-\overline{x}_1,\overline{a}_2-\overline{a}_1,0)|,
  $$
  where the constant $D_3$ can be made arbitrarily small by decreasing $\varepsilon$ and the size of the box. The function $h^u_{\varepsilon,L^2_A(\R^+)^d}$ satisfies the Lipschitz condition
  $$
  \|h^u_{\varepsilon,L^2_A(\R^+)^d}(\overline{x}_2,\overline{a}_2) - h^u_{\varepsilon,L^2_A(\R^+)^d}(\overline{x}_1,\overline{a}_1)\|_{L^2_A( \R^+)^d}\leq D_4 |(\overline{x}_2-\overline{x}_1,\overline{a}_2-\overline{a}_1,0)|,
  $$
  for some constant $D_4>0$.
\end{lem}			
\begin{proof}
	Existence of $h^s_\varepsilon$ and its Lipschitz condition follows by Lemma \ref{lem:stable} and existence of $h^u_{\varepsilon}$ and its Lipschitz condition follows by Lemma \ref{lem:transport_eps} with $G_\varepsilon = S^\varepsilon(t)$ for a certain $t$. As in Lemma \ref{lem:box} we need to show that it is possible to choose $\delta_{k_1}, \delta_{k_2}, \delta_c > 0$ for which $\text{range}(h^u_{\varepsilon,k_2})\subset  B(0,\delta_{k_2})$ and $\text{range}(h^s) \subset B(0,\delta_{k_1})$. The proof is similar as of Lemma \ref{lem:box} with additional necessity to take into account $\varepsilon$. 
	
	The disks constructed in  Lemmas \ref{lem:stable} and \ref{lem:transport_eps} satisfy: $h^s_\varepsilon:B_{L^2_A( \R^+)^d}(0,R)\times B(0,\delta_1)\to B(0,\delta_2)$ with $\delta_2 \leq c_1(\varepsilon)+c_2(\delta_1)\delta_1$ and $h_u:B(0,\delta_3)\to B(0,\delta_3)\times L^2_A(\R^+)^d$ with $\delta_4 \leq c_3(\varepsilon) + c_4(\delta_3)\delta_3$, where $c_1, c_2, c_3, c_4$ are functions which tend to zero as the argument tends to zero. 
	
	As in the proof of Lemma \ref{lem:box} we decrease $\delta_1$ and $\delta_3$ to have $\delta_1=\delta_3=\delta$,  $c_2(\delta)\leq \frac{1}{4}$ and $c_4(\delta)\leq \frac{1}{4}$. We next decrease $\varepsilon$ to get $c_1(\varepsilon) \leq \frac{\delta}{4}$ and $c_3(\varepsilon) \leq \frac{\delta}{4}$. Now we can pick $\delta_{k_1} = \delta_{k_2} = \delta_c = \frac{\delta}{2}$. 
	For $(\overline{x},\overline{a}) \in B(0,\delta_{k_1})\times B(0,c)$ we have $|\overline{x}|^2+|\overline{a}|^2\leq \frac{\delta^2}{2} < \delta^2$, so 
	$|h^u_{\varepsilon,k_2}(\overline{x},\overline{a})|\leq c_1(\varepsilon)+c_2(\delta)\delta \leq \frac{\delta}{2} = \delta_{k_1}$, and, likewise, for $(\overline{a},\overline{y}) \in B(0,c)\times B(0,\delta_{k_2})$ we have $|\overline{a}|^2+|\overline{y}|^2\leq \frac{\delta^2}{2} < \delta^2$, so 
	$|h_s(\eta,\overline{a},\overline{y})|\leq c_3(\varepsilon)+c_4(\delta)\delta \leq \frac{\delta}{2} = \delta_{k_2}$. The proof is complete. 
	\end{proof}	 			
We pass to the proof that the manifolds for $\varepsilon>0$ intersect. In the constructed box, we must find the intersection of the unstable manifold of $(0,e^\varepsilon_i)$ with the stable manifold of $(0,e^\varepsilon_j)$. We will consider the section of the box at $\overline{a}=0$, thus fixing to zero the variable in the space that is tangent to the manifold intersection.  
We first consider the mapping
$$
B(0,\delta_{k_1})\ni \overline{x} \mapsto (h^u_{\varepsilon,L^2_A(\R^+)^d}(\overline{x},0),h^u_{\varepsilon,k_2}(\overline{x},0))\in L^2_A( \R^+)^d\times B(0,\delta_{k_2}).
$$
We want to compose it with the mapping
$$
 \overline{B}_{L^2_A( \R^+)^d}(0,R)\times B(0,\delta_{k_1}) \ni (\eta,\overline{y}) \mapsto h^s_\varepsilon(\eta,0,\overline{y}) \in B(0,\delta_{k_1}).
$$
If we are able to prove that the composition of the above mappings has a fixed point then this fixed point corresponds to the intersection point of the manifolds for $\varepsilon>0$. The first thing we must show is that $h^u_{\varepsilon,L^2_A(\R^+)^d}(\overline{x},0) \in \overline{B}_{L^2_A( \R^+)^d}(0,R)$, which is required for the composition to make sense. This, together with the Banach fixed point argument is done in the next Lemma.  
\begin{lem}\label{lem:intersection}
	There exists $\varepsilon_0>0$ such that for every $\varepsilon\in [0,\varepsilon_0]$ the unstable manifold of $(0,e^\varepsilon_i)$ intersects with the stable manifold of $(0,e^\varepsilon_j)$.
\end{lem}
\begin{proof}
	Take $\overline{x}_1, \overline{x}_2 \in B(0,\delta_{k_1})$. We have
	$$
	|h^u_{\varepsilon,k_2}(\overline{x}_2,0) - h^u_{\varepsilon,k_2}(\overline{x}_1,0)|\leq D_3 |(\overline{x}_2-\overline{x}_1,0,0)|,
	$$
	and
	  $$
	\|h^u_{\varepsilon,L^2_A(\R^+)^d}(\overline{x}_2,0) - h^u_{\varepsilon,L^2_A(\R^+)^d}(\overline{x}_1,0)\|_{L^2_A( \R^+)^d}\leq D_4 |(\overline{x}_2-\overline{x}_1,0,0)|.
	$$
	We already know that the range of $h^u_{\varepsilon,k_2}$ is in the ball $B(0,\delta_{k_2})$. But we still need to show that the range of $h^u_{\varepsilon,L^2_A(\R^+)^d}$ is a subset of $\overline{B}_{L^2_A(\R^+)^d}(0,R)$. To this end, we remind that the point of intersection $z$ of the stable and unstable manifolds at $\varepsilon=0$ is given by $S^0(t)p$ for $p=x_p+y_p$ in the local unstable manifold of $e_i$ ($y_p$ is the function of $x_p$). Let us estimate the norm of the memory variable at $\varepsilon=0$ corresponding to this solution as $s\to \infty$. We have
	$$
	\eta^t(r) = S^0(t-r)p-S^0(t)p.
	$$
	The norm of this memory variable is given by
	\begin{align*}
	& \|\eta^t\|_{L^2_A( \R^+)^d}^2 = \int_0^\infty (A(r)(S^0(t-r)p-S^0(t)p),(S^0(t-r)p-S^0(t)p))\, dr \\
	& \qquad = \int_0^{r_0} (A(r)(S^0(t-r)p-S^0(t)p),(S^0(t-r)p-S^0(t)p))\, dr \\
	& \qquad \ \ + \int_{r_0}^\infty (A(r)(S^0(t-r)p-S^0(t)p),(S^0(t-r)p-S^0(t)p))\, dr.
	\end{align*}
	In the following calculation by $C$ we will denote a generic constant. Fix $\gamma>0$ and   $r_0$. For these values we can find $t$ such that for $r\in [0,r_0]$ we have $t-r \geq t-r_0 \geq t_0$ with $t_0$ sufficiently large to guarantee that $S^0(s)p \in B(e_j,\gamma)$ for $s\geq t_0$. We have
	$$
	\|\eta^t\|_{L^2_A( \R^+)^d}^2 \leq 4\gamma^2\int_0^{r_0}\|A(r)\|dr + C\int_{r_0}^\infty\|A(r)\|dr \leq C\gamma^2 + Ce^{-Cr_0}.
	$$
	whence
	$$
	\|\eta^t\|_{L^2_A( \R^+)^d} \leq C\gamma + Ce^{-Cr_0}.
	$$
	Now $h^u_{\varepsilon,L^2_A(\R^+)^d}(\overline{x},0) $ is an image by $S^\varepsilon_\eta(t)$ of a certain point $(\xi,x+y)\in L^2_A(\R^+)\times \R^d$ in a local unstable manifold of $(0,e^\varepsilon_i)$ such that the distance $|x-x_p|$ does not exceed $\varepsilon_1$ (where we can choose arbitrarily small $\varepsilon_1$ and this choice determines $\varepsilon$).
	From \eqref{eq:close_1} we deduce that
	$$
	\|h^u_{\varepsilon,L^2_A(\R^+)^d}(\overline{x},0) - \eta^t\|_{L^2_A( \R^+)^d} \leq Ce^{Ct}(\|\xi-\eta^0\|_{L^2_A( \R^+)^d} + |x-x_p| + |y-y_p| \varepsilon),
	$$
	where $\eta^0$ is the memory variable corresponding to the total solution passing through $p$ with $\varepsilon=0$. By Theorem \ref{thm:unstable} the last quantity can be estimated from above as follows
	$$
	\|h^u_{\varepsilon,L^2_A(\R^+)^d}(\overline{x},0) - \eta^t\|_{L^2_A( \R^+)^d} \leq  Ce^{Ct}(|x-x_p|+\varepsilon) \leq Ce^{Ct}(\varepsilon_1 +\varepsilon).
	$$
	Now we estimate $\|h^u_{\varepsilon,L^2_A(\R^+)^d}(\overline{x},0)\|_{L^2_A( \R^+)^d}$. We have
	$$
	\|h^u_{\varepsilon,L^2_A(\R^+)^d}(\overline{x},0)\|_{L^2_A( \R^+)^d} \leq \|h^u_{\varepsilon,L^2_A(\R^+)^d}(\overline{x},0) - \eta^t\|_{L^2_A( \R^+)^d} + \|\eta^t\|_{L^2_A( \R^+)^d} \leq  Ce^{Ct}\varepsilon_1 +Ce^{Ct}\varepsilon + C\gamma + Ce^{-Cr_0}.
	$$
	We need the last quantity to be less than $R$. We fix $r_0$ and $\gamma$ so that each of two two last terms is no larger than $\frac{R}{4}$. This forces us to choose $t$. Now choose $\varepsilon_1$ such that the first term is no larger than $\frac{R}{4}$. Finally, if necessary, decrease $\varepsilon$ so that the second term does not exceed $\frac{R}{4}$.
	We come back to the calculation of the Lipschitz constant for the fixed point mapping. We have
\begin{align*}
& |h^s_\varepsilon(h^u_{\varepsilon,L^2_A(\R^+)^d}(\overline{x}_2,0),0,h^u_{\varepsilon,k_2}(\overline{x}_2,0)) - h^s_\varepsilon(h^u_{\varepsilon,L^2_A(\R^+)^d}(\overline{x}_1,0),0,h^u_{\varepsilon,k_2}(\overline{x}_1,0))|\\
& \qquad \leq D_1 |h^u_{\varepsilon,k_2}(\overline{x}_2,0) - h^u_{\varepsilon,k_2}(\overline{x}_1,0) | + D_2E\|h^u_{\varepsilon,L^2_A(\R^+)^d}(\overline{x}_2,0) - h^u_{\varepsilon,L^2_A(\R^+)^d}(\overline{x}_1,0)\|_{L^2_A( \R^+)^d}\\
& \qquad  \leq (D_1 D_3 + D_2E D_4) |\overline{x}_2-\overline{x}_1|.
\end{align*}
The constants $D_1, D_3$ can be made arbitrarily small by decreasing $\varepsilon$ and scaling down  the size of the box. Moreover $E$ can be made arbitrarily small by decreasing $\varepsilon$ (cf. Lemma \ref{lem:cone:par}). We decrease these constants such that $D_1 D_3 + D_2E D_4 < 1$. Then the constructed mapping is a contraction and hence it has a fixed point which is the sought intersection of the manifolds.
\end{proof}
	\begin{figure}
	\centering
	\includegraphics[width=0.7\linewidth]{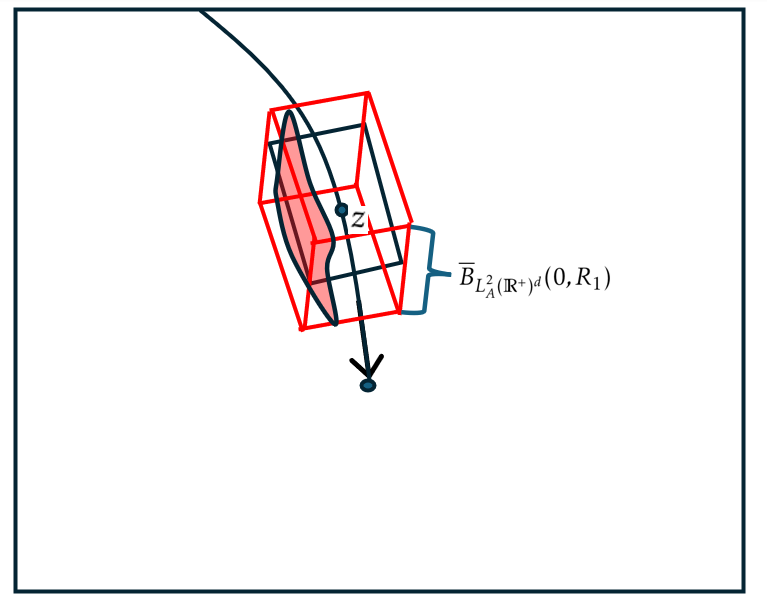}
	\caption{Illustration of the part of the proof of Lemma \ref{lem:boxeps} for the stable manifold. We prove that, if we take a Cartesian product of the box from Lemma \ref{lem:box} with a ball $\overline{B}_{L^2_A}(\R^+)^d$, then the intersection of the local stable manifold of $(0,e_i^\varepsilon)$ with this new box is a disk over $(\eta,\overline{a},\overline{y})$. The  variables $(\overline{a},\overline{y})$ are the same as in the box for $\varepsilon=0$ constructed in Lemma \ref{lem:box}, $\eta\in \overline{B}_{L^2_A}(\R^+)^d$, and the value $\overline{x}$ also stays in the same box as for $\varepsilon=0$. The Lipschitz constants both for $(\overline{a},\overline{y})$ and $\eta$ arguments can be made arbitrarily small.}
	\label{fig:stable}
\end{figure}

\begin{figure}
	\centering
	\includegraphics[width=0.7\linewidth]{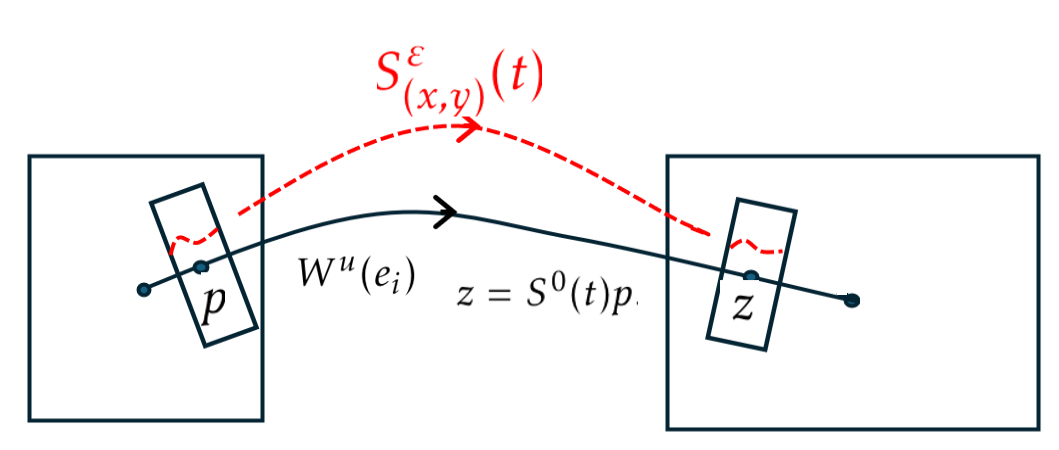}
	\caption{Illustration of the part of the proof of Lemma \ref{lem:boxeps} for the unstable manifold. We prove (using the Brouwer degree argument) that, if we take the image of the certain horizontal disk being the fragment of the local unstable manifold for $\varepsilon>0$ (red dashed line on the left) and transport it forwards, we obtain the horizontal disk $\overline{y}(
		\overline{a},\overline{x})$ in the local variables in the box centered around $z$ (red dashed line on the right) with arbitrarily small Lipschitz constant.  The memory variable $\eta$,  which is the second variable in the image of the transported disk, needs to satisfy the Lipschitz condition  with the constant not necessarily small (it is obtained in the final part of the proof of Lemma \ref{lem:transport_eps}) and needs to belong to the ball $\overline{B}_{L^2_A}( \R^+)^d (0,R)$ (it is obtained in Lemma \ref{lem:intersection}).}
	\label{fig:unstable}
\end{figure}
As a consequence of the above result, we have the following theorem.
\begin{thm}\label{thm:lower}
	There exists $\varepsilon_0>0$ such that  if there exists the connection between the equilibria $e_i$ and $e_j$ via the ODE \eqref{eq:ODE} then for every $\varepsilon\in (0,\varepsilon_0]$ there exists the connection between the corresponding equilibria $(0,e_i^\varepsilon)$ and $(0,e_j^\varepsilon)$ via the system \eqref{main2}--\eqref{eq:main}.
\end{thm}

% appendix

\section{Appendix 1: Asymptotic compactness}

We consider  the semigroup
$$
S^\varepsilon(t):L^2_A(\R^+)^d\times \R^d   \to L^2_A(\R^+)^d\times \R^d,
$$
given by the solutions of \eqref{eq:main}--\eqref{main2}, namely
as $S^\varepsilon(t)(\eta^0,x_0) = (\eta^t,x(t))$. We do not expect the compactness of $S^\varepsilon(t)$ for a finite $t$. Instead we prove the following lemma on asymptotic compactness
\begin{lem}\label{lemma:compact}
Assume that we have the estimate $|x(t)|\leq C(|x_0|, \|\eta_0\|)$ for the function $C$ nondecreasing with respect to both  arguments (this a priori estimate follows from the Lyapunov function in Lemma \ref{lem:LapFunc1}).
	Assume that $\{\eta^{0,n},x_0^n\}$ is a sequence of initial data bounded in $L^2_A(\R^+)^d\times \R^d$ and $t_n\to \infty$. Then $S^\varepsilon(t_n)( \eta^{0,n},x_0^n)$ is relatively compact.
	\end{lem}
\begin{proof}
Denote $(\eta^{t,n},x^n(t)) =  S^\varepsilon(t)(\eta^{0,n},x_0^n)$. Then $|x^n(t_n)|$ is bounded, so it has a convergent subsequence. We denote this subsequence by the same index $n$, without renumbering. Then $x^n(t_n)\to \xi$ in $\R^d$. We need to show the relative compactness of $\eta^{{t_n},n}$, that is of
$
x^n(t_n-s) - x^n(t_n).
$ in $L^2_A(\R^+)^d$. Observe that
$$
0 \leq \int_0^\infty (A(s)(x^n(t_n)-\xi),x^n(t_n)-\xi)\, ds\leq \int_0^\infty \|A(s)\|\, ds |x^n(t_n)-\xi|^2 \to 0.
$$
It is enough to prove the relative compactness of $[0,\infty) \ni s \to x^n(t_n-s)\in \R^d$ in the space  $L^2_A(\R^+)^d$. We first demonstrate the relative compactness on $L^2_A(0,T)^d$ for every $T$. Note that the continuity and positive definiteness of $[0,T]\ni s \mapsto A(s)$ implies that the norms $L^2_A(0,T)^d$ and $L^2(0,T)^d$ are equivalent. We are in position to use the Kolmogorov-Riesz-Frechet theorem \textbf{mozna cos zaczytowac}  which states that the set $B\subset L^2(0,T)^d$ is relatively compact if and only if it is bounded in that space and
$$
\lim_{h\to 0}\sup_{u\in B}\int_0^{T-h}|u(s+h)-u(s)|^2\, ds = 0.
$$

In our case we need to show that
$$
\lim_{h\to 0}\sup_{n}\int_0^{T-h}|x^n(t_n-s-h)-x^n(t_n-s)|^2\, ds = 0,
$$
or
$$
\lim_{h\to 0}\sup_{n}\int_0^{T-h}\int_{t_n-s-h}^{t_n-s}|(x^n)'(r)|^2\, dr\, ds = 0.
$$
It is enough that the result is obtained for $n\geq n_0$ where $n_0$ may depend on $T$. We chose $n_0$ sufficiently large such that $t_n\geq T$. Then
$$
\int_0^{T-h}\int_{t_n-s-h}^{t_n-s}|(x^n)'(r)|^2\, dr\, ds \leq h \int_{t_n-T}^{t_n}|(x^n)'(r)|^2 \, dr.
$$
But, cf. \eqref{eq:main} and
Lemma \ref{lem:LapFunc1},
\begin{align*}
& \int_{t_n-T}^{t_n}|(x^n)'(r)|^2 \, dr\leq 3\int_{t_n-T}^{t_n}|f(x^n(s))|^2\, ds + \varepsilon C\int_{t_n-T}^{t_n}|x^n(s)|^2\, ds+  \varepsilon C \int_{t_n-T}^{T_n}\|\eta^{s,n}\|^2\, ds \\
&\ \  \leq T C(\|\eta^{0,n}\|,|x_0|),
 \end{align*}
 and the assertion follows. By the diagonal argument we can construct a subsequence, still denoted by $n$, which converges in $L^2_A(0,T)^d$ for every $T$. We denote the limit by $\eta$. We also have that $\eta^{t_n,n}(r)\to \eta(r)$ for almost every $r\geq 0$. As $\eta^{t_n,n}(r) = x^n(t_n-r)-x^n(r)$, it follows that $|\eta^{t_n,n}(r)|\leq E$ for a constant $E>0$ and $r\in [0,t_n]$, whence $|\eta(r)|\leq E$ for a.e. $r \geq 0$. We claim that this subsequence actually converges in $L^2_A(\R^+)^d$. To get this assertion we need to show that for every $\delta>0$ there exists $n_\delta$ such that
 for every $n\geq n_\delta$ there holds
 $$
\|\eta^{t_n,n}-\eta\|= \int_0^\infty (A(r)(\eta^{t_n,n}(r)-\eta(r)),(\eta^{t_n,n}(r)-\eta(r)))\, dr \leq \delta
 $$
 We choose $T_\delta$ such that
 $$
 \int_{T_\delta}^\infty \|A(r)\|\, dr \leq \frac{\delta}{12 E^2}.
 $$
 For large $T_\delta$ we split the integral into three parts
 \begin{align*}
 & \int_0^\infty (A(r)(\eta^{t_n,n}(r)-\eta(r)),\eta^{t_n,n}(r)-\eta(r))\, dr \\
 & = \int_0^{T_\delta} (A(r)(\eta^{t_n,n}(r)-\eta(r)),\eta^{t_n,n}(r)-\eta(r))\, dr \\
 & + \int_{T_\delta}^{t_n}(A(r)(\eta^{t_n,n}(r)-\eta(r)),\eta^{t_n,n}(r)-\eta(r))\, dr + \int_{t_n}^\infty (A(r)(\eta^{t_n,n}(r)-\eta(r)),\eta^{t_n,n}(r)-\eta(r))\, dr.
  \end{align*}
  From convergence in $L^2_A(0,T)^d$ for every $T$ it follows that we can find $n_\delta$ such that the first integral in no greater than $\delta/3$. Norm of the second integral is majorized as follows
  $$
  \left|\int_{T_\delta}^{t_n}(A(r)(\eta^{t_n,n}(r)-\eta(r)),\eta^{t_n,n}(r)-\eta(r))\, dr\right| \leq 4E^2 \int_{T_\delta}^{t_n} \|A(r)\|\, dr \leq4 E^2\int_{T_\delta}^\infty \|A(r)\|\, dr \leq \frac{\delta}{3}.
  $$
    To deal with the last integral let us compute    \begin{align*}
    & \int_{t_n}^\infty (A(r)(\eta^{t_n,n}(r)-\eta(r)),\eta^{t_n,n}(r)-\eta(r))\, dr \\
    & \qquad \qquad \leq 2  \int_{t_n}^\infty (A(r)\eta^{t_n,n}(r),\eta^{t_n,n}(r))\, dr + 2\int_{t_n}^\infty (A(r)\eta(r),\eta(r))\, dr \\
    & \qquad \qquad \leq 2  \int_{0}^\infty (A(r+t_n)\eta^{0,n}(r),\eta^{0,n}(r))\, dr + 2E^2\int_{t_n}^\infty \|A(r)\|\, dr .
    \end{align*}
    From \ref{asm1} we deduce
    \begin{align*}
    & \int_{t_n}^\infty (A(r)(\eta^{t_n,n}(r)-\eta(r)),\eta^{t_n,n}(r)-\eta(r))\, dr \\
    & \qquad \qquad \leq 2 e^{-Ct_n} \int_{0}^\infty (A(r)\eta^{0,n}(r),\eta^{0,n}(r))\, dr + \frac{\delta}{6} = 2 e^{-Ct_n} \|\eta^{0,n}\|^2 + \frac{\delta}{6}.
    \end{align*}
    We can find $n_\delta$ large enough, that the right-hand side of the last bound is no greater than $
    \delta/3$ and the proof is complete.
\end{proof}

\section{Appendix 2: Graph transform for existence of local stable and unstable manifolds.}
We begin with the definitions and properties of h-sets, isolating blocks, and cone conditions adapted for the problems with the distributed delay.
\begin{df}
	A family of  mappings $\{ S(t) \}_{t\geq 0}$ will be called a $C^0$ semiflow on $X$ if
	\begin{itemize}
		\item $[0,\infty) \times X \ni (t,x) \to S(t)x$ is continuous,
		\item $S(0) = I_X$, the identity,
		\item $S(t+s)x = S(t)(S(s)x)$ for every $s,t\geq 0$ and $x\in X$
	\end{itemize}
\end{df}	
\begin{df}
	A $C^0$ semiflow on $X$ $\{ S(t) \}_{t\geq 0}$ is asymptotically compact if
	for a bounded sequence $\{ x_n \} \subset  X$ and a sequence $t_n\to \infty$ the sequence  $S(t_n)x_n$ is relatively compact.
\end{df}
\begin{df} \label{df:wB}
	For a bounded set $B$ we define its $\omega$-limit set as
	$$
	\omega(B) = \{ x\in X\, :\ x = \lim_{n\to\infty} S(t_n)x_n\ \ \textrm{for sequences}\ \ t_n\to \infty \ \ \textrm{and} \ \ \{ x_n \}\subset B\}.
	$$
\end{df}
The following result is well known.
\begin{lem}\label{as:com}
	If a a $C^0$ semiflow is assymptotically compact, then for every nonempty bounded set $B\subset X$, the set $\omega(B)$ is nonempty, compact, connected, invariant, and
	$$
	\lim_{t\to\infty}\mathrm{dist}\, (S(t)B, \omega(B)) = 0.
	$$
\end{lem}

\begin{df}\label{h-set}
	Let $X$ be a Banach space.
	The set $A \subset X$ is called an \textit{h-set (hyperbolic set)}  if there exist the linear closed subspaces $X_1, X_2$ of $X$ with $X = X_1 \oplus X_2$ and $\textrm{dim}\, X_1  < \infty$, $\textrm{dim}\, X_1 = s + u$, with $s, u \in \mathbb{N}$, $u=u_1+2u_2$ and $s=s_1+2s_2$, the numbers $\{a_k\}_{k=1}^{s_1+s_2+u_1+u_2}$ with $a_k>0$ and an affine bijective mapping $L:\R^{\textrm{dim}\, X_1}\to X_1$ such that
	$$
	A = L \left( N_u \times N_s\right) \oplus \overline{B}_{X_2}(0,r),
	$$
	where
	$$
	N_u = \prod_{k=1}^{u_1}[-a_k, a_k] \times \prod_{k=u_1+1}^{u_1+u_2}\{(x,y)\in \mathbb{R}^2\, :\ x^2+y^2\leq a_k^2\},
	$$
	and
	$$
	N_s = \prod_{k=u_1+u_2+1}^{u_1+u_2+s_1}[-a_k, a_k] \times \prod_{k=u_1+u_2+s_1+1}^{u_1+u_2+s_1+s_2}\{(x,y)\in \mathbb{R}^2\, :\ x^2+y^2\leq a_k^2\},
	$$
\end{df}

\medskip

We also define
$$
N_{u,\varepsilon} = \prod_{k=1}^{u_1}[-a_k-\varepsilon, a_k+\varepsilon] \times \prod_{k=u_1+1}^{u_1+u_2}\{(x,y)\in \mathbb{R}^2\, :\ x^2+y^2\leq (a_k+\varepsilon)^2\}.
$$
If an element $x$ belongs to an h-set $A$ we can represent it uniquely as
$$
x = L((x_u,x_s)) + y,
$$
where $y \in \overline B_{X_2}(0,r)$, $x_u \in N_u \subset \mathbb{R}^u$ and $x_s \in N_s \subset \mathbb{R}^s$. We will use the notation $P_u x = x_u$, $P_{s}x = x_s$ and $P_{X_2}x = y$.
For an h-set $A$ we define its exit set as $$
A_{exit} = L \left( \partial N_u \times N_s\right) \oplus \overline{B}_{X_2}(0,r).
$$
and its $\varepsilon$ exit extension as
$$
A^\varepsilon = L \left( N_{u,\varepsilon} \times N_s\right) \oplus \overline{B}_{X_2}(0,r).
$$
Note that $P_u, P_{s}$ and $P_{X_2}$ make sense for elements of $A^\varepsilon$. An equivalent norm on $X$ will be denoted by $\|x\|_X = |P_u x| + |P_{s}x|+\|P_{X_2}x\|_{X_2}$, where by $|\cdot |$ we denote an euclidean norm on $\R^s$ or $\R^u$.
\begin{df}\label{isolation}
	Let $X$ be a Banach space and let $\{ S(t) \}_{t\geq 0}$ be a $C^0$ semiflow of mappings $S(t):X\to X$. An h-set $A$ is an \textit{isolating block} with respect to this semiflow if there exists $\varepsilon > 0$ and  $t(\varepsilon)>0$ such that for every $s\in (0,t(\varepsilon)]$
	\begin{itemize}
		\item[(A1)] $S(s)A\subset A^\varepsilon$,
		\item[(A2)] $\left[S(s)\left(A_{exit}\right)\right] \cap A = \emptyset$.
	\end{itemize}
	
\end{df}
Condition (A1) implies that if via the evolution $S(t)$ we leave an isolating block, we have to stay in $A^\varepsilon$ within the short time interval, while condition (A2) implies that if we are on the exit set of $A$, then, although we stay in $A^\varepsilon$,  we cannot reenter $A$ is a short time.
\begin{df}\label{df:17}
	The h-set  $A \subset X$ is called an \textit{h-set with cones} if there exist three continuous quadratic  forms $\alpha:\mathbb{R}^u \to \R$, $\beta:\mathbb{R}^{s} \to \R$ and $\gamma: X_2\to \R$ with
	\begin{align*}
		& m_\alpha |x|^2 \leq \alpha(x) \leq M_{\alpha}|x|^2\ \ \textrm{for every}\ \ x\in \mathbb{R}^u,\\
		& m_\beta |x|^2 \leq \beta(x) \leq M_{\beta}|x|^2\ \ \textrm{for every}\ \ x\in \mathbb{R}^{s},\\
		& m_\gamma \|y\|_{X_2}^2 \leq \gamma(y) \leq M_{\gamma}\|y\|_{X_2}^2\ \ \textrm{for every}\ \ y\in X_2,
	\end{align*}
	such that for every $x_1, x_2\in A$ satisfying $x_1\neq x_2$ the function
	$$
	t\mapsto \alpha(P_u (S(t)x_1-S(t)x_2)) - \beta(P_{s} (S(t)x_1-S(t)x_2)) - \gamma(P_{X_2} (S(t)x_1-S(t)x_2))
	$$
	is strictly increasing as long as both $S(t)x_1$ and $S(t)x_2$ stay in $A$.
	For short we will write, for $x\in X$
	$$
	Q(x) = \alpha(P_u(x)) - \beta(P_{s}(x)) - \gamma(P_{X_2}(x)).
	$$
\end{df}
Consider two points $x_1, x_2\in A$. If
$
Q(x_1-x_2) > 0,
$
then we will say that $x_1$ is in the positive cone of $x_2$ (and, equivalently, $x_2$ is in the positive cone of $x_1$), and if
$Q(x_1-x_2) < 0$
then we will say that $x_1$ is in the negative cone of $x_2$ (and, equivalently, $x_2$ is in the negative cone of $x_1$).
\begin{df}
	Let $\{ S(t) \}_{t\geq 0}$ be a $C^0$ semiflow on $X$. A point $x_0 \in X$ is an equilibrium if $S(t)x_0 = x_0$ for every $t\geq 0$. Let $A$ be and h-set such that $x_0\in A$ is an equilibrium.
	We define its local stable and unstable sets
	$$
	W^s_{loc,A}(x_0) = \{ x\in A\,:\ S(t)x\in A\ \ \textrm{for every}\ \ t\geq 0\ \ \textrm{and}\ \ \lim_{t\to \infty} S(t) x = x_0\},
	$$
	\begin{align*}
		& W^u_{loc,A}(x_0) = \{ x\in A\,:\ \textrm{there exists the function}\ \ u:(-\infty,0]\to A\ \ \textrm{such that} \\
		& \ \ \ \  u(0) = x, \lim_{s\to -\infty}u(s) = x_0\ \ \textrm{and for every}\ \  s\in (-\infty,0]\ \ \textrm{and}\ \ t\in [0,-s] \ \ \textrm{we have}\ \ S(t)u(s) = u(s+t)\},
	\end{align*}
\end{df}
We provide the theorem of the existence of a unique fixed points and local stable and unstable manifolds inside the isolating h-set with cones. Its proof is a version of Hadamard's proof of the existence of local stable and unstable manifolds and is based on a concept of the graph transform method.
\begin{thm}\label{graph_t}
	Let $A$ be an isolating h-set with cones for an asymptotically compact $C^0$ semiflow $\{ S(t) \}_{t\geq 0}$. Then there exist:
	\begin{itemize}
		\item a unique equilibrium $x_0$ in $A$,
		\item a Lipschitz continuous mapping
		$$
		F_s:L \left( \prod_{k=1}^{u} \{ 0 \}\ \times N_s\right) \oplus \overline{B}_{X_2}(0,r) \to A,
		$$
		with $P_s F_s(L(x_u,x_s)+y) = x_s$ and $P_{X_2} F_s(L(x_u,x_s)+y) = y$ such that $\textrm{im}\, F_s =  W^s_{loc,A}(x_0)$,
		\item a Lipschitz continuous mapping
		$$
		F_u:N_u \to A,
		$$
		with $P_u F_u(L(x_u,x_s)+y) = x_u$ such that $\textrm{im}\, F_u =  W^u_{loc,A}(x_0)$.
	\end{itemize}
\end{thm}
\begin{proof}
	
	\textbf{Step 1. Graph transform.}
	Consider a function  $h:N_u \to A$ with $P_u(h(x)) = x$ for every $x\in N_u$ such that for every $x_1, x_2\in N_u$ with $x_1\neq x_2$ the point $h(x_1)$ is in the positive cone of $h(x_2)$. We will call such function the horizontal disk. For every $x\in N_u$ consider $s\in (0,t(\varepsilon)]$ and observe that
	$$P_u S(s) (h(x)) \in N_{u,\varepsilon}.$$
	
	Choose $s\in (0,t(\varepsilon))$. We should show that for every $ x\in N_u$ there exists a unique $z\in N_u$  such that $S(r)h(z) \in A$ for every $r\in (0,s]$ and $ P_u S(s) (h(z)) = x$.
	We start from the proof of uniqueness. For the sake of contradiction assume that $P_u S(s) (h(z_1)) = P_u S(s) (h(z_2))$. By (A1) and (A2) we can use the cone condition whence
	\begin{align*}
		& \ \ 0 \geq - \beta(P_{s} (S(s)h(z_1)-S(s)h(z_2))) - \gamma(P_{X_2}(S(s)h(z_1)-S(s)h(z_2))) \\
		& \ \ \ \  = Q(S(s)h(z_1)-S(s)h(z_2)) > Q(h(z_1)-h(z_2)),
	\end{align*}
	which is a contradiction with the fact that $h$ is a horizontal disk.
	
	To prove the existence consider the map
	$$
	\Phi_s:N_u\to N_{u,\varepsilon}
	$$
	defined by
	$$
	N_u \ni{x} \mapsto   P_u S(s)(h(x)) \in N_{u,\varepsilon}
	$$
	and the map $\Psi_s:\mathbb{R}^u \to \mathbb{R}^u$
	$$
	\mathbb{R}^u\ni{x} \to e^{s} x \in \mathbb{R}^u.
	$$
	Define the homotopy
	$$
	f_r(x) = \Psi_{(1-r)s}(\Phi_{rs}(x)) \quad \textrm{for}\quad r\in [0,1].
	$$
	From (A1), (A2) and the fact the $\Psi_s$ is expanding and $\Psi_0$ is the identity we obtain that $N_u \cap \Psi_{(1-r)s}(\Phi_{rs}(\partial N_u))) = \emptyset$.
	This implies that
	$$ \mathrm{deg}(\Phi_1,\textrm{int}\ N_u,x) = \mathrm{deg}(\Psi_1,\textrm{int}\ N_u,x) \neq 0,$$
	for every $x\in N_u$. In consequence we get the needed existence.
	Define the set
	$$
	N_u \supset N_u(t) = \{ x \in N_u :  S(s)(h(x)) \in A \ \mathrm{for}\ s\in [0,t]\}
	$$ and a mapping
	$$
	N_u \ni P_u S(t)h(y)  \mapsto S(t)h(y) \in A\ \ \textrm{for some}\ \ y\in N_u(t).
	$$
	We have to prove that this mapping is a horizontal disk. This fact holds from the observation that the cone condition and the fact that $h$ is a horizontal disk imply
	\begin{align*}
		& Q(S(t)h(x_1)-S(t)h(x_2)) > Q(h(x_1)-h(x_2))> 0.
	\end{align*}
	Our aim is to prove that
	$$
	\bigcap_{t\geq 0} N_u(t) \neq \emptyset.
	$$
	This is a deceasing family of sets which are nonempty, bounded and closed and hence compact. Their intersection is nonempty and there exists $x\in N_u$ such that $S(t)(h(x)) \in A$ for every $t\geq 0$.

	\textbf{Step 2. Existence of unique equilibrium.} In this step we will prove that  there exists a unique $z_0 \in A$ such that if $S(t)z \in A$ for every $t\in [0,\infty)$ then $\lim_{t\to \infty}S(t)z = z_0$.
	Take $z\in A$ such that $S(t)z\in A$ for every $t\geq 0$ and let $\overline{z} \in \omega(z)$. We will show that a cone condition allows us to construct a Lyapunov function, and we will use the invariance principle. Let $S(t_n)z\to \overline{z}$. Note that $S(t)S(t_n)z \to S(t)\overline{z}$.
	The function
	$[0,\infty) \ni s\to Q(S(s)z - S(s)S(t)z)$ is nondecreasing and bounded from above. Hence $\lim_{s\to \infty} Q(S(s)z - S(t)S(s)z) = Q_0$. There holds $Q(\overline{z}-S(t)\overline{z}) = Q_0$. Assume that $\overline{z}\neq S(t)\overline{z}$. Then $Q(S(r)\overline{z}-S(r)S(t)\overline{z}) > Q_0$ for $r>0$.
	But
	$$ Q(S(r)S(t_n)z - S(r)S(t_n)S(t)z) = Q(S(t_n+r)z - S(t)S(t_n+r)z) \to Q_0$$
	and, simultaneously
	$$
	Q(S(r)S(t_n)z - S(r)S(t_n)S(t)z) \to Q(S(r)\overline{z}-S(r)S(t)\overline{z}) > Q_0,
	$$
	a contradiction. Hence $\overline{z}= S(t)\overline{z}$. Hence every $\overline{z} \in \omega(z)$ is an equilibrium. An immediate observation that uses the cone condition implies
	that the  equilibrium in $A$ must be unique. Hence $\omega(z) = \{z_0\}$ and $S(t)z \to z_0$ as $t\to \infty$.

	\textbf{Step 3. Local stable manifold.} We prove that for any horizontal disk $h:N_u\to A$ the point $h(x)$ such that its trajectory stays in $A$ is unique. We will denote such point $x_h\in N_u$. Indeed assume that
	there are two such points $h({x}_1)$ and $h({x}_2)$. Then both $S(t)(h({x}_1)) \to z_0$ and $S(t)(h({x}_2)) \to z_0$ as $t\to \infty$. Hence
	$$
	Q(h({x}_1) - h({x}_2)) < Q(S(t)(h({x}_1)) - S(t)(h({x}_2))) \to 0.
	$$
	On the other hand $Q(h({x}_1)-h({x}_2))>0$, a contradiction.
	For $z\in L \left( \prod_{k=1}^{u} \{ 0 \}\ \times N_s\right) \oplus \overline{B}_{X_2}(0,r)$ given by $z = L((0,P_s z)) + P_{x_2}z$ define the horizontal disk $h_z(x) = L((x,P_s z)) + P_{X_2}z$. There exists a unique point in this disk $x_{h_z}$ such that its trajectory stays in $A$ for all $t$. We denote $F_s(z) = x_{h_z}$. Graph of $F_s$ is a local stable manifold of the unique eqilibrium $z_0$. We prove that $F_s$ is Lipschitz. If $z_1\neq z_2$ then $F_s(z_1)$ and $F_s(z_2)$ cannot stay mutually in their positive cones, otherwise their trajectories could not converge to the same point (hence the map $F_s$ is a vertical disk). This means that
	$$
	\alpha(P_u(F_s(z_1) - F_s(z_2))) \leq \beta(P_s (z_1-z_2)) + \gamma(P_{X_2}(z_1-z_2)),
	$$
	which is enough to assert that $F_s$ is Lipschitz.
	
	\textbf{Step 4. Local unstable manifold}. Consider a horizontal disk $h:N_u \to A$ and the map $N_u(t)\ni x \to S(t)(h(x))\in A$. As it was established in Step 1, for every $t \geq 0$ there exists a horizontal disk with the image equal to the image of this map. Fix $t > 0, x \in  N_u$ and consider the sequence
	$$a_k := \left[S(kt)(h(N_u(kt)))\right] \cap \left[L \left( \{x\} \times N_s\right) \oplus \overline{B}_{X_2}(0,r)\right].$$ For every $k \in \mathbb{N}$ the intersection has exactly one point, so this sequence is well defined. By applying the diagonal argument to this sequence, we can find
	$$w_x \in L \left( \{x\} \times N_u\right) \oplus \overline{B}_{X_2}(0,r)$$ with infinite backward orbit in $A$. Indeed: since each $a_k$ has the backward orbit in $A$ with the length at least $t$, we can consider the sequence $\{b_k\}\subset A$ such that $S(t)b_k = a_k$, for all $k \in \mathbb{N}$. By the asymptotic compactness we can  pick a convergent subsequence of $b_k$ and, abusing the notation, we consider the  corresponding subsequence of $a_k$ without renumbering it. By the continuity of $S(t)$ we have $S(t) \lim b_k = \lim a_k$, so the subsequence $a_k$ has limit with backward orbit in $A$ of length at least $t$. We set $w_1 = a_1$. We take the subsequence of ${a_k}$ consisting of points which have the backward orbits in $A$ of time length at least $2t$ and we do not renumber it. We repeat the procedure to obtain the subsequence with the limit having the backward orbit in $A$ with time length $2t$ and take as $w_2$ the first element of this new subsequence. Then we continue the argument for time intervals of length $lt$ for every $l \in \mathbb{N}$ and each time we set $w_l = a_1$, the first element of the new subsequence. By construction, the limit of this diagonal sequence $w_x$ has infinite backward orbit $\{o_k\}_{k \in \mathbb{Z}_{\leq 0}}, o_{k+1} = S(t)o_k$ in $A$ for $k \in \mathbb{Z_-}, o_0 = w_x$. Since $V(z) := Q(z - z_0)$ is a Lyapunov function, it holds that $\lim_{k \to -\infty} o_k = z_0$. Thus $w_x \in W^u_{loc,A}(z_0)$. Assume that for some $z  \in L \left( \{x\} \times N_s\right) \oplus \overline{B}_{X_2}(0,r)$ such that $z\neq w_x$ for some $x$ there exists an infinite backward orbit $o'_k$ in $A$. Then
	$0 > Q(z - w_x) > \lim_{k\to -\infty} Q(o_k - o'_k) = Q(z_0 - z_0) = 0$, a contradiction.
	We define $F_u: N_u \ni x \mapsto w_x \in A$. This is the local unstable manifold, and by the argument analogous to the one in the step 3, it is a Lipschitz function.
\end{proof}

\section{Appendix 3: $C^0$ dependence of local unstable and stable manifolds on parameter.}
\subsection{Cone condition with parameter} Consider the family $\{S_\delta\}_{\delta \in [0,\Delta]}$ of semiflows on the space $X$ and a set $A\subset X$ which is an isolating h-set for every $\delta \in [0,\Delta]$.
\begin{df}\label{df:hconespar}
	Let $\{S_\delta(t)\}_{t\geq 0}$ given for  $\delta \in [0,\Delta]$ be $C^0$ semiflows and let $A \subset X$ be an h-set with cones for every $\delta \in [0,\Delta]$. We say that this set is a parameterized h-set with cones if there exist three continuous quadratic  forms $\alpha:\mathbb{R}^u \to \R$, $\beta:\mathbb{R}^{s} \to \R$ and $\gamma: X_2\to \R$
	\begin{align*}
		& m_\alpha |x|^2 \leq \alpha(x) \leq M_{\alpha}|x|^2\ \ \textrm{for every}\ \ x\in \mathbb{R}^u,\\
		& m_\beta |x|^2 \leq \beta(x) \leq M_{\beta}|x|^2\ \ \textrm{for every}\ \ x\in \mathbb{R}^{s},\\
		& m_\gamma \|x\|_{X_2}^2 \leq \gamma(x) \leq M_{\gamma}\|x\|_{X_2}^2\ \ \textrm{for every}\ \ x\in X_2,
	\end{align*}
	and a positive constant $L\in \mathbb{R}$ such that:
	\begin{itemize}
		\item[(i)] for every $x_1, x_2\in A$ and every $\delta_1,\delta_2\in [0,\Delta]$ if the function
		\begin{align*}
			&	[0,\infty) \ni t\mapsto  L|\delta_1-\delta_2|^2 + \alpha(P_u (S_{\delta_1}(t)x_1-S_{\delta_2}(t)x_2))\\
			& \qquad  - \beta(P_{s} (S_{\delta_1}(t)x_1-S_{\delta_2}(t)x_2)) - \gamma(P_{X_2} (S_{\delta_1}(t)x_1-S_{\delta_2}(t)x_2)) = \widehat{Q}(t)
		\end{align*}
		satisfies $\widehat{Q}(0) \geq 0$ then $\widehat{Q}(t)\geq 0$ as long as both $S_{\delta_1}(t)x_1$ and $S_{\delta_2}(t)x_2$ stay in $A$,
		\item[(ii)] for every $x_1, x_2\in A$ and every $\delta_1,\delta_2\in [0,\Delta]$ if the function
		\begin{align*}
			&[0,\infty) \ni t\mapsto  \alpha(P_u (S_{\delta_1}(t)x_1-S_{\delta_2}(t)x_2)) \\
			& \qquad - \beta(P_{s} (S_{\delta_1}(t)x_1-S_{\delta_2}(t)x_2)) - \gamma(P_{X_2} (S_{\delta_1}(t)x_1-S_{\delta_2}(t)x_2)) - L|\delta_1-\delta_2|^2  = \overline{Q}(t)
		\end{align*}
		satisfies $\overline{Q}(0) \geq 0$ then $\overline{Q}(t) \geq 0$ as long as both $S_{\delta_1}(t)x_1$ and $S_{\delta_2}(t)x_2$ stay in $A$,
		\item[(iii)] for every given $\delta\in [0,\Delta]$ the function
		$$
		[0,\infty) \ni t\mapsto \alpha(P_u (S_{\delta}(t)x_1-S_{\delta}(t)x_2)) - \beta(P_{s} (S_{\delta}(t)x_1-S_{\delta}(t)x_2)) - \gamma(P_{X_2} (S_{\delta}(t)x_1-S_{\delta}(t)x_2))
		$$
		is strictly increasing for every $x_1\neq x_2$ as long as both trajectories $S_\delta(t)x_1$ and $S_\delta(t)x_2$ stay in $A$.
	\end{itemize}
	If there exists a parameterized h-set with cones then the same $\alpha, \beta, \gamma$ can be used in the definition of an h-set with cones for every $\delta\in [0,\Delta]$, so every $S_\delta$ must have a unique equilibrium $x_0^\delta\in A$ and a local stable and unstable manifolds $W^s_{loc,A}(x_0^\delta), W^u_{loc,A}(x_0^\delta)$ given by the images of the Lipshitz functions
	$$
	F_u^\delta:N_u \to A, \qquad F_s^\delta :L \left( \prod_{k=1}^{u} \{ 0 \}\ \times N_s \right) \oplus \overline{B}_{X_2}(0,r) \to A,
	$$
\end{df}

\subsection{Lipshitz continuous dependence of local unstable manifolds on parameter} In the proof of the Lipshitz continuous dependence of local unstable manifolds on  parameter we will say that the pairs $(\delta_1,x_1)$ and $(\delta_2,x_2)$ belong mutually to their positive cones if
$$L|\delta_1-\delta_2|^2+\alpha(P_u (x_1-x_2))>  \beta(P_{s} (x_1-x_2)) + \gamma(P_{s_2} (x_1-x_2)) ,
$$ so we link the variable $\delta$ with the unstable variable $x_u$. We prove the following result.
\begin{thm}\label{thm:unstable}
	Let $A$ be an isolating parameterized h-set with cones with a constant $L>0$ for asymptotically compact $C^0$ semiflows $\{ S_\delta(t) \}_{t\geq 0}$ for $\delta \in [0,\Delta]$ and let $F_u^\delta$ be the Lipschitz functions such that $im F_u^\delta = W^u_{loc,A}(x_0^\delta)$. Then there exists a constant
	$C > 0$ such that for every $\delta_1, \delta_2 \in [0,\delta]$ and every $x_1,x_2 \in N_u$ we have
	$$
	\|F_u^{\delta_1}(x_1) - F_u^{\delta_2}(x_2)\|_X \leq C ( |\delta_1-\delta_2| + |x_1-x_2|).
	$$
\end{thm}
\begin{proof}
	The proof follows the lines of Steps 1 and 4 in the proof of Theorem \ref{graph_t}, where we additionally treat the extra variable $\delta$ (which is  constant in time) as one of unstable variables. We provide the details of the proof for the completeness of the exposition.
	
	\textbf{Step 1. Graph transform in extended variables.}
	Define
	$$
	N_{cu}=  [0,\Delta] \times N_u,
	$$
	and consider the function $h:N_{cu} \to A$ with $P_u(h(\delta,x)) = x$ such that for every $(\delta_1,x_1), (\delta_2,x_2)\in N_{cu}$ with $(\delta_1,x_1)\neq (\delta_2,x_2)$ the point $(\delta_1,h(\delta_1,x_1))$ is in the positive cone of $(\delta_2,h(\delta_2,x_2))$. Proceeding exactly as in the proof of Theorem \ref{graph_t}, for every $\delta\in [0,\Delta]$ and $t>0$ there exists the nonempty and compact set $N_{cu}(t,\delta) \subset N_u$ such that
	$$
	N_{cu} \supset \bigcup_{\delta\in [0,\Delta]}\{\delta\}\times N_{cu}(t,\delta) = \{ (\delta,x) \in N_{cu} :  S_\delta(s)(h(\delta,x)) \in A \ \mathrm{for}\ s\in [0,t]\}
	$$ and the mapping
	$$
	N_{cu} \ni (\delta,P_u S_\delta(t)h(z))  \mapsto S_\delta(t)h(z) \in A\ \ \textrm{for some}\ \ z\in N_{cu}(t,\delta)
	$$
	is a horizontal disk, i.e. any two points in its graph belong mutually to their positive cones.

	\textbf{Step 2. Local unstable manifold in extended variables}. We proceed as in step 4 of the proof of Theorem \ref{graph_t}. From the previous step, by evolving the horizontal disk $h:N_{cu} \to A$ by the family of semigroups $\{S_\delta\}_{\delta\in [0,\Delta]}$ we obtain horizontal disks for every $t>0$. Fix $t > 0, (\delta,x) \in  N_{cu}$ and consider the sequence obtained by intersecting the horizontal disk with the vertical segment
	$$a_k(\delta,x) := \left[S_\delta(kt)(h(N_{cu}(kt,\delta)))\right] \cap \left[L \left( \{x\} \times N_s\right) \oplus \overline{B}_{X_2}(0,r)\right].$$ As we have shown in step 4 in the proof  of Theorem \ref{graph_t} this sequence has a convergent subsequence and the limit $w_{x,\delta}$ has an infinite backward trajectory via $S_\delta$ convergent backward in time to the unique equilibrium $z_0^\delta$ of $S_\delta$ in $A$. Moreover, for every $\delta\in [0,\Delta]$ the limit $w_{x,\delta}$ is the unique point among the points $z$ with $P_uz=x$ with the infinite backward trajectory in $A$. This uniqueness implies that the whole sequence $a_k(\delta,x)$ converges to $w_{x,\delta}$. We can define the mapping $F_{cu}:N_{cu}\ni (\delta,x)\to w_{x,\delta}\in A$. For every $\delta \in [0,\Delta]$ we have $\textrm{im}\, F_{cu}(\delta,\cdot) = W^u_{loc,A}(z^\delta_0)$.
	To show that $F_{cu}$ is Lipschitz observe that for every $(\delta_1,x_1),(\delta_2,x_2)\in N_{cu}$ the points $a_k(\delta_1, x_1)$ and $a_k(\delta_2,x_2)$  belong to the same horizontal disk so they also belong to each other's positive cones. Hence, for every $k$ we have
	$$
	0 < L|\delta_1-\delta_2|^2 + \alpha(x_1-x_2) - \beta(P_{s} (a_k(\delta_1, x_1)-a_k(\delta_2,x_2)))) - \gamma(P_{X_2} (a_k(\delta_1, x_1)-a_k(\delta_2,x_2)))),
	$$
	and passing to the limit with $k\to \infty$ we obtain
	$$
	\beta(P_{s} (w_{x_1,\delta_1}-w_{x_2,\delta_2})) + \gamma(P_{X_2} (w_{x_1,\delta_1}-w_{x_2,\delta_2})) \leq L|\delta_1-\delta_2|^2 + \alpha(x_1-x_2),
	$$
	which leads to the required Lipschitz condition.
	
\end{proof}

\subsection{Lipshitz continuous dependence  of local stable manifolds on parameter} In the proof of the Lipshitz continuous dependence of local stable manifolds on  parameter the key role will be played by the cone condition given in item (ii) of Definition \ref{df:hconespar}. We will now say that the pairs $(\delta_1,x_1)$ and $(\delta_2,x_2)$ belong mutually to their positive cones if
$$\alpha(P_u (x_1-x_2))>  \beta(P_{s} (x_1-x_2)) + \gamma(P_{X_2} (x_1-x_2)) + L |\delta_1-\delta_2|^2.
$$ We prove the following result.
\begin{thm}
	Let $A$ be an isolating parameterized h-set with cones with a constant $L>0$ for asymptotically compact $C^0$ semiflows $\{ S_\delta(t) \}_{t\geq 0}$ for $\delta \in [0,\Delta]$ and let $F_s^\delta$ be the Lipschitz functions such that $im F_s^\delta = W^s_{loc,A}(x_0^\delta)$. There exists a constant
	$C > 0$ such that for every $\delta_1, \delta_2 \in [0,\delta]$ and every $z_1,z_2 \in L \left( \prod_{k=1}^{u} \{ 0 \}\ \times N_s\right) \oplus \overline{B}_{X_2}(0,r)$ we have
	$$
	\alpha(P_u(F^{\delta_1}_s(z_2) - F^{\delta_2}_s(z_1))) \leq  \beta(P_s(z_1-z_2)) + \gamma(P_{X_2}(z_1-z_2)) + L |\delta_1-\delta_2|^2.
	$$
\end{thm}
\begin{proof}
	Again the proof follows the lines of Steps 1 and 3 in the proof of Theorem \ref{graph_t}.
	
	\textbf{Step 1. Graph transform in extended variables.}
	As in Step 1 of the proof of Theorem \ref{thm:unstable} we define
	$$
	N_{cu}=  [0,\Delta] \times N_u,
	$$
	and consider the function $h:N_{cu} \to A$ with $P_u(h(\delta,x)) = x$ such that for every $(\delta_1,x_1), (\delta_2,x_2)\in N_{cu}$ with $(\delta_1,x_1)\neq (\delta_2,x_2)$ the point $(\delta_1,h(\delta_1,x_1))$ is in the positive cone of $(\delta_2,h(\delta_2,x_2))$, now with respect to $\overline{Q}$. Again, evolving the graph of this function we obtain a family of horizontal disks in extended variables parameterized by time.

	\textbf{Step 2. Local stable manifold in extended variables}. Exactly as in step 3 of the proof of Theorem \ref{graph_t}, for $z\in L \left( \prod_{k=1}^{u} \{ 0 \}\ \times N_s\right) \oplus \overline{B}_{X_2}(0,r)$ given by $z = L((0,P_s z)) + P_{X_2} z$ define the horizontal disk $h_z(\delta,x) = L((x,P_s z)) + P_{X_2}z$. This disk, after time $t$ transforms to the horizontal disk $h_{z,t}(\delta,x)$. Let us define the mapping
	$$
	[0,\Delta]\times L \left( \prod_{k=1}^{u} \{ 0 \}\ \times N_s \right) \oplus \overline{B}_{X_2}(0,r) \ni (\delta,z) \mapsto f_t(\delta,z) = L(x(t,\delta,z),P_s z) + P_{X_2}z \in A,
	$$
	where $x(t,\delta,z) \in N_u$ is such a point that $P_u(h_{z,t}(\delta,x(t,\delta,z))) = 0.$ We prove that this mapping is a vertical disk, that is, that
	$$\alpha(x(t,\delta_1,z_1)-x(t,\delta_2,z_2)) \leq  \beta(P_s(z_1-z_2)) + \gamma(P_{X_2}(z_1-z_2)) + L |\delta_1-\delta_2|^2.
	$$
	Indeed, if the opposite inequality holds
	$$\alpha(x(t,\delta_1,z_1)-x(t,\delta_2,z_2)) >  \beta(P_s(z_1-z_2)) + \gamma(P_{X_2}(z_1-z_2)) + L |\delta_1-\delta_2|^2,
	$$
	then points $(\delta_1,f_t(\delta_1,z_1))$ and $(\delta_2,f_t(\delta_2,z_2))$ belong mutually to their positive cones, whence, after time $t$, we should have, that
	\begin{align*}
		& \alpha(0-0) \geq  \beta(P_{s}(h_{z_1,t}(\delta_1,x(t,\delta_1,z_1))-h_{z_2,t}(\delta_2,x(t,\delta_2,z_2)))) \\
		& \qquad + \gamma(P_{X_2}(h_{z_1,t}(\delta_1,x(t,\delta_1,z_1))-h_{z_2,t}(\delta_2,x(t,\delta_2,z_2)))) + L |\delta_1-\delta_2|^2,
	\end{align*}
	which would mean that $\delta_1=\delta_2=\delta$ and $h_{z_1,t}(\delta,x(t,\delta,z_1)) = h_{z_2,t}(\delta,x(t,\delta,z_2)).$ But this means that
	$$\alpha(x(t,\delta,z_1)-x(t,\delta,z_2)) >  \beta(P_s(z_1-z_2)) + \gamma(P_{X_2}(z_1-z_2)),
	$$
	i.e. $Q(f_t(\delta,z_1)-f_t(\delta,z_2))>0$, whence, after time $t$
	\begin{align*}
		& 0 = \alpha(0-0) >  \beta(P_{s}(h_{z_1,t}(\delta,x(t,\delta,z_1))-h_{z_2,t}(\delta,x(t,\delta,z_2)))) \\
		& \qquad + \gamma(P_{X_2}(h_{z_1,t}(\delta,x(t,\delta,z_1))-h_{z_2,t}(\delta,x(t,\delta,z_2)))) = 0,
	\end{align*}
	a contradiction. We prove that for every $(\delta,z)$ there holds
	$$
	\lim_{t\to \infty} f_t(\delta,z) = F_s^\delta(z).
	$$
	Indeed, for a given fixed $z$ and $\delta$ the stable part of $f_t(\delta,z)$ is constant in time and equal to $z$ and the unstable part given by $x(t,\delta,z)$ belongs to the sets $N_u(t)$ (depending also on $\delta$ and $z$) given in Step 1 of the proof of Theorem \ref{graph_t}, i.e. those points in the horizontal disk $h_z(\delta,\cdot)$ whose trajectory stays in $A$ for time at least $t$. The sets $N_u(t)$ are a decreasing family of nonempty and compact sets, whose intersection is a singleton given by $P_s(F_s^\delta(z))$. We can pass to the limit with $t$ to infinity in the vertical disk condition
	$$\alpha(P_u(f_t(\delta_1,z_1)-f_t(\delta_2,z_2))) \leq  \beta(P_s(z_1-z_2)) + \gamma(P_{X_1}(z_1-z_2)) + L |\delta_1-\delta_2|^2,
	$$
	which yields
	$$\alpha(P_u(F^{\delta_1}_s(z_2) - F^{\delta_2}_s(z_1))) \leq  \beta(P_s(z_1-z_2)) + \gamma(P_{X_2}(z_1-z_2)) + L |\delta_1-\delta_2|^2,
	$$
	the assertion of the theorem.
\end{proof}

\section{Appendix 4: $C^1$ smoothness of local stable and unstable manifolds}

\subsection{Fibre contraction theorem}
The following result is known as the fiber contraction theorem \cite[Theorem 1.2]{Hirsch_Pugh}
\begin{thm}\label{thm:fibre}
	Let $(X,\varrho_X), (Y, \varrho_Y)$ be complete metric spaces and let $f:X\to X$ and $g:X\times Y \to  Y$ be continuous maps such that
	\begin{align*}
		&\varrho_X(f(x_1),f(x_2)) \leq \lambda_1  \varrho_X(x_1,x_2)\ \ \textrm{for every}\ \ x_1,x_1\in X,\\
		&\varrho_Y(g(x,y_1),g(x,y_2)) \leq \lambda_2  \varrho_Y(y_1,y_2)\ \ \textrm{for every}\ \ x\in X,y_1,y_2\in Y,
	\end{align*}
	where $\lambda_1,\lambda_2\in (0,1)$. Then there exists a unique pair $(x_\infty,y_\infty)\in X\times Y$ such that $f(x_\infty)=x_\infty$, $g(x_\infty,y_\infty) = y_\infty$. Moreover $(x_\infty,y_\infty)$ is attracting.
\end{thm}
The mapping $X\times Y \in (x,y) \mapsto \Lambda(x,y) = (f(x),g(x,y)) \in X\times Y$ in the above theorem is called a fibre contraction.
\begin{thm}\label{thm:param}
	Suppose we have a family of fibre contractions $\Lambda^\varepsilon$ depending on the parameter $\varepsilon \in [0,\varepsilon_0]$ with constants $\lambda_1, \lambda_2$ such that $\Lambda^\varepsilon(x,y) = \Lambda(\varepsilon,x,y)$ is continuous. Then, for their fixed points, we have
	\begin{align*}
	& \lim_{\varepsilon\to 0}\varrho_X(x_\infty^\varepsilon,x_\infty^0) = 0,\\
	& \lim_{\varepsilon\to 0}\varrho_Y(y_\infty^\varepsilon,y_\infty^0) = 0.
	\end{align*}
	\end{thm}  
\begin{proof}
	We have
	\begin{align*}
	& \varrho_X(x_\infty^\varepsilon,x_\infty^0) = \varrho_X(f^\varepsilon(x_\infty^\varepsilon),f^0(x_\infty^0)) \leq \varrho_X(f^\varepsilon(x_\infty^\varepsilon),f^\varepsilon(x_\infty^0))+ \varrho_X(f^\varepsilon(x_\infty^0),f^0(x_\infty^0))\\
	& \ \ \leq \lambda_1 \varrho_X(x_\infty^\varepsilon,x_\infty^0) + \varrho_X(f^\varepsilon(x_\infty^0),f^0(x_\infty^0)).
	\end{align*}
This means that
$$
\varrho_X(x_\infty^\varepsilon,x_\infty^0) \leq \frac{1}{1-\lambda_1} \varrho_X(f^\varepsilon(x_\infty^0),f^0(x_\infty^0)),
$$
and the first desired convergence follows. Next,
	\begin{align*}
	& \varrho_Y(y_\infty^\varepsilon,y_\infty^0) = \varrho_Y(g^\varepsilon(x_\infty^\varepsilon,y_\infty^\varepsilon),g^0(x_\infty^0,y_\infty^0))\\
	& \ \ \ \leq \varrho_Y(g^\varepsilon(x_\infty^\varepsilon,y_\infty^\varepsilon),g^\varepsilon(x_\infty^\varepsilon,y_\infty^0))+ \varrho_Y(g^\varepsilon(x_\infty^\varepsilon,y_\infty^0),g^0(x_\infty^0,y_\infty^0))\\
	& \ \ \leq \lambda_2 \varrho_Y(y_\infty^\varepsilon,y_\infty^0) + \varrho_Y(g^\varepsilon(x_\infty^\varepsilon,y_\infty^0),g^0(x_\infty^0,y_\infty^0)).
\end{align*}
Hence
$$
\varrho_Y(y_\infty^\varepsilon,y_\infty^0) \leq \frac{1}{1-\lambda_2} \varrho_Y(g^\varepsilon(x_\infty^\varepsilon,y_\infty^0),g^0(x_\infty^0,y_\infty^0)),
$$
and the proof is complete by continuity.
	\end{proof}
\subsection{Isolation and cone conditions.}
If $Z$ is a normed space then for a linear map $A:Z\to Z$ we define $m(A)$ as the largest constant $L\geq 0$ such that $\|Ax\|\geq L\|x\|$ for every $x\in Z$. If $A$ is singular (has nonzero kernel), than $m(A) = 0$. If $A:\R^n\to \R^n$ is a nonsingular matrix, then $m(A) = \|A^{-1}\|^{-1}$. 

In this section we construct the local stable and unstable manifolds for an equilibrium of a map $f:Z\to Z$, where $Z$ is a Banach space which has a Cartesian product structure $Z=\mathcal{X}\times \mathcal{Y}$. We assume that $z_0$ is a hyperbolic fixed point for $f: Z \to Z$, which of class $C^1$ and that $Z\ni z=(x,y)\in \mathcal{X}\times \mathcal{Y}$, where $x$ is unstable direction and $y$ is the stable direction.  
\subsubsection{Isolation conditions.} Assume that the equilibrium $z_0$ belongs to the set $N=\overline{B}_u(0,r_u) \times \overline{B}_s(0,r_s)  \subset \mathcal{X} \times \mathcal{Y}$,
and $N$ is an isolating $h$-set, cf. Definition \ref{isolation}, that is 
\begin{itemize}
	\item[(1)] $f(N) \subset N^\varepsilon$,
	\item[(2)] $f(N_{exit})\cap N = \emptyset,$
\end{itemize}
where $N^\varepsilon = \overline{B}_u(0,r_u+\varepsilon) \times \overline{B}_s(0,r_s)$ and $N_{exit} = \overline{B}_u(0,r_u+\varepsilon) \times \overline{B}_s(0,r_s)$. We stress that throughout this section the space $\mathcal{Y}$ can be infinite dimensional, but we require in this section that $u=\dim \mathcal{X} < \infty$.

\begin{rem}
	We will apply the results of this section to the framework defined in Appendix 2, where in place of the ball $\overline{B}_s(0,r_s)$ we consider the "box-like" set $N_s$ and in place of $\overline{B}_u(0,r_u)$ we consider the "box-like set $N_u \times \overline{B}_{X_2}(0,r)$. These sets are given in Definition \ref{h-set}. While, for the sake of notation simplicity, the argument of this section is done for balls in stable and unstable spaces, it also works for these more general sets $N_s$ and $N_u \times \overline{B}_{X_2}(0,r)$.   
\end{rem}

\begin{rem}
	The coordinates with which we work in Appendix 2 are the original coordinates of the system. Hence, the points in the set $N$ are represented there as $x=L((x_u,x_s))+y$, where $L$ is a nonsingular matrix, $(x_u,x_s) \in \R^{u+s} = \R^d$ and $y\in X_2=L^2_A(\R_+)^d$.  Here, we work with transformed coordinates in the finite dimensional variable, such that the stable and the unstable variables are separated. Thus, $\mathcal{X}=\R^u$ and $\mathcal{Y}=\R^s\times L^2_A(\R_+)^d$. This allows us to focus on the graph transform and avoid technicalities associated with the coordinate systems. 
\end{rem}
	
The isolation conditions for the constructed sets are verified in Sections \ref{subsec:cont-iso-block} and \ref{sec:isoeps}.  

\subsubsection{Cone conditions.} Cone conditions in these section will be formulated in different way than in Appendix 2, where in Definition \ref{df:17} the cone condition is given via quadratic forms. Here they assume the form of the bounds on the derivative of the stable and unstable parts of the mapping $f$ with respect to the stable and unstable components of the argument variable. The verification that they hold for the considered system will be the content of  Appendix 5. We consider the following set which we call a cone 
\begin{equation}
	C_u=\{ (x,y) \ : \  \|y\| \leq L \|x\| \}
\end{equation}
for some $L>0$. We define some constants

\begin{align}
	  &\xi= m\left(\frac{\partial f_x}{\partial x} \right) - L\left\|\frac{\partial f_x}{\partial y} \right\|, \label{eq:xi} \\
	&\mu=\frac{1}{L}\left\|\frac{\partial f_y}{\partial x}\right\|  +  \left\|\frac{\partial f_y}{\partial y} \right\|  \label{eq:mu},\\
	&\beta = \frac{\mu}{\xi} L \left\|\frac{\partial f_x}{\partial y}\right\|  +  \left\|\frac{\partial f_y}{\partial y} \right\|,  \label{eq:beta} \\
	& \xi_1 = m\left( \frac{\partial f_x}{\partial x}\right) -\frac{1}{L} \left\|\frac{\partial f_y}{\partial x}\right\|, \label{eq:xi1} \\
	&\mu_1 =   \left\|\frac{\partial f_y}{\partial y}\right\| + L \left\|\frac{\partial f_x}{\partial y}\right\| . \label{eq:mu1}
\end{align}

Note that for $z_1=(x_1,y_1)$, $z_2=(x_2,y_2)$, such that $z_1 - z_2 \in C_u$ we have
\begin{eqnarray}
  \|f_x(z_1) - f_x(z_2)\| &\geq& \left(m\left(\frac{\partial f_x}{\partial x}[z_1,z_2] \right) - L\left\|\frac{\partial f_x}{\partial y}[z_1,z_2] \right\|\right) \|x_1 - x_2\|, \label{eq:dfx} \\
  \|f_y(z_1) - f_y(z_2)\| &\leq& \left(\left\|\frac{\partial f_y}{\partial x}[z_1,z_2]\right\|  +  L\left\|\frac{\partial f_y}{\partial y}[z_1,z_2] \right\| \right) \|x_1 - x_2\|
  \label{eq:dfy}
\end{eqnarray}

For $z_1 - z_2 \notin C_u$ we obtain
\begin{equation}
  \|f_y(z_1)-f_y(z_2)\| \leq \mu \|y_1 - y_2\|. \label{eq:contr-not-cu}
\end{equation}

If $\mu\leq \xi$ then the graph transform for unstable manifold is well defined and if $\beta < 1$, then the graph transform for unstable manifold is a contraction (Thm.~\ref{thm:T-contr}) and the same holds also for $C^1$ the graph transform if $\beta < \min\{1,\xi,\xi^2\}$ (Thm.~\ref{thm:DT-contr}).

If $\mu_1\leq \xi_1$ then the graph transform for the stable manifold is well defined and if $\xi_1 >1$, then the graph transform for stable manifold is a contraction (Thm.~\ref{thm:S-contrakcja}). The same holds for the $C^1$-graph transform (Thm.~\ref{thm:25}) if $\xi_1 > \max\{1,\mu,\mu^2\}$.

Hence, in order to show that the local stable and unstable manifolds are $C^1$, it is enough to show the following five inequalities
\begin{equation}\label{eq:bounds_constants}
	\xi > 1,\ \  \mu < 1,\ \   \beta < 1,\ \  \xi_1 > 1, \ \ \mu_1 < 1.
\end{equation}
In order to get them it suffices to show that
\begin{align*}
	& (1) \ \ m\left(\frac{\partial f_x}{\partial x}\right) > 1,\ \  \left\|\frac{\partial f_y}{\partial y}\right\| < 1,\\
	& (2) \left\|\frac{\partial f_x}{\partial y}\right\|\  \textrm{can be made arbitrarily small by decreasing, if necessary, the set}\ N\ \textrm{and parameter}\ \varepsilon. 
\end{align*}
Note that the value $ \left\|\frac{\partial f_y}{\partial x}\right\|$ does not have to be small. Indeed, it appears in the constants $\mu$ and $\xi_1$ and is always multiplied by $\frac{1}{L}$. So, if only $m\left(\frac{\partial f_x}{\partial x}\right) > 1$, we can always choose $L$ large enough so that $\xi>1$. Likewise, if only $\left\|\frac{\partial f_y}{\partial y}\right\| < 1$, we can always choose $L$ large enough n order to guarantee that $\mu<1$. So, once (1) is satisfied, and $ \left\|\frac{\partial f_y}{\partial x}\right\|$ is found, we choose large $L$ to guarantee that $\mu<1$ and $\xi_1>1$, and then, for this $L$ we decrease $ \left\|\frac{\partial f_x}{\partial y}\right\|$ to guarantee that $\xi>1$, $\beta<1$, and $\mu_1<1$.

\subsection{Fixed point procedure for unstable manifold and its derivative}\label{fix_unstable}
\subsubsection{Graph transform for the unstable manifold}

\begin{df}
\label{def:hd-cc}
For a continuous map $y:\overline{B}_u(0,r_u) \to \overline{B}_s(0,r_s)$
we will say that $(x,y(x))$ is \emph{horizontal disk satisfying cone condition} if
\begin{equation}
  \|y(x_1)-y(x_2)\| \leq  L\|x_1-x_2\|. \label{eq:hordisk}
\end{equation}
\end{df}

\begin{df}
Let $H \subset C^{0}(\overline{B}_u(0,r_u),\overline{B}_s(0,r_s))$ be given defined as follows: $h \in H$ if and only if $h$ is a horizontal disk satisfying cone condition.
\end{df}
Observe that $H$ is closed. Assume that $(x,y(x))$ is an unstable manifold of $z_0$. Then we have
\begin{equation}
  f_y(x,y(x))=y(f_x(x,y(x))).
\end{equation}
We are in position to define the graph transform $\mathcal{T}:H\to H$ by the formula 
\begin{equation}
	f_y(x,h(x))=\mathcal{T}(h)(f_x(x,h(x))).
\end{equation}
The proof that for every $\overline{x}\in 
\overline{B}_u(0,r_u)$ there exists $x\in \overline{B}_u(0,r_u)$ such that $\overline{x} = f_x(x,h(x))$ uses the isolation conditions, and proceeds with the use of the Brouwer degree analogously as in the proof of Theorem \ref{graph_t}. The proof that the value of the graph transform $\mathcal{T}(h)$ given by $f_y(x,h(x))$ is defined uniquely as well as that  $\mathcal{T}(h)$ maps horizontal disks to horizontal disks is given in Theorem \ref{thm:T-contr}. This implies that the graph transform  $\mathcal{T}(h)$ is well defined. 

Define, implicitly, the mapping $G(h)$ as $G(h)(x) = \overline{x}$ such that $x = f_x(\overline{x},h(\overline{x}))$. In other words, $G(h)(x)$ satisfies the following implicit equation
\begin{equation}
	f_x(G(h)(x),h(G(h)(x)))=x.   \label{eq:impl-G}
\end{equation}
    Observe that using the map $G$  we can write the graph transform as follows
\begin{equation}
	\mathcal{T}(h) (x) = f_y(G(h)(x),h(G(h)(x))).  \label{eq:TG}
\end{equation}

%If  $x \mapsto f_x(x,y_u(x))$ is invertible, then
%\begin{equation}
%  y_u(x)=f_y \circ (Id,y_u) \circ (f_x \circ (Id,y_u))^{-1}.
%\end{equation}
%This is a fixed point equation for function $y_u$, which is called a \emph{graph transform} and is defined as follows:
%\begin{itemize}
% \item given   $h \in H$ (later we will assume $C^1$),
%     we see that the map $f \circ (id, h)$ is a horizontal disk satisfying cone condition,
%     hence $\overline{B}_u(0,r_u) \ni x \mapsto f_x(x,h(u))$ is invertible on $\overline{B}_u(0,r_u)$, in particular from (\ref{eq:dfx}) we obtain
%     \begin{equation}
%       \|f_x(x_1,h(x_1)) - f_x(x_2,h(x_2))\| \geq \xi \|x_1 - x_2\|,  \label{eq:h-xi}
%     \end{equation}
%     hence for the inverse map we have
%     \begin{equation}
%       \|(f_x\circ (Id, h))^{-1}(x_1) -  (f_x\circ (Id, h))^{-1}(x_2)\| \leq \frac{1}{\xi} \|x_1 - x_2\|  \label{eq:fx-inv-e}
%     \end{equation}
%
%
%For $y=h(x) \in H$  let us define
%\begin{equation}
%  G(h)=(f_x\circ (Id, h))^{-1},  \label{eq:G-def}
%\end{equation}
%which maps $H$ to maps $C^0(\overline{B}_u(0,r_u),\overline{B}_u(0,r_u))$
%
%Now (\ref{eq:fx-inv-e}) can be rewritten as
%\begin{equation}
%   \|G(h)(x_1) -  G(h)(x_2)\| \leq \frac{1}{\xi} \|x_1 - x_2\|. \label{eq:Gh-lip}
%\end{equation}
%
%
%  \item we define
%    \begin{equation}
%    \mathcal{T}(h) = f_y \circ (Id, h) \circ  (f_x\circ (Id, h))^{-1}. \label{eq:graphtran}
%    \end{equation}
%    Observe that using map $G$  we can write the graph transform as follows
%\begin{equation}
%   \mathcal{T}(h) (x) = f_y(G(h)(x),h(G(h)(x))).  \label{eq:TG}
%\end{equation}
%
%\end{itemize}

\begin{rem}
	The following lemma, which implies the uniform convergence of the graph transform, is proved in \cite{Geom,Zgliczynski}.\begin{lem}
There exists $K, \mu$, such that for any natural $m$ and
for any horizontal disks $h_1,h_2$ we have
\begin{equation}
\|\mathcal{T}^m(h_1) - \mathcal{T}^m(h_2)\| \leq K \mu^m
\end{equation}
\end{lem}
\end{rem}

%\subsection{Estimates on G in $C^0$-norm}

\begin{lem}
\label{lem:G-estm}
Let $\xi>0$. Then the mapping $G$ is well defined, and, assuming that $h_1,h_2 \in H$, we have
  \begin{equation}
    \|G(h_1)(x_1) - G(h_2)(x_2)\| \leq      \frac{1}{\xi} \left\| \frac{\partial f_x}{\partial y}  \right\|  \|h_1 - h_2\| + \frac{1}{\xi}\|x_1-x_2\|.
  \end{equation}
\end{lem}
\begin{proof}
Let us fix $x_1, x_2 \in \overline{B}_u(0,r_u)$ and  let us denote $\overline{x}_i=G(h_i)(x_i)$.  By definition of $G$ we have $f_x(\overline{x}_i,h_i(\overline{x}_i))=x_i$, hence
\begin{align*}
  & \|x_1 - x_2\|= \|f_x(\overline{x}_1,h_1(\overline{x}_1)) - f_x(\overline{x}_2,h_2(\overline{x}_2))\|\\
   &\geq m\left(\frac{\partial f_x}{\partial x}\right) \|\overline{x}_1 - \overline{x}_2\| - \left\| \frac{\partial f_x}{\partial y}  \right\| \cdot \|h_1(\overline{x}_1) - h_2(\overline{x}_2)\|
\end{align*}
But
$$
	\|h_1(\overline{x}_1) - h_2(\overline{x}_2)\| \leq \|h_1(\overline{x}_1) - h_1(\overline{x}_2)\| +   \|h_1(\overline{x}_2) - h_2(\overline{x}_2)\| 
	\leq L \|\overline{x}_1-\overline{x}_2\| + \|h_1 - h_2\|,
$$
and hence
$$
\|x_1-x_2\| \geq \left(m\left(\frac{\partial f_x}{\partial x}\right) - L \left\| \frac{\partial f_x}{\partial y}\right\|\right)\|\overline{x}_1 - \overline{x}_2\| - \left\| \frac{\partial f_x}{\partial y}  \right\| \cdot \|h_1 - h_2\|,
$$
which immediately implies the assertion. 
\end{proof}

The following estimate is crucial for proving that the graph transform is a contraction.

\begin{thm}
\label{thm:T-contr}

Let $\xi>0$. For any $h_1,h_2 \in H$ and $x_1, x_2 \in \overline{B}_u(0,r_u)$ the following estimate holds
\begin{equation}
 \| \mathcal{T}(h_1)(x_1) - \mathcal{T}(h_2)(x_2) \| \leq \beta \|h_1 - h_2\| + L\frac{\mu}{\xi}\|x_1-x_2\|.
\end{equation}
\end{thm}
\begin{proof}
Assume that $h_1,h_2 \in H$ and $x_1, x_2 \in \overline{B}_u(0,r_u)$. We have
\begin{align*}
 & \|\mathcal{T}(h_1)(x_1) - \mathcal{T}(h_2)(x_2)\|=  \|f_y(G(h_1)(x_1),h_1(G(h_1)(x_1))) - f_y(G(h_2)(x_2),h_2(G(h_2)(x_2))) \| \\
 &  \leq \left\|\frac{\partial f_y}{\partial x} \right\| \cdot \|G(h_1)(x_1)-G(h_2)(x_2)\| + \left\|\frac{\partial f_y}{\partial y} \right\|\cdot \|h_1 (G(h_1)(x_1)) - h_2 (G(h_2)(x_2))\|. 
\end{align*}
Since
\begin{align}
&  \|h_1 (G(h_1)(x_1)) - h_2  (G(h_2)(x_2))\| \nonumber\\
& \qquad \qquad \leq \|h_1  (G(h_1)(x_1)) - h_1  (G(h_2)(x_2))\|  
+ \|h_1 (G(h_2)(x_2)) - h_2  (G(h_2)(x_2))\|  \nonumber \\
& \qquad \qquad \leq L \|G(h_1)(x_1) - G(h_2)(x_2)\| + \|h_1 - h_2\|,\label{eq:gdiff}
\end{align}
 we obtain
 $$
  \|\mathcal{T}(h_1)(x_1) - \mathcal{T}(h_2)(x_2)\| \leq \left(\left\|\frac{\partial f_y}{\partial x} \right\| + L  \left\|\frac{\partial f_y}{\partial y} \right\|\right) \|G(h_1)(x_1)-G(h_2)(x_1)\| + \left\|\frac{\partial f_y}{\partial y} \right\|\|h_1 - h_2\|.
 $$
We are in position to use  Lemma~\ref{lem:G-estm}, whence 
\begin{align*}
   &\|\mathcal{T}(h_1)(x_1) - \mathcal{T}(h_2)(x_2)\|\leq 
  L\mu \left( \frac{1}{\xi} \left\| \frac{\partial f_x}{\partial y}  \right\|  \|h_1 - h_2\| + \frac{1}{\xi}\|x_1-x_2\|\right)
   + \left\|\frac{\partial f_y}{\partial y} \right\| \cdot \|h_1 - h_2\|  \\
   &\leq \left(  \frac{\mu}{\xi} L \left\|\frac{\partial f_x}{\partial y}\right\|  +  \left\|\frac{\partial f_y}{\partial y} \right\|  \right) \cdot \|h_1 - h_2\| + L\frac{\mu}{\xi} \|x_1-x_2\|,
\end{align*}
and the proof is complete.
\end{proof}
We easily deduce the following two results
\begin{thm}\label{incl1}
	If $\mu\leq \xi$ then $\mathcal{T}(H)\subset H$.
	\end{thm}
\begin{thm}\label{contr2}
	If $\beta < 1$ then $\mathcal{T}$ is a contraction.
	\end{thm}

\subsubsection{Graph transform for the derivative of unstable manifold}\label{graph_der}

Let us fix  $ h \in H \cap C^1$. We first differentiate $G$ with respect to $x$, we will denote the differentiation symbol by $D$. By applying such differentiation  with respect to $x$ to \eqref{eq:impl-G} we obtain
\begin{eqnarray}
	\left(\frac{\partial f_x}{\partial x}(G(h)(x),h(G(h)(x))) + \frac{\partial f_x}{\partial y}(G(h)(x),h(G(h)(x))) Dh(G(h)(x)) \right)D(G(h))(x)=I. \label{eq:dGdx_a}
\end{eqnarray}
Setting $z(h)(x) = (G(h)(x),h(G(h)(x)))$ the above equality can be rewritten in a simpler way as
\begin{eqnarray}
	\left(\frac{\partial f_x}{\partial x}(z(h)(x)) + \frac{\partial f_x}{\partial y}(z(h)(x)) Dh(G(h)(x)) \right)D(G(h))(x)=I. \label{eq:dGdx}
\end{eqnarray}
Observe that if $\|Dh(G(h)(x))\| \leq L$, then the condition $\xi>0$ implies that the matrix in the parethesis is invertible and we have 
\begin{equation}
	D(G(h))(x) =  \left(\frac{\partial f_x}{\partial x}(z(h)(x)) + \frac{\partial f_x}{\partial y}(z(h)(x)) Dh(G(h)(x)) \right) ^{-1}.  \label{eq:DxGhM-def}
\end{equation}
Let us differentiate the graph transform $\mathcal{T}$ with respect to $x$, we use  formula (\ref{eq:TG})
\begin{eqnarray}
	D (\mathcal{T}(h))(x)&=&   \frac{\partial f_y}{\partial x}(z(h)(x)) D(G(h))(x)  \nonumber \\
	& & +
	\frac{\partial f_y}{\partial y}(z(h)(x)) Dh(G(h)(x)) \cdot D(G(h))(x) \\ \label{eq:gt-diff}
	&=& \left( \frac{\partial f_y}{\partial x}(z(h)(x)) + \frac{\partial f_y}{\partial y}(z(h)(x)) Dh(G(h)(x))  \right)\cdot D(G(h))(x) \nonumber
\end{eqnarray}
We deduce that
$$
D (\mathcal{T}(h))(x) = \left( \frac{\partial f_y}{\partial x}(z(h)(x)) + \frac{\partial f_y}{\partial y}(z(h)(x)) Dh(G(h)(x))  \right)  \left(\frac{\partial f_x}{\partial x}(z(h)(x)) + \frac{\partial f_x}{\partial y}(z(h)(x)) Dh(G(h)(x)) \right) ^{-1}.
$$
In other words
$$
D (\mathcal{T}(h))(x)  \left(\frac{\partial f_x}{\partial x}(z(h)(x)) + \frac{\partial f_x}{\partial y}(z(h)(x)) Dh(G(h)(x)) \right) = \left( \frac{\partial f_y}{\partial x}(z(h)(x)) + \frac{\partial f_y}{\partial y}(z(h)(x)) Dh(G(h)(x))  \right).
$$

This motivates the implicit definition of extended graph transform $\mathcal{U}$ acting on $(h,M)$, where $h \in H$ and
$M:\overline{B}_u(0,r_u) \to \mbox{Lin}(\mathcal{X},\mathcal{Y})$ with the  $C^0$-norm
\begin{equation}\label{eq:u}
 \mathcal{U} (h,M)(x)\left(\frac{\partial f_x}{\partial x}(z(h)(x)) + \frac{\partial f_x}{\partial y}(z(h)(x)) M(G(h)(x)) \right) = \left( \frac{\partial f_y}{\partial x}(z(h)(x)) + \frac{\partial f_y}{\partial y}(z(h)(x)) M(G(h)(x))  \right).
\end{equation}
We have the following lemma that is a consequence of the implicit function theorem
\begin{lem}
	Let $h\in C^1(\overline{B}_u(0,r_u),\overline{B}_s(0,r_s))$ with $\|Dh\|\leq L$ and let $\mu \leq \xi$ and $\beta<1$. Then the graph transform $\mathcal{T}(h)$ is continuously differentiable and 
	$D(\mathcal{T}(h)) = \mathcal{U}(h,Dh).
	$
	\end{lem}
\begin{proof}
	The fact that $\mu \leq \xi$ imples that the matrix $\frac{\partial f_x}{\partial x}(z(h)(x))+\frac{\partial f_x}{\partial y}(z(h)(x))Dh(G(h))(x)$ is invertible. Then by the implicit function theorem $G(h)$ is differentiable with a derivative given by \eqref{eq:DxGhM-def}, and the assertion follows from differentiation of \eqref{eq:TG}. 
	\end{proof}
We will consider the mapping
$$
(h,M)\mapsto (\mathcal{T}(h),\mathcal{U}(h,M)),
$$
and we will prove that it is the fiber contraction.
%Observe  that from (\ref{eq:xi}) we have the following bound on $ D_xG (h,M)$
%\begin{equation}
%	\| D_xG (h,M)\| \leq \frac{1}{\xi}, \quad h \in H, \|M\| \leq L.  \label{eq:bound-DxG}
%\end{equation}

\subsubsection{A priori bound for $\mathcal{U}(h,M)$}

\begin{lem}\label{lem:res}
  Assume that $\mu\leq \xi$. If $h \in H$ and $\|M\| \leq L$, then $\|\mathcal{U}(h,M)\| \leq L$.
\end{lem}
\begin{proof}
Note that 
\begin{equation}\label{eq:ksi}
m\left(\frac{\partial f_x}{\partial x}(z(h)(x)) + \frac{\partial f_x}{\partial y}(z(h)(x)) M(G(h)(x))\right) \geq m\left(\frac{\partial f_x}{\partial x}\right) - L\left\|\frac{\partial f_x}{\partial y}\right\| = \xi.
\end{equation}
This means that 
\begin{equation}\label{eq:ksi_gie}
\left\|\left(\frac{\partial f_x}{\partial x}(z(h)(x)) + \frac{\partial f_x}{\partial y}(z(h)(x)) M(G(h)(x)) \right)^{-1}\right\| \leq \frac{1}{\xi}. 
\end{equation}
Therefore 
\begin{eqnarray*}
   \|\mathcal{U}(h,M)\| \leq \left( \left\|\frac{\partial f_y}{\partial x} \right\|   +  \left\|\frac{\partial f_y}{\partial y}\right\| \cdot \|M\| \right) \frac{1}{\xi}
   \leq L \left(\frac{1}{L}\left\|\frac{\partial f_y}{\partial x} \right\| + \left\|\frac{\partial f_y}{\partial y}\right\|  \right) \cdot  \frac{1}{\xi} \leq L \frac{\mu}{\xi} \leq L.
\end{eqnarray*}
The proof is complete.
\end{proof}

\subsubsection{A priori bound for the difference  $\mathcal{U}(h_1,M_1)-\mathcal{U}(h_2,M_2)$}

Denote 
$$
F(h,M)(x) = \left(\frac{\partial f_x}{\partial x}(z(h)(x)) + \frac{\partial f_x}{\partial y}(z(h)(x)) M(G(h)(x))) \right) ^{-1}
$$
In the first step we will estimate the difference between $F$ at two distinct points.  
\begin{lem}
\label{lem:LipDGdx-estm}
Let $\mu\leq \xi$. Assume that, for $i\in \{1,2\}$ we have $h_i
 \in H$ and $\|M_i\| \leq L$ and
\begin{equation}
  \|M_i(x_1) - M_i(x_2)\| \leq L_M \|x_1-x_2\|\, \ \textrm{for every}\ \ x_1,x_2\in \overline{B}_u(0,r_u).
\end{equation}
Then
\begin{equation}
 \left\|F(h_1,M_1)(x_1)- F(h_2,M_2)(x_2)\right\| \leq  \left(C_1  + \frac{L_M  \left\|\frac{\partial f_x}{\partial y}\right\|}{\xi^3} \right) \|x_1-x_2\| + C_2 \|h_1-h_2\| + \frac{1}{\xi^2}\left\|\frac{\partial f_x}{\partial y}\right\|\|M_1-M_2\|.
\end{equation}
where $C_1=C_1(N,f,Df,D^2f,L)$ does not depend on $L_M$, and $C_2=C_2(N,f,Df,D^2f,L, L_M)$.
\end{lem}
\begin{proof}
To shorten the notation we will write $z_i = z(h_i(x_i))$ and $G_i=G(h_i)(x_i)$. 
We first observe that using Lemma \ref{lem:G-estm}  
\begin{align}\label{eq:Mdiff}
	& \|M_1(G_1)-M_2(G_2)\|\leq L_M\|G_1-G_2\| + \|M_1-M_2\| \leq L_M\|G_1-G_2\| + \|M_1-M_2\| \nonumber \\
	& \ \leq L_M\frac{1}{\xi}\left\|\frac{\partial f_x}{\partial y}\right\|\|h_1-h_2\| + L_M \frac{1}{\xi}\|x_1-x_2\|+ \|M_1-M_2\|
\end{align}
From  the definition of $F$ it follows that 
\begin{align*}
&\left(\frac{\partial f_x}{\partial x}(z_1) + \frac{\partial f_x}{\partial y}(z_1)M_1(G_1) \right) F(h_1,M_1)(x_1)\\
& =
 \left(\frac{\partial f_x}{\partial x}(z_2) + \frac{\partial f_x}{\partial y}(z_2)M_2(G_2) \right)F(h_2,M_2)(x_2).
 \end{align*}
This means that
\begin{align*}
& \left(\frac{\partial f_x}{\partial x}(z_1) + \frac{\partial f_x}{\partial y}(z_1)M_1(G_1)   -\frac{\partial f_x}{\partial x}(z_2) + \frac{\partial f_x}{\partial y}(z_2)M_2(G_2) \right) F(h_2,M_2)(x_2)\\
& \ = \left(\frac{\partial f_x}{\partial x}(z_1) + \frac{\partial f_x}{\partial y}(z_1)M_1(G_1) \right)  \cdot \left( F(h_1,M_1)(x_1)- F(h_2,M_2)(x_2)\right)
\end{align*}

From (\ref{eq:ksi}) and assumption $\|M\| \leq L$ we have
\begin{eqnarray}\label{eq:estF}
  \xi^2 \left\| F(h_1,M_1)(x_1)- F(h_2,M_2)(x_2)\right\| \leq  
   \left\|\frac{\partial f_x}{\partial x}(z_1) - \frac{\partial f_x}{\partial x}(z_2)     \right\|   +\left\|
   \frac{\partial f_x}{\partial y}(z_1)M_1(G_1) - \frac{\partial f_x}{\partial y}(z_2)M_2(G_2) \right\|  
\end{eqnarray}
We estimate both terms separately, using \eqref{eq:gdiff} and Lemma \ref{lem:G-estm}
\begin{align}\label{diff}
  & \left\|\frac{\partial f_x}{\partial x}(z_1) - \frac{\partial f_x}{\partial x}(z_2)     \right\| =  \left\|\frac{\partial f_x}{\partial x}(G_1,h_1(G_1)) - \frac{\partial f_x}{\partial x}(G_2,h_2(G_2))     \right\| \\
 & \   \leq  \left\| \frac{\partial^2 f_x}{\partial x^2}\right\| \cdot \|G_1 - G_2\|  +  \left\| \frac{\partial^2 f_x}{\partial x \partial y}\right\| \cdot \|h_1(G_1)  - h_2(G_2) \|\nonumber\\
  &\ \leq \left( \left\| \frac{\partial^2 f_x}{\partial x^2}\right\|  +  \left\| \frac{\partial^2 f_x}{\partial x \partial y}\right\| L  \right) \cdot \|G_1 - G_2\| + \left\| \frac{\partial^2 f_x}{\partial x \partial y}\right\| \|h_1-h_2\|\nonumber\\
   & \leq \left( \left\| \frac{\partial^2 f_x}{\partial x^2}\right\|  +  \left\| \frac{\partial^2 f_x}{\partial x \partial y}\right\| L  \right)  \frac{1}{\xi} \|x_1 - x_2\| + \left(
   \left\| \frac{\partial^2 f_x}{\partial x \partial y}\right\|+\frac{1}{\xi}\left\|\frac{\partial f_x}{\partial y}\right\|
   \left( \left\| \frac{\partial^2 f_x}{\partial x^2}\right\|  +  \left\| \frac{\partial^2 f_x}{\partial x \partial y}\right\| L  \right)
   \right)
    \|h_1-h_2\|.\nonumber
\end{align}
Before we estimate the second term observe that analogous computations give
\begin{align*}
& \left\| \frac{\partial f_x}{\partial y}(z_1) - \frac{\partial f_x}{\partial y}(z_2) \right\| \\
& \  \leq \left( \left\| \frac{\partial^2 f_x}{\partial y \partial x}\right\|  +  \left\| \frac{\partial^2 f_x}{\partial y^2}\right\| L  \right)  \frac{1}{\xi} \|x_1 - x_2\| + \left(\left\| \frac{\partial^2 f_x}{ \partial y^2}\right\|+\frac{1}{\xi}\left\|\frac{\partial f_x}{\partial y}\right\|\left( \left\| \frac{\partial^2 f_x}{\partial y \partial x}\right\|  +  \left\| \frac{\partial^2 f_x}{\partial y^2}\right\| L  \right)\right) \|h_1-h_2\|.  
\end{align*}
We use this last bound to estimate the second term in \eqref{eq:estF}. 
\begin{align}\label{diffM}
 & \left\|
   \frac{\partial f_x}{\partial y}(z_1)M_1(G_1) - \frac{\partial f_x}{\partial y}(z_2)M_2(G_2) \right\|  \\
  &\ \leq  \left\|
   \frac{\partial f_x}{\partial y}(z_1)M_1(G_1) - \frac{\partial f_x}{\partial y}(z_1)M_2(G_2) \right\| 
   +  \left\|
   \frac{\partial f_x}{\partial y}(z_1)M_2(G_2) - \frac{\partial f_x}{\partial y}(z_2)M_2(G_2) \right\|\nonumber\\
 &  \leq \left\|\frac{\partial f_x}{\partial y}\right\| \cdot \|M_1(G_1) - M_2(G_2)\| + L \left\| \frac{\partial f_x}{\partial y}(z_1) - \frac{\partial f_x}{\partial y}(z_2) \right\| \nonumber \\
     & \leq \left\|\frac{\partial f_x}{\partial y}\right\| L_M\frac{1}{\xi}\left\|\frac{\partial f_x}{\partial y}\right\|\|h_1-h_2\| + L_M\left\|\frac{\partial f_x}{\partial y}\right\| \frac{1}{\xi}\|x_1-x_2\|+ \left\|\frac{\partial f_x}{\partial y}\right\| \|M_1-M_2\|\nonumber\\
    & +  L \left( \left\| \frac{\partial^2 f_x}{\partial y \partial x}\right\|  +  \left\| \frac{\partial^2 f_x}{\partial y^2}\right\| L  \right)  \frac{1}{\xi} \|x_1 - x_2\| + L \left(\left\| \frac{\partial^2 f_x}{ \partial y^2}\right\|+\frac{1}{\xi}\left\|\frac{\partial f_x}{\partial y}\right\|\left( \left\| \frac{\partial^2 f_x}{\partial y \partial x}\right\|  +  \left\| \frac{\partial^2 f_x}{\partial y^2}\right\| L  \right)\right) \|h_1-h_2\|.\nonumber
\end{align}
Combining the above estimates we obtain the assertion of the lemma.
\end{proof}

\begin{lem}
\label{lem:LipT2-estm}
Let $\mu\leq \xi$. Assume that, for $i\in \{1,2\}$ we have $h_i
\in H$ and $\|M_i\| \leq L$ and
\begin{equation}
	\|M_i(x_1) - M_i(x_2)\| \leq L_M \|x_1-x_2\|\, \ \textrm{for every}\ \ x_1,x_2\in \overline{B}_u(0,r_u).
\end{equation}
Then
\begin{equation}
 \left\|\mathcal{U}(h_1,M_1)(x_1)- \mathcal{U}(h_2,M_2)(x_2)\right\| \leq  \left(C_1 + \frac{\beta}{\xi^2} L_M \right) \cdot \|x_1-x_2\| + C_2\|h_1-h_2\|+\frac{\beta}{\xi}\|M_1-M_2\|. 
\end{equation}
where $C^1=C(N,f,Df,D^2f,L)$ does not depend on $L_M$.
\end{lem}
\begin{proof}
To shorten some formulas we will use the following notation $F_i=F(h_i,M_i)(x_i)$, $z_i=z(h_i)(x_i)$  and $G_i=G(h_i)(x_i)$ for $i\in \{1,2\}$. We will also denote by $C_1$ a generic constant dependent on $N, f, Df, D^2f, L$ and by $C_2$ a generic constant dependent on $N, f, Df, D^2f, L, L_M$. 
From the definition \eqref{eq:u} of $\mathcal{U}$ we have
\begin{align*}
  &\mathcal{U}(h_1,M_1)(x_1)- \mathcal{U}(h_2,M_2)(x_2) = \left(\frac{\partial f_y}{\partial x}(z_1) + \frac{\partial f_y}{\partial y}(z_1)M(G_1) \right) \cdot  (F_1 - F_2) \\
  &\ \ \ \ + \left( \left(\frac{\partial f_y}{\partial x}(z_1) - \frac{\partial f_y}{\partial x}(z_2) \right)
     + \left(\frac{\partial f_y}{\partial y}(z_1)M(G_1)
    - \frac{\partial f_y}{\partial y}(z_2)M(G_2) \right) \right) F_2
\end{align*}
For the first term from Lemma~\ref{lem:LipDGdx-estm}  we obtain the bound
\begin{align*}
 &  \left\|\left(\frac{\partial f_y}{\partial x}(z_1) + \frac{\partial f_y}{\partial y}(z_1)M(G_1) \right) \cdot  (F_1 - F_2)\right\|\\
 &  \leq 
  \left(\left\|\frac{\partial f_y}{\partial x}\right\| + \left\|\frac{\partial f_y}{\partial y} \right\| \cdot L \right) \cdot \left(\left(C_1  + \frac{L_M  \left\|\frac{\partial f_x}{\partial y}\right\|}{\xi^3} \right) \|x_1-x_2\| + C_2 \|h_1-h_2\| + \frac{1}{\xi^2}\left\|\frac{\partial f_x}{\partial y}\right\|\|M_1-M_2\|\right)\\
  &\leq \left( C_1 +  L_M  \frac{\mu L \left\|\frac{\partial f_x}{\partial y}\right\|}{\xi^3}  \right)\cdot \|x_1-x_2\| + C_2\|h_1-h_2\| +  \frac{\mu L}{\xi^2}\left\|\frac{\partial f_x}{\partial y}\right\|\|M_1-M_2\|. 
\end{align*}
where $C=C(N,f,Df,D^2f,L)$ does not depend on $L_M$.
We deal with the second term. Note that by  \eqref{eq:ksi_gie}  we have $\|F_2\|\leq \frac{1}{\xi}$. Moreover, analogously to \eqref{diff}

$$  \left\|\frac{\partial f_y}{\partial x}(z_1) - \frac{\partial f_y}{\partial x}(z_2) \right\| \leq C_1 \|x_1-x_2\|+C_1\|h_1-h_2\|.
  $$
  and
  $$  \left\|\frac{\partial f_y}{\partial y}(z_1) - \frac{\partial f_y}{\partial y}(z_2) \right\| \leq C_1 \|x_1-x_2\|+C_1\|h_1-h_2\|.
  $$
  We deal with the second term analogously as in \eqref{diffM}, namely
  \begin{align*}
  	& \left\|
  	\frac{\partial f_y}{\partial y}(z_1)M_1(G_1) - \frac{\partial f_y}{\partial y}(z_2)M_2(G_2) \right\|  \\
  	&\ \leq  \left\|
  	\frac{\partial f_y}{\partial y}(z_1)M_1(G_1) - \frac{\partial f_y}{\partial y}(z_1)M_2(G_2) \right\| 
  	+  \left\|
  	\frac{\partial f_y}{\partial y}(z_1)M_2(G_2) - \frac{\partial f_y}{\partial y}(z_2)M_2(G_2) \right\|\nonumber\\
  	&  \leq \left\|\frac{\partial f_y}{\partial y}\right\| \cdot \|M_1(G_1) - M_2(G_2)\| + L \left\| \frac{\partial f_y}{\partial y}(z_1) - \frac{\partial f_y}{\partial y}(z_2) \right\| \nonumber \\
  	& \leq C_2\|h_1-h_2\| + \left(C_1+L_M\left\|\frac{\partial f_y}{\partial y}\right\| \frac{1}{\xi}\right)\|x_1-x_2\|+ \left\|\frac{\partial f_y}{\partial y}\right\| \|M_1-M_2\|.\nonumber
  \end{align*} 
Combining all estimates leads us to the bound 
\begin{align*}
& \|\mathcal{U}(h_1,M_1)(x_1)- \mathcal{U}(h_2,M_2)(x_2)\| \leq \left( C_1 +  L_M  \frac{\mu L \left\|\frac{\partial f_x}{\partial y}\right\|}{\xi^3}  \right)\cdot \|x_1-x_2\| + C_2\|h_1-h_2\| \\
& \ \   +  \frac{\mu L}{\xi^2}\left\|\frac{\partial f_x}{\partial y}\right\|\|M_1-M_2\|+ \left(C_1+L_M\left\|\frac{\partial f_y}{\partial y}\right\| \frac{1}{\xi^2}\right)\|x_1-x_2\|+ \frac{1}{\xi}\left\|\frac{\partial f_y}{\partial y}\right\| \|M_1-M_2\|,
\end{align*}
which implies the assertion of the lemma.
\end{proof}
\subsubsection{A priori bound for the Lipschitz constant for $\mathcal{U}(h,M)$}
\begin{thm}
\label{thm:T2-a-priori-bnds}
  Let $\mu\leq \xi$ and assume that $\beta<\min\{1,\xi^2\}$. There exists a constant $L_M$ (depending on $N$, $f$, $Df$, $D^2f$ and $L$), such that  if $h \in H$  and $\|M\| \leq L $ and
   \begin{equation}
  \|M(x_1) - M(x_2)\| \leq L_M \|x_1-x_2\|\ \ \textrm{for every} \ x_1,x_2\in \overline{B}_u(0,r_u),  \label{eq:Lip-M}
\end{equation}
then
\begin{equation}
  \|\mathcal{U}(h,M)(x_1) - \mathcal{U}(h,M)(x_2)\| \leq L_M \|x_1-x_2\|\ \ \textrm{for every} \ x_1,x_2\in \overline{B}_u(0,r_u),  \label{eq:Lip-preserved}
\end{equation}
\end{thm}
\begin{proof}
From Lemma~\ref{lem:LipT2-estm}  it follows that it is enough to have
\begin{equation*}
  \left(C + \frac{\beta}{\xi^2} L_M \right)  \leq L_M.
\end{equation*}
Observe that $\frac{\beta}{\xi^2} < 1$.
Therefore it is enough to take
\begin{eqnarray*}
  L_M \geq \frac{C}{1 - \frac{\beta}{\xi^2} }.
\end{eqnarray*}
\end{proof}

\subsubsection{Graph transform  $(\mathcal{T}, \mathcal{U})$ for the unstable manifold and its derivative has an absorbing fixed point.}

\begin{thm}
\label{thm:DT-contr}
Let $\xi\geq \mu$ and let $\beta<1$. Assume that for $i\in \{1,2\}$ we have $h_i\in H$, $\|M_i\| \leq L$, and
\begin{eqnarray}
	\|M_i(x_1) - M_i(x_2)\| \leq L_M \|x_1 - x_2\|\ \ \textrm{for}\ \ x_1,x_2\in \overline{B}_u(0,r_u).
\end{eqnarray}
 Then there exists constant $C$ depending on $N, f, Df, D^2f, L,L_M$, such that
\begin{eqnarray*}
  \| \mathcal{T}_2(h_1,M_1) -  \mathcal{T}_2(h_2,M_2) \| \leq C \|h_1 - h_2\| + \frac{\beta}{\xi}\|M_1 - M_2\|.
\end{eqnarray*}
\end{thm}
\begin{proof}
	The result follows from Lemma \ref{lem:LipT2-estm} by taking $x_1=x_2$.
	\end{proof}

\begin{thm}\label{thm:15}
	Let $L_M$ be as in Theorem \ref{thm:T2-a-priori-bnds}. Assume that $\beta < \min\{1,\xi,\xi^2\}$ and $\xi\geq \mu$. 
	The mapping $(h,M)\mapsto (\mathcal{T}(h),\mathcal{U}(h,M))$ leads from the set  
	$$
	H \times \left\{ M\in C^0(\overline{B}_u(0,r_u);\mathrm{Lin}(\mathcal{X},\mathcal{Y}))\,:\ \|M\| \leq L,\ \  M\  \mathrm{is}\ L_M-\mathrm{Lipschitz} \right\},
	$$
	into itself and has the unique fixed point which is moreover attracting. 
\end{thm}
\begin{proof}
	The fact that the mapping $(\mathcal{T}, \mathcal{U})$ leads from the above set into itself is a straightforward consequence of Lemma \ref{incl1} and Lemma \ref{lem:res}, as well as Theorem \ref{thm:T2-a-priori-bnds}. The result follows from Theorem \ref{thm:fibre} by Theorem \ref{contr2} and Theorem \ref{thm:DT-contr}. 
\end{proof}

\subsection{Fixed point argument for construction of the stable manifold and its derivative}
\subsubsection{Graph transform for the stable manifold}

Now we will consider vertical cones satisfying cone condition.

\begin{df}
\label{def:vd-cc}
For a continuous map $x:\overline{B}_s(0,r_s) \to \overline{B}_u(0,r_u)$
we will say that $(x(y),y)$ is a \emph{vertical  disk satisfying cone condition} if
\begin{equation}
  \|x(y_1)-x(y_2)\| \leq  \frac{1}{L}\|y_1-y_2\|\ \ \textrm{for every}\ \ y_1,y_2\in \overline{B}_s(0,r_s)\label{eq:verdisk}
\end{equation}
\end{df}

\begin{df}
Let $V \subset C^{0}(\overline{B}_s(0,r_2),\overline{B}_u(0,r_u))$ be defined as follows: $v \in V$ if and only if $h$ is a vertical disk satisfying the cone condition.
\end{df}
Observe that $V$ is closed. Assume that $(x(y),y)$ is an stable manifold of $z_0$. Then for any $y$ there exists $y_0$ such that
\begin{equation}
  f(x(y),y)=(x(y_0),y_0).
\end{equation}
This is equivalent to
\begin{equation}
  f_x(x(y),y)=x(f_y(x(y),y)), \quad y_0=f_y(x(y),y).
\end{equation}
This suggests the following definition of the graph transform, given $v \in V$ we want $f^{-1}(v)$ parameterized as a vertical disk to be its graph transform. Therefore for given $y \in \overline{B}_s(0,r_s)$ we look for $x=\mathcal{S}(v)(y)$ such that point $f(\mathcal{S}(v)(y),y)$ belongs to image of  $v$, i.e. there exists $y_0$ such that
\begin{equation}
  f(\mathcal{S}(v)(y),y)=(v(y_0),y_0),
\end{equation}
which is equivalent to
\begin{equation}
  f_x(\mathcal{S}(v)(y),y)=v( f_y(\mathcal{S}(v)(y),y)).  \label{eq:def-S}
\end{equation}
This is an implicit definition of $\mathcal{S}(v)$. The fact that at least one $x\in \mathcal{S}(v)(y)$ that satisfies the above equation must exist follows from the fact that horizontal disks are mapped to horizontal disks. Uniqueness of $x$ as well as the fact that this graph transform maps vertical disks to vertical disks follows from the next lemma.

\begin{thm}
\label{thm:S-contrakcja}
Let $\mathcal{S}$ satisfy \eqref{eq:def-S}. If $\xi_1 >0$, then
for $v_1,v_2 \in V$ holds
\begin{equation}
  \|\mathcal{S}(v_1)(y_1) - \mathcal{S}(v_2)(y_2)\| \leq \frac{\|v_1- v_2\|}{\xi_1} + \frac{1}{L} \frac{\mu_1}{\xi_1}\|y_1-y_2\|.
\end{equation}
\end{thm}
\begin{proof}
Let us fix $v_1,v_2 \in V$ and  $y_1,y_2 \in \overline{B}_s(0,r_s)$. Let us denote $x_i=\mathcal{S}(v_i)(y_i)$.
Then
\begin{equation}
  f_x(x_i,y_i)=v_i(f_y(x_i,y_i))\ \ \textrm{for} \quad i=1,2.
\end{equation}
Hence, subtracting, we obtain
\begin{eqnarray*}
  f_x(x_1,y_1) - f_x(x_2,y_2)=v_1(f_y(x_1,y_1)) - v_2(f_y(x_2,y_2)).
\end{eqnarray*}
We have
$$
  \|f_x(x_1,y_1) - f_x(x_2,y_2)\| \geq  m\left( \frac{\partial f_x}{\partial x}\right) \|x_1 - x_2\| - \left\|\frac{\partial f_x}{\partial y}\right\|\|y_1 - y_2\|,
$$
and
\begin{align*}
 & \| v_1(f_y(x_1,y_1)) - v_2(f_y(x_2,y_2))\| \leq  \| v_1(f_y(x_1,y_1)) - v_1(f_y(x_2,y_2))\| +  \| v_1(f_y(x_2,y_1)) - v_2(f_y(x_2,y_2))\| \\
&   \leq \frac{1}{L}\|f_y(x_1,y_1) - f_y(x_2,y_2)\| + \|v_1 - v_2\|   \leq \frac{1}{L} \left\|\frac{\partial f_y}{\partial x}\right\| \cdot \|x_1-x_2\| +\frac{1}{L} \left\|\frac{\partial f_y}{\partial y}\right\| \cdot \|y_1-y_2\| + \|v_1 - v_2\|.
\end{align*}
Combining the above inequalities we obtain
\begin{eqnarray*}
   \left(m\left( \frac{\partial f_x}{\partial x}\right) -\frac{1}{L} \left\|\frac{\partial f_y}{\partial x}\right\|\right) \|x_1 - x_2\| \leq \frac{1}{L} \left( \left\|\frac{\partial f_y}{\partial y}\right\| + L\left\|\frac{\partial f_x}{\partial y}\right\|\right) \|y_1-y_2\|+ \|v_1 - v_2\|,
\end{eqnarray*}
which yields the assertion of the theorem.
\end{proof}

\begin{thm}\label{Sitself}
	Assume that $\mu_1\leq \xi_1$. Then 
	\begin{equation}
		\mathcal{S}(V) \subset V.
	\end{equation}
	\end{thm}
\begin{proof}
	Take $v\in V$. The topological argument implies that at least one $x = \mathcal{S}(v)(y)$ exists. 
	Its uniqueness and the fact that  $\mathcal{S}(v)$ is $\frac{1}{L}$-Lipschitz follows from Theorem \ref{thm:S-contrakcja} by taking $v_1=v_2=v$.
	\end{proof}

\begin{thm}\label{Scontraction}
	If $\xi_1>1$, then $\mathcal{S}$ is a contraction on $V$.
	\end{thm}
\begin{proof}
	The result follows by taking $y_1=y_2=y$ in Theorem \ref{thm:S-contrakcja}.
	\end{proof}

\subsubsection{Graph transform for the derivative of stable manifold}

Now we derive the equation for \\
$D \mathcal{S}(v)(y)=\frac{\partial \mathcal{S}(v)}{\partial y}(y)$.  We assume that $v \in C^1$ and differentiate (\ref{eq:def-S}). We obtain
\begin{align*}
&  \frac{\partial f_x}{\partial x}(\mathcal{S}(v)(y),y) D \mathcal{S}(y) +  \frac{\partial f_x}{\partial y}(\mathcal{S}(v)(y),y)\\
 & =      D v(f_y(\mathcal{S}(v)(y)),y)\left(\frac{\partial f_y}{\partial x}(\mathcal{S}(v)(y),y) D \mathcal{S}(y) +  \frac{\partial f_y}{\partial y}(\mathcal{S}(v)(y),y) \right)
\end{align*}
Let us define
\begin{equation}
z(v)(y)=(\mathcal{S}(v)(y),y).
\end{equation}
Observe that $z(v)(y) \in N$ for $y \in \overline{B}_s(0,r_s)$. Let $M=D v$. We can rewrite the above implicit equation as follows
\begin{align*}
  & \frac{\partial f_x}{\partial x}(z(v)(y)) D (\mathcal{S}(v))(y) +  \frac{\partial f_x}{\partial y}(z(v)(y))\\
   & =   M(f_y(z(v)(y)))\left(\frac{\partial f_y}{\partial x}(z(v)(y)) D (\mathcal{S}(v))(y) +  \frac{\partial f_y}{\partial y}(z(v)(y)) \right),
\end{align*}
which becomes
\begin{align*}
 & \left(\frac{\partial f_x}{\partial x}(z(v)(y)) -  M(f_y(z(v)(y))) \frac{\partial f_y}{\partial x}(z(v)(y)) \right) D (\mathcal{S}(v))(y)  \nonumber \\
   & =     M(f_y(z(v)(y)))\frac{\partial f_y}{\partial y}(z(v)(y)) -   \frac{\partial f_x}{\partial y}(z(v)(y)).
\end{align*}
Now we define the extended graph transform acting on $(v,M)$, where
$$v \in V\ \ \textrm{and}\ \   M \in C^0(\overline{B}_s(0,r_s),Lin(\mathcal{Y},\mathcal{X}))$$ by the formula
\begin{align}
&   \mathcal{R}(v,M)(y)= \left(\frac{\partial f_x}{\partial x}(z(v)(y)) -  M(f_y(z(v)(y))) \frac{\partial f_y}{\partial x}(z(v)(y)) \right)^{-1} \cdot  \label{eq:S2-def} \\
   & \qquad \qquad \qquad \qquad   \left( M(f_y(z(v)(y)))\frac{\partial f_y}{\partial y}(z(v)(y)) -   \frac{\partial f_x}{\partial y}(z(v)(y)) \right). \nonumber
\end{align}

\begin{lem}
 Assume that $v\in C^1(\overline{B}_s(0,r_s);\overline{B}_u(0,r_u))$  is such that $\xi_1 > 1$. Then $\mathcal{S}(v)$ is continuously differentiable and $D(\mathcal{S}(v)) = \mathcal{R}(v,Dv)$.
\end{lem}
\begin{proof} The fact that $\xi_1>1$ implies that the jacobian matrix $\frac{\partial f_x}{\partial x}(z(v)(y)) -  M(f_y(z(v)(y))) \frac{\partial f_y}{\partial x}(z(v)(y))$ is invertible for every $y\in \overline{B}_s(0,r_s)$. The assertion follows from the implicit function theorem.
	\end{proof}

\subsubsection{A priori bounds for $\mathcal{R}$}

\begin{lem}
\label{lem:S2-Mbnd}
Assume that $\xi_1 \geq \mu_1$ and $\xi_1>1$. If $v \in V$ and $\|M\| \leq \frac{1}{L}$, then $\mathcal{R}(v,M) \leq \frac{1}{L}$.
\end{lem}
\begin{proof} Denote for simplicity $z = z(v)(y)$. We have
	$$
		\left\|\left(\frac{\partial f_x}{\partial x}(z) -  M(f_y(z)) \frac{\partial f_y}{\partial x}(z) \right) \mathcal{R}(v,M)(y) \right\|=
		\left\| M(f_y(z))\frac{\partial f_y}{\partial y}(z) -   \frac{\partial f_x}{\partial y}(z)\right\|.
		$$
It follows that
	$$
m\left(\frac{\partial f_x}{\partial x}(z) -  M(f_y(z)) \frac{\partial f_y}{\partial x}(z) \right) \left\|\mathcal{R}(v,M)(y) \right\|\leq
\frac{1}{L}\left\| \frac{\partial f_y}{\partial y}\right\| +   \left\|\frac{\partial f_x}{\partial y}\right\|.
$$
We deduce
	$$
\left(m\left(\frac{\partial f_x}{\partial x}\right) -  \frac{1}{L} \left\|\frac{\partial f_y}{\partial x} \right\|\right) \left\|\mathcal{R}(v,M)(y) \right\|\leq
\frac{1}{L}\left\| \frac{\partial f_y}{\partial y}\right\| +   \left\|\frac{\partial f_x}{\partial y}\right\|.
$$
It means that
$$
\xi_1 \|\mathcal{R}(v,M)\| \leq  \frac{1}{L} \left(  \left\|\frac{\partial f_y}{\partial y}\right\| + L \left\|\frac{\partial f_x}{\partial y}\right\| \right)=\frac{\mu_1}{L},
$$
whence the assertion follows.
\end{proof}

\subsubsection{A priori bounds for the difference of two $\mathcal{R}$'s.}
\begin{thm}
	\label{lem:S-M-Lip_2}
	  Assume that, for $i\in \{1,2\}$ we have $v_i \in V$,
	\begin{equation}
		\|M_i\|\leq \frac{1}{L},
	\end{equation} and
	\begin{equation}
		\|M_i(y_1) - M_i(y_2)\| \leq L_M \|y_1 - y_2\|\ \ \textrm{for every}\ \ y_1,y_2\in \overline{B}_s(0,r_s).
	\end{equation}
	Then
	\begin{align*}
& \| \mathcal{R}(v_2,M_2)(y_2) - \mathcal{R}(v_1,M_1)(y_1)\| \leq 	\left(C_1 + L_M \frac{\mu}{\xi_1} \left(\left\|\frac{\partial f_y}{\partial y}\right\|+\frac{1}{L}\frac{\mu_1}{\xi_1}\left\|\frac{\partial f_y}{\partial x}\right\|\right)\right)\|y_1-y_2\| \\
	& \qquad + C_2\|v_1-v_2\| + \frac{1}{\xi_1}\left(\left\|\frac{\partial f_y}{\partial y}\right\|+\frac{1}{L}\frac{\mu_1}{\xi_1}\left\|\frac{\partial f_y}{\partial x}\right\|\right) \|M_1-M_2\|.
	\end{align*}
If, additionally, $\mu_1\leq \xi_1$, then
	\begin{equation}
		\| \mathcal{R}(v_2,M_2)(y_2) - \mathcal{R}(v_1,M_1)(y_1)\|  \leq \left(C_1+ L_M \frac{\mu^2}{\xi_1}\right)\|y_1-y_2\| + C_2\|v_1-v_2\| + \frac{\mu}{\xi_1} \|M_1-M_2\|,
	\end{equation}
	where  $C_1=C_1(N,f,Df,D^2f,L)$ does not depend on $L_M$ and $C_2=C_2(N,f,Df,D^2f,L,L_M)$.
\end{thm}
\begin{proof}
	To shorten some formulas let us denote $R_i=\mathcal{R}(v_i,M_i)(y_i)$ and   $z_i=z(v_i)(y_i)$ for $i=1,2$. Our point of departure is equation (\ref{eq:S2-def}) rewritten below
	for $i=1,2$ as an implicit equation
	\begin{eqnarray*}
		\left(\frac{\partial f_x}{\partial x}(z_i) -  M_i(f_y(z_i)) \frac{\partial f_y}{\partial x}(z_i) \right) R_i=
		M_i(f_y(z_i))\frac{\partial f_y}{\partial y}(z_i) -   \frac{\partial f_x}{\partial y}(z_i).
	\end{eqnarray*}
	Hence we obtain
	\begin{align*}
		& \left(\frac{\partial f_x}{\partial x}(z_2) -  M_2(f_y(z_2)) \frac{\partial f_y}{\partial x}(z_2) \right) R_2 - \left(\frac{\partial f_x}{\partial x}(z_1) -  M_1(f_y(z_1)) \frac{\partial f_y}{\partial x}(z_1) \right) R_1 \\
		&\ \ =  M_2(f_y(z_2))\frac{\partial f_y}{\partial y}(z_2) -   \frac{\partial f_x}{\partial y}(z_2) - M_1(f_y(z_1))\frac{\partial f_y}{\partial y}(z_1) +   \frac{\partial f_x}{\partial y}(z_1)
	\end{align*}
	Our aim is to derive the upper bound for $\|R_1-R_2\|$. From the above equation we obtain
	\begin{align} \label{eq:R2R1-estm}
		&\left(\frac{\partial f_x}{\partial x}(z_2) -  M_2(f_y(z_2)) \frac{\partial f_y}{\partial x}(z_2) \right) (R_2-R_1)  \\
	&\ \ 	 =  \left( \frac{\partial f_x}{\partial x}(z_1) -\frac{\partial f_x}{\partial x}(z_2)\right)R_1    +  \left(M_2(f_y(z_2)) \frac{\partial f_y}{\partial x}(z_2) -  M_1(f_y(z_1)) \frac{\partial f_y}{\partial x}(z_1)  \right) R_1 \nonumber \\
		&\ \ \ \ + M_2(f_y(z_2))\frac{\partial f_y}{\partial y}(z_2) - M_1(f_y(z_1))\frac{\partial f_y}{\partial y}(z_1) +   \frac{\partial f_x}{\partial y}(z_1)-   \frac{\partial f_x}{\partial y}(z_2) = I + II + III + IV \nonumber
	\end{align}
	
	For the lhs of (\ref{eq:R2R1-estm}) we have the estimate
	\begin{align*}
		&\left\| \left(\frac{\partial f_x}{\partial x}(z_2) -  M_2(f_y(z_2)) \frac{\partial f_y}{\partial x}(z_2) \right)(R_2-R_1) \right\|  \\
		&\ \ \geq m \left(\frac{\partial f_x}{\partial x}(z_2) -  M_2(f_y(z_2)) \frac{\partial f_y}{\partial x}(z_2) \right)\left\|R_2-R_1\right\|  \\
		&\ \ \geq \left( m \left(\frac{\partial f_x}{\partial x}\right) -   \frac{1}{L} \left\|\frac{\partial f_y}{\partial x}\right\| \right) \left\|R_2-R_1\right\| = \xi_1  \left\|R_2-R_1\right\|.
	\end{align*}
	In the estimates for the rhs of (\ref{eq:R2R1-estm})  we will have several expressions proportional either to $\|z_1 - z_2\|$, or to $\|M_i(f_y(z_2))-M_i(f_y(z_1))\|$, or to $\|M_1 - M_2\|$. It is important to us to get the explicit constants multiplying the last two terms.	We have
	\begin{align*}
		&\|z_1 - z_2 \| = \|z(v_1)(y_1) - z(v_2)(y_2)\| = \|(\mathcal{S}(v_1)(y_1)-\mathcal{S}(v_2)(y_2), y_1-y_2)\|\\
		&\leq \|\mathcal{S}(v_1)(y_1)-\mathcal{S}(v_1)(y_2)\| + \|\mathcal{S}(v_1)(y_2)-\mathcal{S}(v_2)(y_2)\|  + \|y_1-y_2\| \leq \frac{1}{\xi_1} \|v_1-v_2\| + \left(\frac{1}{L}+1\right)\|y_1-y_2\|,
	\end{align*}
	and
	\begin{align}  \label{eq:Mfy-diff}
	&	\|M_2(f_y(z_2))-M_1(f_y(z_1))\| \leq \|M_2(f_y(z_2))-M_1(f_y(z_2))\| + \|M_1(f_y(z_2))-M_1(f_y(z_1))\| \\
	& \leq \|M_2-M_1\| + L_M \|f_y(\mathcal{S}(v_1)(y_1),y_1) - f_y(\mathcal{S}(v_2)(y_2),y_2)\| \nonumber \\
		&\nonumber \leq \|M_1-M_2\| + L_M \left( \left\| \frac{\partial f_y}{\partial x} \right\|\cdot \|\mathcal{S}(v_1)(y_1)-\mathcal{S}(v_2)(y_2)\| + \left\| \frac{\partial f_y}{\partial y} \right\|\cdot  \|y_1 -y_2\|\right) \\
		&\nonumber  \leq  \|M_1-M_2\| + L_M \left(\left( \frac{1}{L}\left\| \frac{\partial f_y}{\partial x} \right\| + \left\| \frac{\partial f_y}{\partial y} \right\|\right)   \|y_1 -y_2\| + \frac{1}{\xi_1}\left\|\frac{\partial f_y}{\partial x}\right\|\|v_1-v_2\|\right)\\
		& = L_M \mu \|y_1 -y_2\| + \|M_1-M_2\| + \frac{L_M}{\xi_1}\left\|\frac{\partial f_y}{\partial x}\right\|\|v_1-v_2\|. \nonumber
	\end{align}
	We are in position to estimate all terms on rhs of  (\ref{eq:R2R1-estm}). We first estimate the term $I$, whence we obtain
	$$
	\|I\| \leq \left\|\frac{\partial f_x}{\partial x}(z_1) -\frac{\partial f_x}{\partial x}(z_2)\right\| \|R_1\| \leq \left(\left\|\frac{\partial^2 f_x}{\partial x^2}\right\| \left(\frac{1}{\xi_1}\|v_1-v_2\| + \frac{1}{L}\|y_1-y_2\|\right) + \left\|\frac{\partial^2 f_x}{\partial x \partial y}\right\| \|y_1-y_2\| \right) \frac{1}{L}\frac{\mu_1}{\xi_1}.
	$$
	Now we estimate the term $IV$. We get
	$$
	\|IV\| \leq \left\| \frac{\partial f_x}{\partial y}(z_1)-   \frac{\partial f_x}{\partial y}(z_2) \right\|
	\leq \left\| \frac{\partial^2 f_x}{\partial y\partial x}\right\|\left(\frac{1}{\xi_1}\|v_1-v_2\| + \frac{1}{L}\|y_1-y_2\|\right) +  \left\| \frac{\partial^2 f_x}{\partial y ^2}\right\|\|y_1-y_2\|. 	$$
	Next, we deal with the term $III$. We obtain
	\begin{align*}
	& 	\| III \| \leq \left \| M_2(f_y(z_2)) - M_1(f_y(z_1)) \right\| \left\| \frac{\partial f_y}{\partial y}(z_2)\right\| + \left\|M_1(f_y(z_1))\right\| \left\|\frac{\partial f_y}{\partial y}(z_2) - \frac{\partial f_y}{\partial y}(z_1) \right\| \\
	& \leq L_M \mu \left\|\frac{\partial f_y}{\partial y}\right\|\|y_1 -y_2\| + \left\|\frac{\partial f_y}{\partial y}\right\|\|M_1-M_2\| + \frac{L_M}{\xi_1}\left\|\frac{\partial f_y}{\partial x}\right\|\left\|\frac{\partial f_y}{\partial y}\right\| \|v_1-v_2\| \\
	& \ \ \ + \frac{1}{L}\left\| \frac{\partial^2 f_y}{\partial y\partial x}\right\|\left(\frac{1}{\xi_1}\|v_1-v_2\| + \frac{1}{L}\|y_1-y_2\|\right) +  \frac{1}{L}\left\| \frac{\partial^2 f_y}{\partial y ^2}\right\|\|y_1-y_2\|.
	\end{align*}
Finally, we estimate the last term $II$, whence
\begin{align*}
	 &\|II\| \leq \left\|M_2(f_y(z_2)) \frac{\partial f_y}{\partial x}(z_2) -  M_1(f_y(z_1)) \frac{\partial f_y}{\partial x}(z_1)  \right\| \|R_1\| \\
	  &\ \ \leq \frac{1}{L} \frac{\mu_1}{\xi_1}\left( \left\|M_2(f_y(z_2))  -  M_1(f_y(z_1))\right\| \left\| \frac{\partial f_y}{\partial x}(z_2)\right\| +  \|M_1(f_y(z_1))\| \left\|\frac{\partial f_y}{\partial x}(z_2) -   \frac{\partial f_y}{\partial x}(z_1)\right\|\right)\\
	 &\ \  \leq \frac{1}{L} \frac{\mu_1}{\xi_1}\left\|\frac{\partial f_y}{\partial x}\right\| L_M \mu \|y_1 -y_2\| + \frac{1}{L}\frac{\mu_1}{\xi_1} \left\|\frac{\partial f_y}{\partial x}\right\|\|M_1-M_2\| + \frac{L_M}{\xi_1}\frac{1}{L}\frac{\mu_1}{\xi_1} \left\|\frac{\partial f_y}{\partial x}\right\|^2\|v_1-v_2\|\\
	 & \ \ \ + \frac{1}{L^2} \frac{\mu_1}{\xi_1}\left(\left\|\frac{\partial^2 f_y}{\partial x^2}\right\| \left(\frac{1}{\xi_1}\|v_1-v_2\| + \frac{1}{L}\|y_1-y_2\|\right) + \left\|\frac{\partial^2 f_y}{\partial x \partial y}\right\| \|y_1-y_2\| \right).	
	\end{align*}
Adding all four estimates we obtain
\begin{align*}
& \xi_1 \|R_1-R_2\| \leq \left(D_1(N,f,Df,D^2f,L) + L_M \mu \left(\left\|\frac{\partial f_y}{\partial y}\right\|+\frac{1}{L}\frac{\mu_1}{\xi_1}\left\|\frac{\partial f_y}{\partial x}\right\|\right)\right)\|y_1-y_2\| \\
& \qquad + D_2(N,f,Df,D^2f,L,L_M)\|v_1-v_2\| + \left(\left\|\frac{\partial f_y}{\partial y}\right\|+\frac{1}{L}\frac{\mu_1}{\xi_1}\left\|\frac{\partial f_y}{\partial x}\right\|\right) \|M_1-M_2\|,
\end{align*}
which implies the assertion. 	\end{proof}

\subsubsection{A priori bounds for Lipschitz constant for $\mathcal{R}(v,M)$}

\begin{lem}
\label{lem:S-M-Lip}
Assume that $\mu_1\leq \xi_1$.  If $v \in V$, $\|M\| \leq \frac{1}{L}$ and
\begin{equation}
  \|M(y_1) - M(y_2)\| \leq L_M \|y_1 - y_2\|,
\end{equation}
then
\begin{equation}
  \| \mathcal{R}(v,M)(y_1) - \mathcal{R}(v,M)(y_2)\|  \leq \left( C + L_M \frac{\mu^2}{\xi_1} \right) \|y_1 - y_2\|,
\end{equation}
where  $C=C(N, f,Df,D^2f,L)$ does not depend on $L_M$.
\end{lem}
\begin{proof} The result follows by taking $v_1=v_2=v$ and $M_1=M_2=M$ in Lemma \ref{lem:S-M-Lip_2}.
\end{proof}

\begin{thm}
\label{thm:S2-lip-forM}
Assume that $\xi_1 > \max\{ 1, \mu^2\}$ and $\xi_1\geq \mu_1$. There exists a constant $L_M$(depending on $N$, $f$, $Df$, $D^2f$ and $L$), such that  if $v \in V$  and $\|M\| \leq \frac{1}{L} $ and
   \begin{equation}
  \|M(y_1) - M(y_2)\| \leq L_M \|y_1-y_2\|,  \label{eq:S-Lip-M}
\end{equation}
then
\begin{equation}
  \|\mathcal{R}(h,M)(y_1) - \mathcal{R}(h,M)(y_2)\| \leq L_M \|y_1-y_2\|,  \label{eq:S-Lip-preserved}
\end{equation}
\end{thm}
\begin{proof}
We use Lemma~\ref{lem:S-M-Lip}.  It is easy to see that we can take any $L_M$ satisfying
\begin{equation*}
  L_M \geq \frac{C}{1 - \frac{\mu^2}{\xi_1}}.
\end{equation*}
\end{proof}

\subsubsection{Graph transform  $(\mathcal{S}, \mathcal{R})$ for the stable manifold and its derivative has an absorbing fixed point.}
\begin{thm}\label{contraction_stable}
Let $\xi_1 \geq \mu_1$. Assume that  $v_1,v_2 \in V$ and $\|M_1\| \leq \frac{1}{L}$, $\|M_2\| \leq \frac{1}{L}$ and $L_M$ be as in Theorem~\ref{thm:S2-lip-forM} and
\begin{eqnarray}
  \|M_i(y_1) - M_i(y_2)\| \leq L_M \|y_1 - y_2\|\ \ \textrm{for}\ \ i\in \{1,2\}.
\end{eqnarray}
Then there exists a constant $C$ depending on  $f$, $Df$, $D^2f$ (restricted to $N$) and $L$ and $L_M$, such that
\begin{equation}
  \|\mathcal{R}(v_1,M_1) - \mathcal{R}(v_2,M_2)\| \leq C \|v_1 - v_2\| + \frac{\mu}{\xi_1} \|M_1 - M_2\|.
\end{equation}
\end{thm}
\begin{proof}
	The result follows from Lemma \ref{lem:S-M-Lip_2} by taking  $y_1 = y_2.$
	\end{proof}

\begin{thm}\label{thm:25}
	Let $L_M$ be as in Theorem \ref{thm:S2-lip-forM}. Assume that $\xi_1 \geq \mu_1$ and $\xi_1 > \max\{1, \mu, \mu^2\}$.
	The mapping $(v,M)\mapsto (\mathcal{S}(v),\mathcal{R}(v,M))$ leads from the set
$$
V \times \left\{ M\in C^0(\overline{B}_s(0,r_s);\mathrm{Lin}(\mathcal{Y},\mathcal{X}))\,:\ \|M\| \leq \frac{1}{L},\ \  M\ \  \mathrm{is}\ L_M-\mathrm{Lipschitz} \right\},
$$
into itself and has the unique fixed point which is moreover attracting.
\end{thm}
\begin{proof}
	The fact that the mapping $(\mathcal{S}, \mathcal{R})$ leads from the above set into itself is a straightforward consequence of Lemma \ref{lem:S2-Mbnd} and Theorem \ref{thm:S2-lip-forM}, as well as Theorem \ref{Sitself}. The result follows from Theorem \ref{thm:fibre} by Theorem \ref{contraction_stable} and Theorem \ref{Scontraction}.
\end{proof}

\section{Appendix 5: Verification of conditions from Appendix 4.}
In this section we work in local coordinates in the isolating set with cones, we denote these coordinates as $(y_s,y_u)$, where the unstable variable is $y_u$ and the stable one is $y_s$. We need to verify the conditions of Appendix 4, namely that
\begin{align*}
	& (1)\ \ \ m\left(\frac{\partial f_u}{\partial y_u}\right) > 1,\ \  \left\|\frac{\partial f_s}{\partial (y_s,\eta)}\right\| < 1,\\
& (2) \ \ \ \left\|\frac{\partial f_u}{\partial (y_s,\eta)}\right\|\ \ \textrm{can be made arbitrarily small by decreasing, if necessary, the set}\ N\ \textrm{and}\ \varepsilon,
\end{align*}
where $f$ is the mapping that assigns to the initial data the solution after a given time and the derivatives are understood with respect to the initial data.

The equation which we are solving has the following form in the local coordinates
$$
y'(t) = h(y(t)) + \varepsilon T_{\kappa}^{-1} \left(\int_0^\infty M(s)\, ds\right) (x_0+T_\kappa y(t)) + \varepsilon T_\kappa^{-1}\int_0^\infty M(s) \eta^t(s)\, ds.
$$
with
$$
h(y) = T_\kappa^{-1}D f(x_0) T_\kappa y + T_\kappa^{-1} f(x_0 + T_\kappa y) - T_\kappa^{-1}D f(x_0) T_\kappa y
$$
The variable $\eta^t$ is evolving according to the rule
$$
\eta^t(s) = \begin{cases} T_k(y(t-s) - y(t))\ \ \textrm{for}\ \ s\leq t\\
	T_k(y(t-s) - y(t))=T_k y_0 +\eta^0(s-t) - T_\kappa y(t)\ \ \textrm{otherwise}.
\end{cases}
$$
We use Lemma \ref{variational} by which  the derivative of the solution with respect to the initial data is given by the solution of the variational problem, which, after the change of variables to the local variables in the isolating set $N$ has the form
\begin{align}
	& w'(t) = T_\kappa^{-1}D f(x_0) T_\kappa w(t) + T_\kappa^{-1} (Df(x_0 + T_\kappa y(t)) - D f(x_0)) T_\kappa w(t) \nonumber\\
	& \ \ \ \qquad \qquad + \varepsilon T_{\kappa}^{-1} \left(\int_0^\infty M(s)\, ds\right) T_\kappa w(t) + \varepsilon T_\kappa^{-1}\int_0^\infty M(s) \theta^t(s)\, ds.\label{weq}
\end{align}
$$
\theta^t(s) = \begin{cases} T_k(w(t-s) - w(t))\ \ \textrm{for}\ \ s\leq t\\
	T_\kappa w_0+\xi^0(s-t) - T_\kappa w(t)\ \ \textrm{otherwise},
\end{cases}
$$
where $(\xi^0,w_0)$ are the initial data.
 We rewrite \eqref{weq} as
\begin{align}
	& w'(t) = T_\kappa^{-1}D f(x_0) T_\kappa w(t) + T_\kappa^{-1} (Df(x_0 + T_\kappa y(t)) - D f(x_0)) T_\kappa w(t) \nonumber\\
	& \ \ \  + \varepsilon T_\kappa^{-1}\int_0^t M(s) T_\kappa w(t-s)\, ds  + \varepsilon T_\kappa^{-1}\int_t^\infty M(s) \, dsT^\kappa w_0 + \varepsilon T_\kappa^{-1}\int_t^\infty M(s) \xi^0(s-t)\, ds.\label{eqw1}
\end{align}
We can further rewrite the above equation as
\begin{align}
	& w'(t) =  T_\kappa^{-1} Df(x_0 + T_\kappa y(t))  T_\kappa w(t) \nonumber\\
	& \ \ \  + \varepsilon T_\kappa^{-1}\int_0^t M(t-s) T_\kappa w(s)\, ds  + \varepsilon T_\kappa^{-1}\int_t^\infty M(s) \, dsT^\kappa w_0 + \varepsilon T_\kappa^{-1}\int_0^\infty M(s+t) \xi^0(s)\, ds. \label{106}
\end{align}
Assume that $t\in [0,1]$. It follows that
$$
|w(t)| \leq C |w(0)|+\varepsilon C\|\xi^0\| + C\int_0^t |w(s)|\, ds. \
$$
So, the Gronwall lemma implies that
\begin{equation}\label{107}
|w(t)|\leq Ce^{Ct} (|w(0)|+\varepsilon \|\xi^0\|).
\end{equation}
This also implies that
\begin{equation}\label{107b}
\int_0^t|w(s)|\, ds \leq C e^{Ct} (|w(0)|+\varepsilon \|\xi^0\|).
\end{equation}
We project \eqref{eqw1} on the stable and unstable direction of $w$, whence we get the following two equations
\begin{align}
	& w'_s(t) = (T_\kappa^{-1}D f(x_0) T_\kappa)_s w_s(t) + \Pi_sT_\kappa^{-1} (Df(x_0 + T_\kappa y(t)) - D f(x_0)) T_\kappa w(t) \nonumber\\
	& \ \ \  + \varepsilon \Pi_sT_\kappa^{-1}\int_0^t M(s) T_\kappa w(t-s)\, ds  + \varepsilon \Pi_sT_\kappa^{-1}\int_t^\infty M(s) \, dsT^\kappa w(0) + \varepsilon \Pi_sT_\kappa^{-1}\int_t^\infty M(s) \xi^0(s-t)\, ds.\label{eqw1_s}
\end{align}

\begin{align}
	& w'_u(t) = (T_\kappa^{-1}D f(x_0) T_\kappa)_u w_u(t) + \Pi_uT_\kappa^{-1} (Df(x_0 + T_\kappa y(t)) - D f(x_0)) T_\kappa w(t) \nonumber\\
	& \ \ \  + \varepsilon \Pi_uT_\kappa^{-1}\int_0^t M(s) T_\kappa w(t-s)\, ds  + \varepsilon \Pi_uT_\kappa^{-1}\int_t^\infty M(s) \, dsT_\kappa w(0) + \varepsilon \Pi_uT_\kappa^{-1}\int_t^\infty M(s) \xi^0(s-t)\, ds.\label{eqw1_u}
\end{align}
We first study \eqref{eqw1_u} to verify
the first condition of (1) and the condition (2).
\begin{lem}
	There exists $\varepsilon_0>0$ and the isolating block with cones  $N$ such that for every $\varepsilon\in [0,\varepsilon_0]$ on the block $N$  we have  $m\left(\frac{\partial f_u}{\partial y_u}\right) > 1$. Moreover
	$\left\|\frac{\partial f_u}{\partial (y_s,\eta)}\right\|$ can be made arbitrarily small by decreasing, if necessary, the block $N$ and value $\varepsilon$.
	\end{lem}
	\begin{proof}
From \eqref{eqw1_u} we obtain
\begin{align}
	& \frac{d}{dt}|w_u(t)| \geq  m\left((T_\kappa^{-1}D f(x_0) T_\kappa)_u\right) |w_u(t)| - C\delta^2|w(t)|  - \varepsilon C\int_0^t |w(s)|\, ds  - \varepsilon C |w(0)| - \varepsilon C\|\xi^0\|.
\end{align}
Furthermore,
\begin{align}
	& \frac{d}{dt}|w_u(t)| \geq  m\left((T_\kappa^{-1}D f(x_0) T_\kappa)_u\right) |w_u(t)| - C(\delta^2+\varepsilon)e^{Ct}(|w(0)|+\|\xi^0\|).
\end{align}
We denote $m\left((T_\kappa^{-1}D f(x_0) T_\kappa)_u\right) = \lambda_1>0$, hence
$$
\frac{d}{dt}|w_u(t)| - \lambda_1 |w_u(t)| \geq   - C(\delta^2+\varepsilon)e^{Ct}(|w(0)|+\|\xi^0\|)
$$
We estimate $t$ in $e^{Ct}$ by $T$ and we multiply by $e^{-\lambda_1 t}$
$$
e^{-\lambda_1 t}\frac{d}{dt}|w_u(t)| - e^{-\lambda_1 t}\lambda_1 |w_u(t)| \geq   - e^{-\lambda_1 t} Ce^{CT}(\delta^2+\varepsilon)(|w(0)|+\|\xi^0\|)
$$
$$
\frac{d}{dt}e^{-\lambda_1 t}|w_u(t)|\geq   - e^{-\lambda_1 t} Ce^{CT}(\delta^2+\varepsilon)(|w(0)|+\|\xi^0\|)
$$
We integrate from $0$ to $T$, whence
$$
e^{-\lambda_1T}|w_u(T)| - |w_u(0)| \geq  -\frac{1}{\lambda_1}(1-e^{-\lambda_1T}) Ce^{CT}(\delta^2+\varepsilon)(|w(0)|+\|\xi^0\|)
$$
It follows that
$$
|w_u(T)|  \geq  e^{\lambda_1 T}|w_u(0)|  - \frac{e^{(\lambda_1+C)T} C}{\lambda_1}(\delta^2+\varepsilon)(|w(0)|+\|\xi^0\|)
$$
Now, if $w(0) = w_u(0)$ and $\xi^0=0$, then
\begin{equation}\label{eq:cone1}
|w_u(T)|  \geq  \left(e^{\lambda_1 T}-\frac{e^{(\lambda_1+C)T} C}{\lambda_1}(\delta^2+\varepsilon)\right)|w_u(0)|,
\end{equation}
and it is possible to choose $\delta$ and $\varepsilon$ small enough to get the constant in front of $|w_u(0)|$ greater then one. This verifies the first assertion of (1).

On the other hand, coming back to \eqref{eqw1_u}, for a constant $\lambda_2= \left\|(T_\kappa^{-1}D f(x_0) T_\kappa)_u\right\|$ we obtain
$$
\frac{d}{dt}|w_u(t)| \leq \lambda_2 |w_u(t)|+C(\delta^2+\varepsilon)e^{CT}(|w(0)|+\|\xi^0\|).
$$
The Gronwall lemma implies that
$$
|w_u(T)| \leq e^{\lambda_2 T}|w_u(0)|+\frac{e^{(\lambda_2+C)T}}{\lambda_2}C(\delta^2+\varepsilon)(|w(0)|+\|\xi^0\|).
$$
Now, if $w_u(0) = 0$, we obtain
\begin{equation}\label{eq:cone2}
|w_u(T)| \leq \frac{e^{(\lambda_2+C)T}}{\lambda_2}C(\delta^2+\varepsilon)(|w_s(0)|+\|\xi^0\|).
\end{equation}
Equations \eqref{eq:cone1} and \eqref{eq:cone2}
 verify the first assertion of (1) and the condition (2). Indeed, no matter how large $T$ we take we can always find small $\delta$ and $\varepsilon$ such that these assertions hold.
\end{proof}
In the next result we study the variational equation for the stable variable in order to obtain the second assertion from (1). Here, we also need to take into account the evolution of the memory variable.
\begin{lem}
	There exists $\varepsilon_0>0$ and the isolating block with cones  $N$ such that for every $\varepsilon\in [0,\varepsilon_0]$ on the block $N$ we have $\left\|\frac{\partial f_s}{\partial (y_s,\eta)}\right\| < 1$.
\end{lem}
\begin{proof}
For the stable part of \eqref{weq} we denote $\mu\left((T_\kappa^{-1}D f(x_0) T_\kappa)_s\right) = - \lambda_3 < 0$. Hence, \eqref{eqw1_s} implies
\begin{equation}\label{ws}
\frac{d}{dt}|w_s(t)| \leq - \lambda_3 | w_s(t)| + Ce^{CT}(\delta^2+\varepsilon)(|w(0)|+\|\xi^0\|).
\end{equation}
In order to deal with the history variable $\theta$ note that, as in Lemma \ref{lem:eta}, we have
$$
\frac{d}{dt}\|\theta ^t\|^2+C\|\theta^t\|^2 \leq -2\left(\int_0^\infty A(s)\theta^t(s)\, ds,T_kw'(s)\right).
$$
Using \eqref{106} this implies that
$$
\frac{d}{dt}\|\theta ^t\|^2+C\|\theta^t\|^2 \leq C_1\|\theta^t\|\, |w'(t)| \leq C_1\|\theta^t\| \left(|w(t)|+\varepsilon\int_0^t|w(s)| ds+\varepsilon|w_0|+\varepsilon\|\xi^0\|\right),
$$
or
$$
\frac{d}{dt}\|\theta ^t\|\leq - C\|\theta^t\|+ C_1 \left(|w_s(t)|+|w_u(t)|+\varepsilon\int_0^t|w(s)| ds ds+\varepsilon|w_0|+\varepsilon\|\xi^0\|\right).
$$
Using \eqref{107b} this means that
$$
\frac{d}{dt}\|\theta ^t\|\leq - C\|\theta^t\|+ C_1 \left(|w_s(t)|+|w_u(t)|+\varepsilon e^{CT}(|w_0|+\|\xi^0\|)\right). $$
Taking a linear combination of this equation with \eqref{ws} we obtain
$$
\frac{d}{dt}(\|\theta ^t\|+K|w_s(t)|)\leq - C\|\theta^t\|+ (C_1-K\lambda_3)|w_s(t)|+C_1|w_u(t)|+(\varepsilon+\delta^2) Ce^{CT}(|w_0|+\|\xi^0\|). $$
We take $K$ such that $C_1-K\lambda_3<0$. Then for some constant $D>0$ we have
$$
\frac{d}{dt}(\|\theta ^t\|+K|w_s(t)|)\leq - D(\|\theta^t\|+ K|w_s(t)|)+C_1|w_u(t)|+(\varepsilon+\delta^2) Ce^{CT}(|w_0|+\|\xi^0\|). $$
First we take $w_u(0)=0$. Then
$$
\frac{d}{dt}(\|\theta ^t\|+K|w_s(t)|)\leq - D(\|\theta^t\|+ K|w_s(t)|)+(\varepsilon+\delta^2) Ce^{CT}(|w_s(0)|+\|\xi^0\|). $$
After application of the Gronwall lemma we obtain
$$
\|\theta ^T\|+K|w_s(T)|\leq e^{-DT}(\|\xi^0\|+K|w_s(0)|)+(\varepsilon+\delta^2) Ce^{CT}(|w_s(0)|+\|\xi^0\|).
$$
This means that for a given $T$ we can find $\varepsilon$ and $\delta$ small enough such that the second assertion of (1) is satisfied.
\end{proof}

Finally let us take $\xi^0=0$ and $w_s(0)=0$. This leads to the estimate of the value of $\left\|\frac{\partial f_s}{\partial y_u}\right\|$ which corresponds to $\left\|\frac{\partial f_y}{\partial x}\right\|$ present in the constants $\mu$ in \eqref{eq:mu} and $\xi_1$ in \eqref{eq:xi1}. Note that this quantity does not have to be small, it only needs to be bounded on the isolating block. Conditions that $\mu<1$ and $\xi_1>1$ are guaranteed by the selection of appropriately large $L$.   We obtain
$$
\frac{d}{dt}(\|\theta ^t\|+K|w_s(t)|)\leq - D(\|\theta^t\|+ K|w_s(t)|)+e^{\lambda_2T}|w_u(0)|+(\varepsilon+\delta^2) Ce^{CT}|w_u(0)|. $$
This means that
$$
\|\theta ^t\|+K|w_s(t)|\leq e^{\lambda_2T}C|w_u(0)|+(\varepsilon+\delta^2) Ce^{CT}|w_u(0)|.
$$

\section{Appendix 6. Continuous dependence of derivatives of stable and unstable manifolds on parameter.}

We verify the conditions of Theorem \ref{thm:param}, namely that the graph transform mappings for stable and unstable manifolds are continuous functions with respect to parameter $\varepsilon$. This will yield the assertion that their fixed points,  stable and unstable manifolds, are $C^1$ continuous functions of $\varepsilon$. Specifically we need to show that the mappings
$$
(\varepsilon,h,M)\mapsto (\mathcal{T}(\varepsilon,h),\mathcal{U}(\varepsilon,h,M)),
$$
and
$$
(\varepsilon,v,M)\mapsto (\mathcal{S}(\varepsilon,v),\mathcal{R}(\varepsilon,v,M)),
$$
are continuous. The arguments is analogous to the arguments of Appendix 4, we need consider mappings with additional dependence on $\varepsilon$, namely $(f_x^\varepsilon)(x,y),f_y^\varepsilon(x,y))$, hence in all estimates we obtain extra terms depending of the difference $f^{\varepsilon_1}-f^{\varepsilon_2}$ or its derivatives. As the derivations of the estimates  closely follow the lines of the ones from Appendix 4, we skip the proofs, presenting only the results.  We make the standing assumptions that for every $(x,y)\in N$ and every $\varepsilon_1, \varepsilon_2\in [0,\varepsilon_0]$ we have
$$
\|f_x^{\varepsilon_1}(x,y)-f_x^{\varepsilon_2}(x,y)\|\leq K |\varepsilon_1-\varepsilon_2|,
$$
and
$$
\|f_y^{\varepsilon_1}(x,y)-f_y^{\varepsilon_2}(x,y)\|\leq K |\varepsilon_1-\varepsilon_2|,
$$
moreover
$$
\left\|\frac{\partial f_y^{\varepsilon_1}(x,y)}{\partial x}-\frac{\partial f_y^{\varepsilon_2}(x,y)}{\partial y}\right\|\leq K |\varepsilon_1-\varepsilon_2| \ \ \textrm{and}\ \ \left\|\frac{\partial f_y^{\varepsilon_1}(x,y)}{\partial x}-\frac{\partial f_y^{\varepsilon_2}(x,y)}{\partial y}\right\|\leq K |\varepsilon_1-\varepsilon_2|,
$$
and
$$
\left\|\frac{\partial f_x^{\varepsilon_1}(x,y)}{\partial x}-\frac{\partial f_x^{\varepsilon_2}(x,y)}{\partial y}\right\|\leq K |\varepsilon_1-\varepsilon_2| \ \ \textrm{and}\ \ \left\|\frac{\partial f_x^{\varepsilon_1}(x,y)}{\partial x}-\frac{\partial f_x^{\varepsilon_2}(x,y)}{\partial y}\right\|\leq K |\varepsilon_1-\varepsilon_2|,
$$
with a constant $K>0$. The estimates for the difference of functions follow from Lemma \ref{lem:difference} and for the difference of derivatives follow from Lemma \ref{lem_diff}.

Constants $\xi, \mu, \beta, \xi_1, \mu_1$ now depend on $\varepsilon$. We will denote the new constants as $\xi^\varepsilon, \mu^\varepsilon, \beta^\varepsilon, \xi_1^\varepsilon, \mu_1^\varepsilon$. Arguments of Appendix 5 demonstrate that the bounds \eqref{eq:bounds_constants} hold independently on $\varepsilon$, and moreover $\xi^\varepsilon >0$ and $\xi_1^{\varepsilon}>0$ for every $\varepsilon\in [0,\varepsilon_0]$. These bounds are used in the proofs of the results in the following parts of this section.

\subsection{Graph transform for the unstable manifold.}
The arguments of this section are obtained analogously to the proofs of Section \ref{fix_unstable}.
The mapping $\mathcal{T}$, the graph transform with parameter, is now given by
\begin{equation}
	\mathcal{T}(\varepsilon,h) (x) = f^\varepsilon_y(G(\varepsilon,h)(x),h(G(\varepsilon,h)(x))),  \label{eq:TGe}
\end{equation}
with $G$ given as $G(\varepsilon,h)(x) = \overline{x}$ such that $x = f^\varepsilon_x(\overline{x},h(\overline{x}))$. Proceeding analogously as in the proof of Lemma \ref{lem:G-estm} we obtain the next result
\begin{lem}
	\label{lem:G-estm_e}
	Let $\xi^{\varepsilon_1}>0$ . Then, assuming that $h_1,h_2 \in H$, we have
	\begin{equation}
		\|G(\varepsilon_1,h_1)(x) - G(\varepsilon_2,h_2)(x)\| \leq \frac{K}{\xi^{\varepsilon_1}}|\varepsilon_1-\varepsilon_2| +      \frac{1}{\xi^{\varepsilon_1}} \left\| \frac{\partial f^{\varepsilon_1}_x}{\partial y}  \right\|  \|h_1 - h_2\|.
	\end{equation}
\end{lem}
\begin{proof}
Let us fix $x \in \overline{B}_u(0,r_u)$ and  let us denote $\overline{x}_i=G(\varepsilon_i,h_i)(x)$.  By definition of $G$ we have $f^{\varepsilon_i}_x(\overline{x}_i,h_i(\overline{x}_i))=x$, hence
\begin{align*}
	& 0= \|f_x^{\varepsilon_1}(\overline{x}_1,h_1(\overline{x}_1)) - f^{\varepsilon_2}_x(\overline{x}_2,h_2(\overline{x}_2))\|\\
	&\geq  \left(\frac{\partial f^{\varepsilon_1}_x}{\partial x}\right) \|\overline{x}_1 - \overline{x}_2\| - \left\| \frac{\partial f^{\varepsilon_1}_x}{\partial y}  \right\| \cdot \|h_1(\overline{x}_1) - h_2(\overline{x}_2)\|-\|f_x^{\varepsilon_1}(\overline{x}_2,h_1(\overline{x}_2)) - f^{\varepsilon_2}_x(\overline{x}_2,h_2(\overline{x}_2))\|\\
	& \ \geq \left(\frac{\partial f^{\varepsilon_1}_x}{\partial x}\right) \|\overline{x}_1 - \overline{x}_2\| - \left\| \frac{\partial f^{\varepsilon_1}_x}{\partial y}  \right\| \cdot (\|h_1-h_2\|+L\|\overline{x}_1-\overline{x}_2\|)-K|\varepsilon_1-\varepsilon_2|,
\end{align*}
and the assertion follows exactly as in Lemma \ref{lem:G-estm}.
\end{proof}
The proof of the next result follows the lined of the proof of Theorem \ref{thm:T-contr}.
\begin{thm}
	\label{thm:T-contr_e}
For any $h_1,h_2 \in H$ and $x \in \overline{B}_u(0,r_u)$ the following estimate holds
	\begin{equation}
		\| \mathcal{T}(\varepsilon_1,h_1)(x) - \mathcal{T}(\varepsilon_2,h_2)(x) \| \leq \beta^{\varepsilon_1} \|h_1 - h_2\| + K\left(1+L\frac{\mu^{\varepsilon_1}}{\xi^{\varepsilon_1}}\right)|\varepsilon_1-\varepsilon_2|.
	\end{equation}
\end{thm}
In order to get the estimate for the derivative of the graph transform first define analogously to the notation of Section \ref{graph_der},
$z(\varepsilon,h)(x) = (G(\varepsilon,h)(x),h(G(\varepsilon,h)(x)))$ and $$
F(\varepsilon,h,M)(x) = \left(\frac{\partial f^\varepsilon_x}{\partial x}(z(\varepsilon,h)(x)) + \frac{\partial f^\varepsilon_x}{\partial y}(z(\varepsilon,h)(x)) M(G(\varepsilon,h)(x))) \right) ^{-1}.
$$
The argument that follows the lines of the proof of Lemma \ref{lem:LipDGdx-estm} allows us to deduce the following result.
\begin{lem}
	\label{lem:LipDGdx-estm_e}
 Assume that, for $i\in \{1,2\}$ we have $h_i
	\in H$ and $\|M_i\| \leq L$ and
	\begin{equation}
		\|M_i(x_1) - M_i(x_2)\| \leq L_M \|x_1-x_2\|\, \ \textrm{for every}\ \ x_1,x_2\in \overline{B}_u(0,r_u).
	\end{equation}
	Then
	\begin{equation}
		\left\|F(\varepsilon_1,h_1,M_1)(x)- F(\varepsilon_2,h_2,M_2)(x)\right\| \leq  C_1^{\varepsilon_1} \|h_1-h_2\| + \frac{1}{\xi^{\varepsilon_1}\xi^{\varepsilon_2}}\left\|\frac{\partial f^{\varepsilon_1}_x}{\partial y}\right\|\|M_1-M_2\| + C_2^{\varepsilon_1}|\varepsilon_1-\varepsilon_2|.
	\end{equation}
	where  $C_1^{\varepsilon_1}=C(\varepsilon_1,N,f^{\varepsilon_1},Df^{\varepsilon_1},D^2f^{\varepsilon_1},L, L_M)$ and $C_2^{\varepsilon_1}=C(\varepsilon_1,N,f^{\varepsilon_1},Df^{\varepsilon_1},L, L_M,K)$.
\end{lem}
The proof of the next result uses Lemma \ref{lem:LipDGdx-estm_e} and  follows the lines of the proof of Lemma \ref{lem:LipT2-estm}.
\begin{thm}
	\label{lem:LipT2-estm_e}
Assume that, for $i\in \{1,2\}$ we have $h_i
	\in H$ and $\|M_i\| \leq L$ and
	\begin{equation}
		\|M_i(x_1) - M_i(x_2)\| \leq L_M \|x_1-x_2\|\, \ \textrm{for every}\ \ x_1,x_2\in \overline{B}_u(0,r_u).
	\end{equation}
	Then
	\begin{equation}
		\left\|\mathcal{U}(\varepsilon_1,h_1,M_1)(x)- \mathcal{U}(\varepsilon_2,h_2,M_2)(x)\right\| \leq   C_1^{\varepsilon_1}\|h_1-h_2\|+\frac{\beta^{\varepsilon_1}}{\xi^{\varepsilon_2}}\|M_1-M_2\|+C_2|\varepsilon_1-\varepsilon_2|.
	\end{equation}
	where $C_1^{\varepsilon_1}=C(\varepsilon_1,N,f^{\varepsilon_1},Df^{\varepsilon_1},D^2f^{\varepsilon_1},L,L_M)$ and $C_2^{\varepsilon_1}=C(\varepsilon_1,N,f^{\varepsilon_1},Df^{\varepsilon_1},L,L_M,K)$.
	\end{thm}
Theorems \ref{thm:T-contr_e} and \ref{lem:LipT2-estm_e} imply the desired $C^1$ continuity of the graph transform map \\
$(\varepsilon,h,M)\mapsto (\mathcal{T}(\varepsilon,h),\mathcal{U}(\varepsilon,h,M))
$ for the unstable manifold.

\subsection{Graph transform for the stable manifold.}
The graph transform with parameter for the stable manifold is defined in the following way:  given $y \in \overline{B}_s(0,r_s)$ we look for $x=\mathcal{S}(\varepsilon,v)(y)$ such that point $f^{\varepsilon}(\mathcal{S}(\varepsilon,v)(y),y)$ belongs to image of  $v$, i.e. there exists $y_0$ such that
\begin{equation}
	f^\varepsilon(\mathcal{S}(\varepsilon,v)(y),y)=(v(y_0),y_0).
\end{equation}
The next result is proved anaogously to Theorem \ref{thm:S-contrakcja}.
\begin{thm}
	\label{thm:S-contrakcja_e}
Let $y\in \overline{B}_s(0,r_s)$. For $v_1,v_2 \in V$ we have
	\begin{equation}
		\|\mathcal{S}(\varepsilon_1,v_1)(y) - \mathcal{S}(\varepsilon_2,v_2)(y_2)\| \leq \frac{\|v_1- v_2\|}{\xi_1^{\varepsilon_1}} + \frac{1}{\xi_1^{\varepsilon_1}}K\left(1+\frac{1}{L}\right)|\varepsilon_1-\varepsilon_2|.
	\end{equation}
\end{thm}
Now, the graph transform for the derivative of the stable manifold is given by the formula
\begin{align}
	&   \mathcal{R}(\varepsilon,v,M)(y)= \left(\frac{\partial f^{\varepsilon}_x}{\partial x}(z(\varepsilon,v)(y)) -  M(f^\varepsilon_y(z(\varepsilon,v)(y))) \frac{\partial f^{\varepsilon}_y}{\partial x}(z(\varepsilon,v)(y)) \right)^{-1} \cdot  \label{eq:S2-def_e} \\
	& \qquad \qquad \qquad \qquad   \left( M(f^\varepsilon_y(z(\varepsilon,v)(y)))\frac{\partial f^\varepsilon_y}{\partial y}(z(\varepsilon,v)(y)) -   \frac{\partial f^\varepsilon_x}{\partial y}(z(\varepsilon,v)(y)) \right), \nonumber
\end{align}
with $z(\varepsilon,v)(y)=(\mathcal{S}(\varepsilon,v)(y),y)$. The following result is proved analogously to Theorem \ref{lem:S-M-Lip_2}, taking into account the additional terms that come from the difference between $f^{\varepsilon_1}$ and $f^{\varepsilon_2}$ and their derivatives.
\begin{thm}
	\label{lem:S-M-Lip_2_e}
	Assume that, for $i\in \{1,2\}$ we have $v_i \in V$,
	\begin{equation}
		\|M_i\|\leq \frac{1}{L},
	\end{equation} and
	\begin{equation}
		\|M_i(y_1) - M_i(y_2)\| \leq L_M \|y_1 - y_2\|\ \ \textrm{for every}\ \ y_1,y_2\in \overline{B}_s(0,r_s).
	\end{equation}
	Then
	\begin{equation}
		\| \mathcal{R}(\varepsilon_2,v_2,M_2)(y) - \mathcal{R}(\varepsilon_1,v_1,M_1)(y)\|  \leq C^{\varepsilon_1}_1\|v_1-v_2\| + \frac{\mu^{\varepsilon_1}}{\xi^{\varepsilon_2}_1} \|M_1-M_2\| + C_2^{\varepsilon_1}|\varepsilon_1-\varepsilon_2|,
	\end{equation}
	where  $C^{\varepsilon_1}_1=C(\varepsilon_1,N,f^{\varepsilon_1},Df^{\varepsilon_1},D^2f^{\varepsilon_1},L,L_M)$ and $C^{\varepsilon_1}_2=C(\varepsilon_1,N,f^{\varepsilon_1},Df^{\varepsilon_1},L,L_M,K)$.
\end{thm}
Theorems \ref{thm:S-contrakcja_e} and \ref{lem:S-M-Lip_2_e} imply the desired $C^1$ continuity of the graph transform map \\
$(\varepsilon,h,M)\mapsto (\mathcal{S}(\varepsilon,h),\mathcal{R}(\varepsilon,h,M)),
$ for the stable manifold.

\section*{Statements}
On behalf of all authors, the corresponding author states that there is no conflict of interest. We do not analyse or generate any datasets, because our work proceeds within a theoretical and mathematical approach.

This work was supported by National Science Center (NCN) of Poland under project No. UMO2016/22/A/ST1/00077. Work of PK was also partially supported by Ministerio de Ciencia e Innovaci\'{o}n of Kingdom of Spain under project No. PID2024-156228NB-I00 and by	FAPESP, Brazil grant 2020/1407.

\end{document}